%% file: exp.tex
\documentclass{amsart}

\usepackage{amssymb}
\usepackage{tikz}
\usepackage{verbatim}
\usepackage{hyperref}
\usepackage[english]{babel}
\usepackage{wasysym}
\usepackage{todonotes}
\usepackage{xcolor}

\input macros.tex

\author{Manuel Bodirsky and Bertalan Bodor}
\address{Institute of Algebra, TU Dresden}
  \thanks{The authors have received funding from the European Research Council under the European Community's Seventh Framework Programme (FP7/2007-2013 Grant Agreement no. 257039, CSP-Infinity).}

\begin{document}

\title{Structures with Small Orbit Growth}
\maketitle

\begin{abstract}
Let $\kexpp$ be the class 
of all structures $\fa$ such that the automorphism group of $\fa$ 
has
at most $c n^{d n}$ orbits
in its componentwise action on the set of $n$-tuples with pairwise distinct entries, 
 for some constants $c,d$ with $d < 1$.
We show that $\kexpp$ is precisely the class of finite covers of first-order reducts of unary structures, and also that $\kexpp$ is precisely the class of first-order reducts of finite covers of unary structures. 
It follows that the class of first-order reducts of finite covers of unary structures is closed under taking model companions and model-complete cores, which is
an important property when studying the constraint satisfaction problem for structures from
$\kexpp$. We also show that Thomas' conjecture holds
for $\kexpp$: all structures in
$\kexpp$ have finitely many first-order reducts up to first-order interdefinability. 
\end{abstract}

\section{Introduction}
A \emph{first-order reduct} of a structure $\fa$ is a relational structure with the same domain as $\fa$ whose relations are first-order definable 
over $\fa$.
Simon Thomas conjectured that every 
homogeneous structure $\fa$ with finite relational signature has only finitely many first-order reducts up to first-order interdefinability~\cite{RandomReducts}. The conjecture 
has been verified for many famous homogeneous structures $\fa$: e.g.~for 
the ordered rationals~\cite{Cameron5}, the countably infinite random graph~\cite{RandomReducts}, the homogeneous universal $K_n$-free graphs~\cite{Thomas96}, 
the expansion of $({\mathbb Q};<)$ by a constant~\cite{JunkerZiegler}, 
the universal homogeneous partial order~\cite{Poset-Reducts},
and the random ordered graph~\cite{42},
and many more~\cite{agarwal,AgarwalKompatscher,BodJonsPham,BBPP18}. 
If we drop the assumption that the signature
of the homogeneous structure $\fa$ is relational, 
then the conjecture of Thomas is false even if we keep the assumption that $\fa$ is $\omega$-categorical: already the countable atomless Boolean algebra has infinitely many first-order reducts~\cite{BodorCameronSzabo}. 

Thomas' conjecture highlights our limited understanding of 
the class of homogeneous structures $\fa$ with finite relational signature.
One approach to widen our understanding is to study homogeneous
structures for some fixed signature;
for example, classifications exist for the class 
of all homogeneous tournaments~\cite{Lachlan}, homogeneous undirected graphs~\cite{Henson}, homogeneous partial orders~\cite{Schmerl}, general homogeneous digraphs~\cite{Cherlin},
homogeneous permutations~\cite{CameronPermutations}, 
and homogeneous coloured multipartite graphs
\cite{MultipartiteJenkinsonTrussSeidel,MultipartiteLockettTruss}. However, 
already the class of homogeneous 3-uniform hypergraphs appears to be very difficult~\cite{AkhtarLachlan95}. 
If we impose additional
assumptions, e.g., that the age of $\fa$ can be described by finitely many forbidding substructures,
we might hope for systematic understanding and effectiveness results for various questions. 
However, it is not clear how to use this assumption for proving that
$\fa$ has finitely many first-order reducts. 

Another approach to understanding the class of homogeneous structures, followed in this paper,
is to 
start with the most symmetric structures in this class. 
Symmetry can be measured by the number of orbits $o^i_n(G)$ of the diagonal action of the automorphism group $G = \aut(\fa)$ 
on tuples from $A^n$ that have pairwise distinct entries. By the theorem of Engeler, Ryll-Nardzeski, and Svenonius, these orbits are in one-to-one correspondence with the model-theoretic types of $n$ pairwise distinct elements in $\fa$. 
Alternatively, we might count the number of orbits $o^s_n(G)$ of the action of $G$ on $n$-element subsets of $A$. The investigation of both of these measures has been pioneered by Cameron; see~\cite{Oligo} for an introduction to the subject.
The sequence $o^i_n(G)$ is linked to \emph{labeled} enumeration problems, which are the most intensively studied counting problems in enumerative combinatorics, while $o^s_n(G)$ is linked to \emph{unlabeled} enumeration problems.
Many structural results about $G$ are available when we impose restrictions on $o^s_n(G)$; see, e.g.,~\cite{Macpherson-Orbits,Macpherson-GrowthRates,Macpherson-RapidGrowth}.
The present article, in contrast, focuses on restricting $o^i_n(G)$. 

A structure $\fa$ is finite if and only if
$o^i_n(\aut(\fa))$ is eventually 0. 
It is a by-product of our results that the class
$\kexp$ of all structures $\fa$ where $o^i_n(\aut(\fa))$ grows at most exponentially  equals the class of first-order reducts of unary structures; by a \emph{unary structure} we mean any at most 
countable structure with finitely many unary relations. 
Our main result pushes this further:
we study the class $\kexpp$ of structures $\fa$ such that $o^i_n(\aut(\fa))$ is bounded by $cn^{dn}$ for some constants $c,d$ with $d< 1$.
Note that for example the structure $({\mathbb Q};<)$ does not belong to $\kexpp$ because $o^i_n({\mathbb Q};<) = n!$. Also, $\kexpp$ contains no structure $\fa$ with a definable equivalence relation with infinitely many infinite classes because $o^i_n(\fa)$ would in this case be
at least as large as the $n$-th Bell number, which
grows asymptotically faster than $cn^{dn}$ (see Lemma~\ref{counting}).
We show that $\kexpp$ contains precisely those structures that are \emph{finite covers} of first-order reducts of unary structures. 

Finite covers in model theory and infinite permutation groups have been studied in the context of classifying totally categorical structures~\cite{AhlbrandtZiegler-Cover,HodgesPillay,HrushovskiTotallyCategorical} and, more generally, for studying 
$\omega$-categorical $\omega$-stable structures
~\cite{CherlinHarringtonLachlan,CherlinHrushovski}.
Finite covers became an important topic in its own~\cite{Evans,EvansPastori,Pastori}; we refer to the survey article for an  introduction~\cite{EvansIvanovMacpherson}.
It follows from our result that the class of finite covers of reducts of unary structures equals 
the class of first-order reducts of finite covers of unary structures. 
Using the terminology of~\cite{EvansIvanovMacpherson}, we show that all finite covers of unary structures \emph{split}, but not necessarily \emph{strongly}. All structures in $\kexpp$ can be expanded to structures that are homogeneous in a finite relational language, and we show that they all satisfy Thomas' conjecture. 
The proof uses a result of Macpherson which implies that structures in $\kexpp$ which have a \emph{primitive} automorphism group must be highly transitive~\cite{Macpherson-Orbits}.

The class $\kexpp$ can be seen as the \emph{`smallest reasonably robust class that contains all finite structures as well as some infinite ones'} (for formalisations of this statement, see Section~\ref{sect:summary}). So whenever
a statement that holds for all finite structures 
needs to be generalised to a class of \emph{`slightly infinite structures'}, it might be a good idea to try to
first prove the statement for $\kexpp$. This is precisely the situation for the constraint satisfaction problem.

\subsection{Complexity of constraint satisfaction}
Let $\fb$ be a structure with finite relational signature. The \emph{constraint satisfaction problem} for $\fb$ is the computational problem of deciding whether a given finite structure $\fa$ with the same signature as $\fb$ has a homomorphism 
to $\fb$. For finite structures $\fb$, Feder and Vardi~\cite{FederVardi} conjectured that 
the computational complexity of $\csp(\fb)$ satisfies a \emph{dichotomy}:  it is either in P or NP-complete.  Using concepts and techniques from universal algebra, Bulatov and Zhuk recently presented independent proofs of this conjecture~\cite{BulatovFVConjecture,ZhukFVConjecture}.

The universal-algebraic approach can also be applied when $\fb$ is countably infinite and $\omega$-categorical. In this case, the computational complexity of $\fb$ is captured by the \emph{polymorphism clone of $\fb$} (see~\cite{BodirskyNesetrilJLC}), which can be seen as a generalisation of the automorphism group of $\fb$: it consists of all homomorphisms from $\fb^n$ to $\fb$, for $n \in \set$. 
Moreover, every $\omega$-categorical structure 
$\fb$ is homomorphically equivalent to an (up to isomorphism unique) structure $\fc$ with the property that the automorphisms of $\fc$ lie dense in the endomorphisms of $\fc$, called the \emph{model complete core} of $\fb$. The model-complete core $\fc$ of $\fb$ is again $\omega$-categorical, and has the same CSP as $\fb$, so that we prefer to analyse $\fc$ rather than $\fb$. 
This simplification of the classification problem is a key step for many results (see, e.g.,~\cite{tcsps-journal},\cite{BodPin-Schaefer-both,BartoPinskerDichotomy}), including the finite-domain classification~\cite{BulatovFVConjecture,ZhukFVConjecture}. 

Therefore, if we want to classify the computational complexity of $\csp(\fb)$ for all structures $\fb$ from a class $\mathcal C$, it is important whether the class $\mathcal C$ is closed under the formation of model-complete cores. 
When $\fc$ is the model-complete core of $\fb$, 
then it is easy to see that 
$o^i_n(\fc) \leq o^i_n(\fb)$; hence, in particular
the classes $\kexp$ and $\kexpp$ are closed under taking model-complete cores. 
This makes these classes attractive goals for
extending the mentioned dichotomy result from finite domains. 

As mentioned before, our results imply that 
every structure in $\kexp$ is a first-order reduct of a unary structure. For those structures, 
it has already been shown that they are 
in P or NP-complete~\cite{BodMot-Unary} (using the mentioned dichotomy for finite-domain CSPs). 
Our main result states that $\kexpp$ is precisely the class of first-order reducts of finite covers of unary structures. For classifying the complexity of the CSP for all structures in this class, our result implies that we can assume without loss of generality that these structures are model-complete cores. We thus see our result
as a first step towards classifying the CSP for first-order reducts of finite covers of unary structures.

 
\subsection{Definable sets with atoms} 
\label{sect:atoms}
In theoretical computer science one is interested
in finite representations of infinite structures;
one approach to this is the framework of 
\emph{definable sets} and \emph{computation with atoms}~\cite{DBLP:journals/corr/BojanczykKL14,DBLP:conf/lics/BojanczykKLT13}. 
This leads to new models of computation over infinite structures with interesting links to long-standing open problems in finite model theory, namely the question whether there is a logic for P and computation in choiceless polynomial time~\cite{BojanczykTorunczyk18}. 

If the `atom structure' is $({\mathbb N};=)$ 
(which is besides $({\mathbb Q};<)$ the most frequently used base structure in this area) then definable sets (in this case also studied under the name \emph{nominal sets}~\cite{GabbayPitts}) correspond precisely to the class ${\mathcal K}_=$ of structures that are first-order interpretable over $({\mathbb N};=)$ in the sense of model theory (for an explicit discussion of the connection, see~\cite{definable-homomorphisms}, Lemma 7 and the remarks thereafter). The class 
${\mathcal K}_=$ 
might appear to be trivial to many model theorists 
(all structures in it are $\omega$-categorical, $\omega$-stable, and they are first-order reducts of homogeneous finitely bounded structures), 
but in fact many questions about this class remain open; see Section~\ref{sect:open} for a small sample of open problems. 
It follows from our results (see Remark~\ref{rem:triv-cover-interpretation}) 
that  $\kexpp \subseteq {\mathcal K}_=$ 
and we can answer for
$\kexpp$ many questions that we cannot answer
for the class ${\mathcal K}_=$ in general.
So our results can also be seen as a first step towards a better understanding of ${\mathcal K}_=$. 
 
\input prelims.tex

\input RU.tex
\input FU.tex
\input FRU.tex
\input Kexp.tex
\input RU2.tex

\input additional.tex

\section{Conclusion and Open Problems}
\label{sect:open}
Our results imply that all structures in
$\kexpp$ are $\omega$-stable (Remark~\ref{rem:triv-cover-interpretation}), that they are first-order reducts of finitely bounded homogeneous structures (Proposition~\ref{lem:finitely-bounded}), and that they satisfy Thomas' conjecture (Corollary~\ref{thomas_strong}). Do $\omega$-stable homogeneous structures with finite relational signature in general satisfy Thomas' conjecture, i.e., do they have finitely many reducts up to interdefinability? 
Note that if we drop the assumption about having a relational signature, then the answer is `no' 
even if we insist on $\fa$ being still $\omega$-categorical and $\omega$-stable (this follows from the example given in ~\cite{BodorCameronSzabo}, which is the expansion of the countably infinite dimensional vector space over the two-element field with one non-zero constant).

Answering the question of the previous paragraph might be very ambitious, so we propose to first study a more concrete and fundamental class of structures. 
Let ${\mathcal K}_=$ be the class
of all structures with
a first-order interpretation over $({\mathbb N};=)$
(which we have already discussed in Section~\ref{sect:atoms}).  Do the structures in
${\mathcal K}_=$ satisfy Thomas' conjecture? 
Is the model companion of a structure in ${\mathcal K}_=$ also in ${\mathcal K}_=$? 
We ask the same question for the model-complete cores of structures in ${\mathcal K}_=$. 


By our results, structures from $\kexpp$ can be represented on a computer as follows. 
First, every trivial covering 
$\fa$ of a structure $\fb \in \unary$ is interdefinable with a homogeneous structure in a finite relational signature $\fc$ (Proposition~\ref{ramsey})
and is finitely bounded (Lemma~\ref{lem:finitely-bounded}). 
So we can represent $\fc$ up to isomorphism by specifying these bounds. 
Second, finite-signature first-order reducts of $\fc$ can be represented by listing formulas for the relations of the reduct (we can assume that these formulas are quantifier-free since $\fc$ is  homogeneous in a finite relational signature and hence has quantifier elimination), and storing these together with the representation for $\fc$. 
We now ask which of the following problems are algorithmically decidable: 
\begin{enumerate}
\item given two structures in $(R \circ F)(\unary)$, 
decide whether they are isomorphic. 
\item  given two structures in $(R \circ F)(\unary)$, 
decide whether they are interdefinable. 
\item given two structures in $(R \circ F)(\unary)$, 
decide whether they are bi-interpretable. 
\end{enumerate}
Szymon Toru\'{n}czyk (personal communication) observed that the first of these questions (about deciding isomorphism of two given structures) is in the larger setting of reducts of finitely bounded homogeneous structures equivalent to an open problem about decidability of first-order definability from~\cite{BPT-decidability-of-definability} (the final open problem mentioned there).

\bibliographystyle{alpha}
\bibliography{local.bib}

\end{document}

%% file: macros.tex
\def\acts{\curvearrowright}

\newtheorem{theorem}{Theorem}[section]
\newtheorem{lemma}[theorem]{Lemma}
\newtheorem{corollary}[theorem]{Corollary}
\newtheorem{proposition}[theorem]{Proposition}

\theoremstyle{definition}
\newtheorem{definition}[theorem]{Definition}

\newtheorem{example}[theorem]{Example}

\theoremstyle{remark}
\newtheorem{remark}[theorem]{Remark}
\theoremstyle{remark}

\DeclareMathOperator{\fin}{fin}
\DeclareMathOperator{\Eq}{Eq}
\newcommand{\sym}{\operatorname{Sym}}
\newcommand{\csp}{\operatorname{CSP}}
\newcommand{\alt}{\operatorname{Alt}}
\newcommand{\aut}{\operatorname{Aut}}
\newcommand{\inv}{\operatorname{Inv}}

\newcommand{\id}{\operatorname{id}}

\newcommand{\fa}{\mathfrak{A}}
\newcommand{\fb}{\mathfrak{B}}
\newcommand{\fc}{\mathfrak{C}}
\newcommand{\fd}{\mathfrak{D}}

\newcommand{\cc}{\mathfrak{A}}
\newcommand{\dd}{\mathfrak{C}}
\newcommand{\ww}{\mathfrak{B}}

\newcommand{\orb}{o}
\newcommand{\oi}{o^i}
\newcommand{\os}{o^s}

\newcommand{\acl}{\operatorname{acl}}

\DeclareMathOperator{\Exp}{exp}

\newcommand{\gexp}{\mathcal{G}_{\Exp}}
\newcommand{\gexpp}{\mathcal{G}_{{\Exp}{+}}}

\newcommand{\kexp}{\mathcal{K}_{\Exp}}
\newcommand{\kexpp}{\mathcal{K}_{{\Exp}{+}}}
\newcommand{\set}{\mathbb{N}}
\newcommand{\sets}{\mathcal{S}}
\newcommand{\unary}{\mathcal{U}}

\newcommand{\struct}{\mathcal{S}}
\newcommand{\group}{\mathcal{G}}
\DeclareMathOperator{\nf}{nf}

\newcommand{\reduct}{R}
\newcommand{\reductfin}{R^{<\infty}}
\newcommand{\fcover}{F}
\newcommand{\mcore}{M}
\newcommand{\const}{C}
\numberwithin{equation}{section}

%% file: prelims.tex

\section{Preliminaries}
If $\sim$ is an equivalence relation on $X$ and $x \in X$, then $[x]_\sim$ denotes the equivalence class of $x$ with respect to $\sim$,
	and 
	$X/{\sim} := \{[x]_{\sim} \mid x \in X\}$ 
	denotes the set of all 
	$\sim$-classes.
	We write $|{\sim}|$ for $|X/{\sim}|$.
	If $\sim_1$ and $\sim_2$ are equivalence relations on $X$ then we say that $\sim_1$ is \emph{finer} than $\sim_2$ (or $\sim_2$ is \emph{coarser} then $\sim_1$) if $\sim_1$ is contained in $\sim_2$ (as binary relations).
	
\subsection{Permutation group notation}
When $G$ is a group we write $H \leq G$ if $H$ is a subgroup of $G$, and $H \triangleleft G$ if $H$ is a normal subgroup of $G$. We write $[G:H]$ for the index of $H$ in $G$. For any set $X$ we write $\sym(X)$ for the group of all permutations of $X$. If $G \leq \sym(X)$ and $x\in X$ then $G_x$ denotes the \emph{stabiliser} of the element $x$. 
	Let $Y\subseteq X$. Then 
	\begin{itemize} 
	\item $G_Y$ denotes the \emph{pointwise stabiliser}, and 
	\item $G_{\{Y\}}$ denotes the \emph{setwise stabiliser} of the set $Y$.
	\item $G|_{Y}$ denotes the restriction of $G$ to $Y$ provided that $Y$ is preserved by $G$. 
	\end{itemize}
	If $Y$ is finite, say $Y=\{x_1,\dots,x_n\}$, then we also use the notation $G_{x_1,\dots,x_n}$ for the pointwise stabiliser of the set $Y$. 
	

	Let $G$ be a permutation group on $X$. 
	An \emph{orbit} of $G$ is a set of the form 
	$\{g(x) \mid g \in G\}$ for some $x \in X$. 
	Then \emph{algebraic closure of $Y \subseteq X$ with respect to $G$} is
	 the union of the finite orbits of $G_Y$, and it is denoted by $\acl_G(Y)$. If $x\in X$, then we use the notation $\acl_G(x)$ instead of $\acl_G(\{x\})$. 
	 It is well-known that $\acl_G$ is a closure operator on the subsets of $X$, and in particular we have $\acl_G(\acl_G(Y)) = \acl_G(Y)$ for all $Y \subseteq X$. 
	If the group $G$ is clear from the context, then we will omit the subscript from this notation.

	An equivalence relation $\sim$ of $X$ is called a \emph{congruence} of a 
	permutation group $G \subseteq \sym(X)$ 
	if $x\sim y$ and $g\in G$ implies $g(x)\sim g(y)$ for all $x,y\in X$ and $g\in G$. In other words, an equivalence relation $\sim$ is a congruence if the corresponding partition is $G$-invariant.
	Every permutation group $G\subseteq \sym(X)$ has two \emph{trivial congruences}, namely $X^2$ and $\{(x,x): x\in X\}$. We call the former the \emph{universal congruence} and the latter the \emph{identity congruence}. 
	If $\sim$ is a congruence of some permutation group $G \subseteq \sym(X)$ then 
	$G$ acts naturally on $X/{\sim}$. The image of this action, as a subgroup of $\sym(X/{\sim})$, is denoted by $G/{\sim}$.

\begin{definition}\label{def:mu-pi}
Let $\pi \colon A \to B$ be a map. We write
$\sim_\pi$ for the equivalence relation $\{(a_1,a_2) \mid \pi(a_1) = \pi(a_2)\}$ on $A$. 
If $G$ is a permutation group on $A$ such that
$\sim_\pi$ is a congruence of $G$, then
$\pi$ gives rise to a homomorphism $\mu_\pi \colon G \to \sym(B)$ defined by
$\mu_\pi(g)(a) := \pi(g(\pi^{-1}(a)))$ 
(this is well-defined since $G$ preserves $\sim_\pi$). 
\end{definition}


\subsection{Direct products}
\label{sect:prod}
Let $I$ be a set. For each $i \in I$, let $A_i$ be a group. Then $\prod_{i \in I} A_i$ denotes the direct product of the $A_i$; i.e., the elements 
have the form $(a_i)_{i \in I}$ for $a_i \in A_i$, 
and group composition is defined point-wise. 
When the $A_i$ are permutation groups 
on disjoint sets $X_i$ for every $i \in I$, 
then $A := \prod_{i \in I} A_i$ acts naturally (intransitively) on $X :=\bigsqcup_{i\in I}X_i$ as follows: for $\alpha \in A$ and $x \in X$, define 
$\alpha(x) := \alpha_i(x)$ if $x \in X_i$.
It is easy to see that if each of the $A_i$ is closed in $\sym(X_i)$, then the permutation group defined by the action of $A$ on $X$ is closed in
$\sym(X)$, and hence is the automorphism
group of some relational structure with domain $X$. 

\subsection{Wreath products}
\label{sect:wr}
	Let $A$ be a group acting on the set $F$,
	and let $Y$ be a set.
	Let $H$ be a group acting on $Y$ and let $X:=F\times Y$. Then there are natural actions of the groups $N := \prod_{y\in Y}A$ and $H$
	on the set $X$, defined as follows.
\begin{enumerate}
\item If $\alpha\in N$ and $(f,y)\in X$, then $\alpha(f,y) := (\alpha_y(f),y)$,
\item If $\beta\in H$ and $(f,y)\in X$, then 
$\beta(f,y) := (f,\beta(y))$.
\end{enumerate}
	Let $G$ be the subgroup of $\sym(X)$ generated by the actions of $N$ and of $H$ on $X$;
	we view $N$ and $H$ as subsets of $G$.
	If $\alpha\in N$ and $\beta\in H$, then $$\beta^{-1}\alpha\beta(f,y)=\beta^{-1}\alpha(f,\beta(y))=\beta^{-1}(\alpha_{\beta(y)}(f),\beta(y))=(\alpha_{\beta(y)}(f),y)$$
	so $\beta^{-1} \alpha \beta \in N$ and 
	$N \triangleleft G$. 
	Then $G = NH$ and $N \cap H = \{\id_X\}$.
	Hence, the group $G$ can be written as the semidirect product $\prod_{y\in Y}A \rtimes H$. The group $G$ is called the \emph{wreath product} of the groups $A$ and $H$ (with its canonical imprimitive action on $X$) and will be denoted by $A\wr H$. 

\subsection{Interdefinability, bi-definability, bi-interpretability}
We write $A$, $B$, $C$ for the domains of the structures
$\fa$, $\fb$, $\fc$, respectively. 
If $G$ is a set of permutations on a set $A$
then $\inv(G)$ denotes the relational structure 
$\fa$ with domain $A$ which carries all relations that are preserved by all permutations of $G$. 
The operations $\aut$ and $\inv$ form a Galois
connection between the set of all relational structures $\fa$ with domain $A$ and the set of sets of permutations $G$ on $A$ (see, e.g., ~\cite{Bodirsky-HDR}). The permutation group $\aut(\inv(G))$ is the 
smallest permutation group containing $G$ that is \emph{closed} 
in $\sym(A)$ equipped with the topology of pointwise convergence. This topology is the restriction of the product topology on $A^A$ where $A$ is taken to be discrete. A permutation group
$G$ on $A$ is closed in $\sym(A)$ if and only if 
$G$ is the automorphism group of a relational structure. 
If $\fa$ is $\omega$-categorical, then the structure 
$\inv(\aut(\fa))$ is the expansion of $\fa$ 
by all relations that can be defined by a first-order formula in $\fa$ (this is a consequence of the proof of the theorem of Ryll-Nardzewski; see~\cite{HodgesLong}).

It follows that $\aut(\fa) \subseteq \aut(\fa')$ if and only if all relations of $\fa'$ are first-order definable (without parameters) over $\fa$; in this case we say that $\fa'$ is a \emph{first-order reduct} of $\fa$. Two structures on the same domain are called \emph{interdefinable} if they are reducts of one another. By the above, if $\fa$ or $\fa'$ is $\omega$-categorical, then $\fa$ and $\fa'$ are interdefinable if and only if $\aut(\fa) = \aut(\fa')$.

Two structures $\fa$ and $\fb$, not necessarily with the same domain, are called \emph{bi-definable} if 
there exists a bijection $f \colon A \to B$ between the domains of $\fa$ and $\fb$ such that 
$\fa$ and $\fb$ are interdefinable after identifying $A$ and $B$ along $f$. 
It follows that two $\omega$-categorical structures $\fa$ and $\fb$ are bi-definable if and only if $\aut(\fa)$ and $\aut(\fa')$ 
are isomorphic as permutation groups. 
For example, the structures 
$({\mathbb Z}; \{0\})$ and 
$({\mathbb Z}; \{1\})$ are bi-definable, but 
not interdefinable.


A \emph{($d$-dimensional) interpretation} of $\fb$ in $\fa$ is
a partial surjective map $I$ from $A^d$ to $B$
such that the pre-image of $B$, of the equality relation on $B$, and of each relation of $\fb$ under $I$ is first-order definable in $\fa$. If
$\fa$ has a $d$-dimensional first-order interpretation $I$ in $\fb$
and $\fb$ has an $e$-dimensional first-order interpretation $J$ in $\fa$ such that the relation
$\{(x,y_{1,1},\dots,y_{d,e}) \mid x = J(I(y_{1,1},\dots,y_{d,1}),\dots,I(y_{1,e},\dots,y_{d,e}))\}$
is first-order definable in $\fb$ and
$\{(x,y_{1,1},\dots,y_{d,e}) \mid x = I(J(y_{1,1},\dots,y_{1,e}),\dots,J(y_{d,1},\dots,y_{d,e}))\}$
is first-order definable in $\fa$,
then $\fa$ and $\fb$ are called \emph{bi-interpretable}. 
By a result of Coquand, Ahlbrandt, and Ziegler~\cite{AhlbrandtZiegler}, two $\omega$-categorical structures 
$\fa$ and $\fb$ are bi-interpretable if and only if
$\aut(\fa)$ and $\aut(\fb)$ are \emph{topologically isomorphic}, i.e., isomorphic via a mapping which is a homeomorphism with respect to the pointwise convergence topology.

\subsection{Orbit growth and some classes of structures}
Let $X$ be a countably infinite set. 
There are three natural counting sequences
attached to a permutation group on $X$, introduced and discussed in general in~\cite{CameronCounting,Oligo}.
\begin{definition}
	Let $G \subseteq \sym(X)$ be a permutation group and let $n\in \set$. Then 
\begin{itemize}
\item $\orb_n(G)$ denotes the \emph{number of $n$-orbits of $G$}, i.e., the number of orbits of the natural action $G\acts X^n$,
\item $\oi_n(G)$ denotes the \emph{number of injective $n$-orbits of $G$}, i.e., the number of orbits of the natural action $G\acts X^{(n)}$ where $X^{(n)} := \big \{(x_1,\dots,x_n)\in X^n \mid x_i\neq x_j \text{ for all distinct $i,j \in \{1,\dots,n\}$}\big \}$,
\item $\os_n(G)$ denotes the \emph{number of orbits of $n$-subsets of $G$}, i.e., the number of orbits of the natural action $G\acts {X\choose n}(=\{Y\subset X^n: |Y|=n \}$).
\end{itemize}
	If $\fa$ is a structure then let $$\orb_n(\fa):=\orb_n(\aut(\fa)),\,\oi_n(\fa):=\oi_n(\aut(\fa)),\,\os_n(\fa):=\os_n(\aut(\fa)).$$
	In the notation above we omit the reference to the group $G$ or the structure $\fa$ if it is clear from the context.
\end{definition}
A permutation group is called \emph{transitive}
if $o_i(G) = 1$ and \emph{highly transitive}
if $o_i(G) = 1$ for all $i \in {\mathbb N}$.

\begin{definition}\label{def:oligo}
	A permutation group $G\subseteq \sym(X)$ is called \emph{oligomorphic} if $\orb_n(G)$ is finite for all $n$.
\end{definition}

Clearly, in Definition~\ref{def:oligo} we could have equivalently required that $\oi_n$ or $\os_n$ are finite for all $n$. By the theorem of Engeler, Ryll-Nardzewski, and Svenonius, a countably infinite relational structure $\fa$ is $\omega$-categorical if and only if $\aut(\fa)$ is oligomorphic (see for instance~\cite{HodgesLong}). 
	In this paper we are particularly interested in the following classes of structures and permutation groups.

\begin{definition}
\begin{itemize}
\item Let $\gexp$ denote the class of those permutation groups $G$ acting on a countable set $X$ for which there is a constant $c$ such that $\oi_n(G)\leq c^n$.
\item Let $\kexp$ denote the class of all countable structures $\fa$ with an automorphism group in $\gexp$. 
\item Let $\gexpp$ denote the class of those permutation groups $G$ acting on a countable set $X$ for which there are constants $c$ and $d<1$ such that $\oi_n(G)\leq cn^{dn}$.
\item Let $\kexpp$ denote the class of all countable structures $\fa$ with an automorphism group in $\gexpp$. 
\end{itemize}
\end{definition}

\begin{remark}
	Note that the conditions $G\in \gexp$ and $G\in \gexpp$ imply that $G$ is oligomorphic, and therefore $\fa\in \kexp$ and $\fa\in \kexpp$ imply that $\fa$ is $\omega$-categorical.
\end{remark}
	
	We write $\set$ not only for the set of natural numbers, but also for the structure with the empty signature whose domain is $\set$. 

\begin{definition}
We write 
\begin{itemize}
\item $\sets$ for the class of all at most countable structures that are first-order interdefinable with a structure having the empty signature;
\item $\unary$ for the class of at most countable structures
that are first-order interdefinable with a structure
having a finite signature of unary relation symbols; 
\item $\unary^*$ for the class of the structures $\fa\in \unary$ such that every orbit of $\aut(\fa)$ is either a singleton or infinite.  
\end{itemize}
\end{definition}

When $\mathcal C$ is a class of structures,
we write $\mathcal C_{\nf}$ for the class consisting of all the structures in $\mathcal C$ that have no finite orbits. Note that $\set \in \sets \subset \unary_{\nf} \subset \unary^* \subset \unary$ and that $(\unary^*)_{\nf} = \unary_{\nf}$.

\subsection{Congruences of oligomorphic groups}
We need the following easy observation about oligomorphic groups.

\begin{proposition}\label{fin_many_cong}
	Every oligomorphic permutation group has finitely many congruences.
\end{proposition}

\begin{proof}
	Every congruence of a permutation group is a union of its 2-orbits. Then the claim follows directly from oligomorphicity.
\end{proof}

\begin{lemma}\label{acl1}
	Let $G$ be an oligomorphic permutation group, and let $\sim$ be a congruence of $G$ which has finite equivalence classes. Then $a\sim b$ implies $b\in \acl_G(a)$.
\end{lemma}

\begin{proof}
	Suppose that $a\sim b$, but $b\not\in \acl(a)$. Then the orbit of $b$ in $G_{a}$ is infinite. Let $b'$ be any element in this orbit. Then by definition $a\sim b'$. Hence the equivalence class of $a$ is infinite, a contradiction.  
\end{proof}

	If $\sim_1$ and $\sim_2$ are congruences, then the inclusion-wise smallest congruence relation that contains both $\sim_1$ and $\sim_2$ is called the equivalence relation \emph{generated} by $\sim_1$ and $\sim_2$.

\begin{lemma}\label{acl3}
	Let $G$ be an oligomorphic permutation group, and let $\sim_1$ and $\sim_2$ be congruences of $G$ with finite classes. Then the
	congruence generated by $\sim_1$ and $\sim_2$ also has finite classes.
\end{lemma}

\begin{proof}
	Let $\sim$ be the congruence generated by $\sim_1$ and $\sim_2$, and suppose that $a\sim b$. Then there exists a sequence $a_0,b_0,a_1,b_1,\dots, a_k,b_k$ with $a_0 = a$ and $b_k = b$ such that $a_i\sim_1 b_i$ for all $i \leq k$ and $b_i \sim_2 a_{i+1}$ for all $i < k$. By Lemma \ref{acl1}, this implies that $b_i\in \acl(a_i)$ and $a_{i+1}\in \acl(b_i)$ for all $i$. Since $\acl$ is a closure operator it follows that $b\in \acl(a)$. Since $G$ is oligomorphic it follows that $\acl(a)$ is finite. Therefore, the equivalence class of $a$ is also finite.
\end{proof} 

\begin{definition}
	Let $G$ be an oligomorphic permutation group. Then
\begin{itemize}
\item $\nabla(G)$ denotes the intersection of all congruences of $G$ with finitely many classes,
\item $\Delta(G)$ denotes the smallest congruence that contains all  congruences of $G$ with finite classes.
\end{itemize}
	If $\fa$ is an $\omega$-categorical structure, then we use the notation $\nabla(\fa):=\nabla(\aut(\fa))$, and $\Delta(\fa):=\Delta(\aut(\fa))$.
\end{definition}

\begin{remark}
	Since $G$ has finitely many congruences it follows that $\nabla(G)$ also has finitely many classes, i.e., it is the \emph{finest} congruence of $G$ with finitely many classes. By Lemma \ref{acl3} it follows that every class of $\Delta(G)$ is finite, i.e., $\Delta(G)$ is the \emph{coarsest} congruence of $G$ with finite classes.
\end{remark}

\begin{remark}
	If $x$ and $y$ are in the same orbit, then their $\Delta$-classes have the same size. If $G$ has finitely many orbits, it follows that there exists some $n \in {\mathbb N}$ such that all elements lie in a $\Delta$-class of size at most $n$. 
\end{remark}

The congruence $\Delta$ has the following
equivalent description. 

\begin{lemma}\label{delta_alt}
	Let $G$ be an oligomorphic permutation group on a countably infinite set $X$. Then $(x,y)\in \Delta(G)$ iff $y\in \acl_G(x)$ and $x\in \acl_G(y)$.
\end{lemma}

\begin{proof}
	Let $\Delta'(G)=\{(x,y) \mid y\in \acl(x) \wedge  x\in \acl(y)\}$. We claim that $\Delta'(G)$ is an equivalence relation. It is clear that $\Delta'(G)$ is reflexive and symmetric. The transitivity follows from the fact that $\acl$ is a closure operator. It is also clear from the definition that $\Delta'(G)$ is preserved by $G$. Hence, $\Delta'(G)$ is a congruence. For any $x\in X$ we have $[x]_{\Delta'(G)} \subseteq \acl_G(x)$, so every class of $\Delta'(G)$ is finite. Therefore, $\Delta'(G)$ is finer than $\Delta(G)$. On the other hand, if $(x,y)\in \Delta(G)$, then $y\in \acl(x)$ and $x\in \acl(y)$, and thus $(x,y)\in \Delta(G)$.
\end{proof}

	We often use the following observation throughout this text.

\begin{lemma}\label{sing_or_inf}
	Let $G$ be an oligomorphic permutation group. Then every class of $\nabla(G)$ is either infinite or a singleton.
\end{lemma}

\begin{proof}
	 If the class of $x \in X$ is finite, then its orbit is also finite. Indeed, let $O$ be the orbit of $x$. Then every class of $\nabla(G)$ in $O$ is of the same size. So if this size is finite, then $O$ is also finite since $\nabla$ has finitely many classes.

	Let $X_{\fin}$ be the union of the finite orbits of $G$. By oligomorphicity it follows that $X_{\fin}$ is finite. Then $\nabla':=\nabla(G)\cap \{(x,x) \mid x\in X_{\fin}\}$ is also a congruence of $\nabla(G)$. Since $X_{\fin}$ is finite the congruence $\nabla'$ has finitely many classes. This implies that $\nabla'=\nabla(G)$ and thus every class of $\nabla(G)$ within $X_{\fin}$ is a singleton.
\end{proof}

\begin{lemma}\label{delta_union_nabla}
	Let $G$ be an oligomorphic permutation group on $X$ 
	and let $\sim$ be a congruence of $G$ with finite classes. Then the congruence generated by $\sim$ and $\nabla(G)$ equals
	$\big\{(x,y) \in X^2 \mid ([x]_\sim, [y]_\sim) \in \nabla(G/{\sim})\big\}$.
	\end{lemma}
\begin{proof}
If $\pi \colon X\rightarrow X/{\sim}$ is the factor map $x \mapsto [x]_\sim$, and $\approx$ is a congruence of $G/{\sim}$, then $$\pi^{-1}(\approx) := \{(x,y)\in X^2 \mid (\pi(x),\pi(y))\in {\approx}\}$$ is a congruence of $G$ which is coarser than $\sim$. 
In fact, $\pi^{-1}$ defines a bijection between the congruences of $G/{\sim}$ and those congruences of $G$ which are coarser than $\sim$. 
The congruence $\pi^{-1}(\nabla(G/{\sim}))$ has finitely many classes since $\nabla(G/{\sim})$ has finitely many classes.  
Hence $\pi^{-1}(\nabla(G/{\sim})$ is the finest congruence of $G$ that is coarser than $\sim$
and has finitely many classes.
So, by definition, it equals the congruence generated by $\sim$ and $\nabla(G)$.
\end{proof}

\subsection{Finite covers}
We now introduce the concept of finite covers that plays a central role in this article. Forming finite covers may be viewed as a way to construct new $\omega$-categorical structures from known ones; 
a more appropriate way is to view them as a way to decompose $\omega$-categorical structures into (hopefully) simpler parts.

\begin{definition}\label{finite_def}
	Let $\cc$ and $\ww$ be structures. A mapping $\pi \colon \cc\rightarrow \ww$ is called a \emph{finite covering map} (or \emph{finite cover}) 
	 if
\begin{enumerate}
\item $\pi$ is surjective,
\item for each $w\in B$ the set $\pi^{-1}(w)$ is finite,
\item $\sim_\pi$ is preserved by $\aut(\cc)$, 
\item the image of $\aut(\fa)$ under $\mu_\pi$ equals $\aut(\fb)$.
\end{enumerate}
	(See Definition~\ref{def:mu-pi}) for the definition of $\sim_\pi$ and $\mu_\pi$.) 
	The sets $\pi^{-1}(w)$, for $w\in B$, are called the \emph{fibers} of the finite covering map 
	$\pi$.
	A structure $\cc$ is called a \emph{finite covering structure of $\ww$} if there is a finite covering map $\pi \colon \cc\rightarrow \ww$.
\end{definition}

\begin{remark}
	A finite covering structure of an $\omega$-categorical structure has an oligomorphic automorphism group, and hence is $\omega$-categorical. 
\end{remark}

\begin{remark}\label{rem:congr}
	Let $\cc$ be an arbitrary structure and 
	let $\sim$ be a congruence of $\aut(\cc)$. 
	If all $\sim$-classes are finite, then $\cc$ is a finite covering structure of the quotient structure $\cc/{\sim}$, where $\cc/{\sim}$ can be any structure such that
	$\aut(\cc/{\sim}) = \aut(\cc)/{\sim}$. 
In fact, 
	every finite covering structure is of this form.  
Indeed, let $\cc$ be a 
	structure. If $\pi \colon \cc\rightarrow \ww$ is a finite covering map, then 
	$\sim_{\pi}$ is a congruence of $\aut(\cc)$ 
	 and there is a natural bijection between $B$ and $A/{\sim}_{\pi}$ defined by $w \mapsto \pi^{-1}(w)$. 
	 Let us identify $B$ and of $A/{\sim}_{\pi}$ along this bijection, and 
	 	 let $\fa/{\sim}_{\pi}$ be any structure
	 such that $\aut(\fa/{\sim}_{\pi}) = \aut(\fb)$. 
	 The image of $\aut(\fa)$
	 under the homomorphism 
	 $\mu_\pi$ equals $\aut(\fb)$,
	 hence $\aut(\fa)/{\sim}_{\pi} = \aut(\fb) = \aut(\fa/{\sim}_\pi)$.
\end{remark}

We present a series of simple examples of finite covers; they illustrate different phenomena of finite covers on which we will comment later, referring back to these examples. 

\begin{example}\label{expl:trivial}
Let $\vec P_1\cdot \omega$ be the directed graph which is an infinite union of directed edges. 
Then $\vec P_1\cdot \omega$ is a finite covering structure of $\set$, with $\pi$ being the projection to the second argument. Also note that 
$\aut(\vec P_1\cdot \omega)$ is topologically isomorphic to $\sym({\mathbb N})$, and that 
$\vec P_1\cdot \omega$ and $({\mathbb N};\neq)$ are bi-interpretable but not bi-definable.
\end{example}

\begin{example}\label{expl:principal}
Let $K_2\cdot \omega$ be the graph which is an infinite union of undirected edges. 
Then $K_2\cdot \omega$ is a finite covering structure of $\set$, with $\pi$ being the projection to the second argument. Identifying the domain of $K_2 \cdot \omega$ with
$\{0,1\} \times {\mathbb N}$ so that 
$(u,n)$ is adjacent to $(v,m)$ iff $n=m$ and $u \neq v$, the automorphism group of $K_2\cdot \omega$ is the wreath product ${\mathbb Z}_2\wr \sym(\omega)$ (see Section~\ref{sect:wr}).
\end{example}

\begin{example}\label{expl:non-free}
Let $K_2\cdot \omega$ be the structure with domain $\{0,1\} \times \set$ from Example~\ref{expl:principal}, and let $\fa$ be the expansion of $K_2\cdot \omega$ by the equivalence relation $\Eq$
defined by $\Eq((u,n),(v,m))$ iff $u=v$. 
Then $\fa$ is a finite covering structure of $\set$
with respect to the covering map $\pi$ 
that maps $(u,n)$ to $n$. Note that $\aut(\fa)$
is isomorphic (as an abstract group) to
the direct product 
${\mathbb Z}_2 \times \sym(\set)$ (see Section~\ref{sect:prod} for direct products and other actions of direct products).
\end{example}

\begin{example}\label{expl:neither-free-nor-trivial}
Let $A := \{0,1,2,3\} \times \set$ and let 
$\pi \colon A \to \set$ be the projection to the second argument. Let $\fa$ be the graph with vertex set $A$ such that $(u,b)$ is adjacent to $(v,c)$ if and only if 
\begin{itemize}
\item $b=c$ and $u=v+1 \mod 4$, or
\item $b \neq c$ and $u=v \mod 2$. 
\end{itemize}
See Figure~\ref{fig:2fibers} for an illustration.
\begin{figure}
\begin{center}
\includegraphics[scale=0.6]{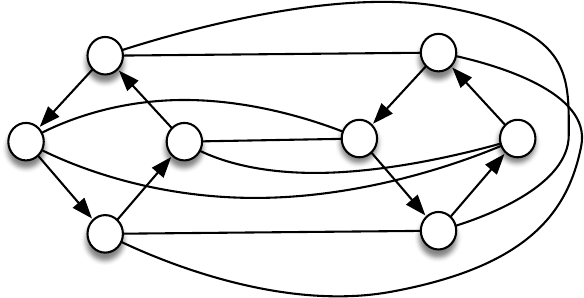}
\end{center}
\caption{An illustration of the subgraph of the structure $\fa$ from Example~\ref{expl:neither-free-nor-trivial} that is induced by 2 fibers.}
\label{fig:2fibers}
\end{figure}
Note that $\fa$ is a finite covering structure of $\set$ with respect to $\pi$. The automorphism group of $\fa$ equals 
$KH$ where 
\begin{itemize}
\item $H = \{\alpha \in \sym(A) \mid \text{if } \alpha(u,v) = (u',v') \text{ then } u=u'\}$ (i.e., $H$ is topologically isomorphic to $\sym(\set)$), and 
\item $K = \{ \alpha \in \prod_{i \in \set} {\mathbb Z}_4 \mid \text{for all } k,l \in \set: \alpha_k {\mathbb Z}_2 = \alpha_l {\mathbb Z}_2\}$
where ${\mathbb Z}_k$ is the cyclic group acting on $\{0,1,\dots,k-1\}$ and 
$\prod_{i \in \set} {\mathbb Z}_4$ is the direct product
in its intransitive action on $A$ (see Section~\ref{sect:prod}). 
\end{itemize}
\end{example}

\begin{example}\label{expl:twisted}
Let $\fb$ be the countable structure which carries an equivalence relation $\Eq$ with three classes $R,S,T$ such that $|S|=|T|$, and a unary relation 
symbol denoting the class $R$.  
Let $A := (\{0,1\} \times R) \cup (\{0\} \times (S \cup T))$. 
We define the structure $\fa$
with domain $A$ and the signature $\{E,F\}$ 
where $E$ and $F$ have arity two, and 
\begin{itemize}
\item $E((u_1,b_1),(u_2,b_2))$ 
holds if and only if $(u_1 = 0$, $b_1 \in R$, and $b_2 \in S$)
or ($u_1 = 1$, $b_1 \in R$, and $b_2 \in T$);  
\item $F((u_1,b_1),(u_2,b_2))$ 
holds if and only if $b_1 = b_2$. 
\end{itemize}
\begin{figure}
\begin{center}
\includegraphics[scale=0.4,angle=90]{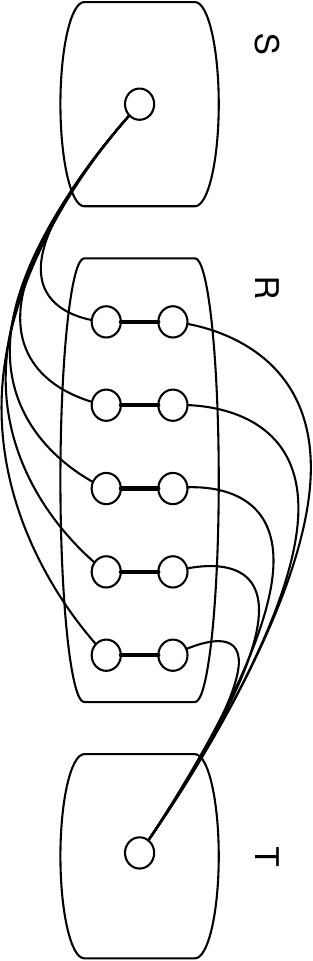}
\end{center}
\caption{An illustration of a subgraph of the structure $\fa$ from Example~\ref{expl:twisted} for the special case $|S| = |T| = 1$.}
\label{fig:twisted}
\end{figure}
See Figure~\ref{fig:twisted}.
Let $\pi \colon A \to B$ be the projection to the second argument. Then $\sim_{\pi} \, = F$ 
and $\pi$ is a finite covering. If $R,S,T$ are countably infinite then $\sim_{\pi} \, = F = \Delta(\fa)$. 
The automorphism group of $\fa$ is
isomorphic to a semidirect product 
$(\sym(R) \times \sym(S)^2) \rtimes {\mathbb Z}_2$.
\end{example}

\begin{definition}
Let $\pi \colon \fa \to \fb$ be a finite covering map, let $b \in B$, and let $S := \pi^{-1}(b)$. 
\begin{itemize}
\item The \emph{fiber group of $\pi$ at $b$} is the group 
$\aut(\fa)_{\{S\}}|_S$. 
\item The \emph{binding group of $\pi$ at $b$}
is the group $K|_S$ where $K$ is the kernel of $\mu_\pi$.
\end{itemize}
\end{definition}

 So the binding group at $b$ is a normal subgroup of the fiber group at $b$. 
If for some $b\in B$ the fiber group and the binding group at $b$ are unequal then $\pi$ is called \emph{twisted}. Example~\ref{expl:twisted} gives an example of a twisted finite cover; Examples~\ref{expl:trivial},~\ref{expl:principal},~\ref{expl:non-free}, and~\ref{expl:neither-free-nor-trivial} are not twisted.

\begin{remark}
The following terminology is not needed for stating or proving our results, but we mention it for a better understanding of the examples of finite covers that we have already presented. 
Let $\pi \colon \fa \to \fb$ be a finite covering map,
and let $B_b$ be the binding group at $b \in B$. 
Then $\pi$ 
 is called \emph{free} if the kernel of $\mu_\pi \colon \aut(\fa) \to \aut(\fb)$ equals 
 $\prod_{b \in B} B_b$. 
Example~\ref{expl:principal}, Example~\ref{expl:trivial},
Example~\ref{expl:neither-free-nor-trivial}, and Example~\ref{expl:twisted}
are free. 
 Example~\ref{expl:non-free} is an example
 of a finite cover which is not free: the binding group at each point is ${\mathbb Z}_2$ and equals the kernel of $\mu_\pi \colon \aut(\fa) \to \aut(\fb)$, which is therefore not equal to $\prod_{b \in B} B_b = {\mathbb Z}^\omega_2$. 
\end{remark}

\subsection{Trivial finite covers}
There are two important notions of triviality for finite covers, intended to describe those finite covers that have an automorphism group which is smallest possible. This is important for our purposes since we will describe general finite covering structures in our class by describing them as certain first-order reducts of trivial finite covers;
and, as we will see, trivial covers are much easier to describe.

\begin{definition}\label{trivial_trans}
	Let $\pi \colon \fa \rightarrow \fb$ be a finite covering map. We say that $\pi$ is  
	\begin{itemize}
	\item a \emph{trivial cover} if the kernel of $\mu_\pi \colon \aut(\fa) \to \aut(\fb)$ is trivial (only contains the identity permutation $\id_A$);
	\item 
	a \emph{strongly trivial cover} if all of its fiber groups are trivial.
\end{itemize}
	A structure $\fa$ is called a \emph{(strongly) trivial covering structure of $\fb$} if there is a finite covering map $\pi \colon \fa \rightarrow \fb$ which is (strongly) trivial.
\end{definition}

It is clear from the definition that $\pi$ is a trivial cover if and only if all of its binding groups are trivial. Hence, if $\pi$ is strongly trivial, then it is also trivial.
Example~\ref{expl:trivial} is an example of a strongly trivial finite covering.
Example~\ref{expl:twisted} is an example of a trivial finite covering which is not a strongly trivial finite covering. 
Examples~\ref{expl:principal},~\ref{expl:non-free}, and~\ref{expl:neither-free-nor-trivial} are examples of non-trivial finite coverings.

Next we give a sufficient condition for a structure $\ww$ under which every trivial cover of $\ww$ is strongly trivial.

\begin{lemma}\label{no_finite_index}
	Let $\fb$ be a structure such that for every $b\in B$ the stabiliser $\aut(\fb)_b$ has no nontrivial finite-index subgroups. Then every trivial cover of $\fb$ is strongly trivial.
\end{lemma}

\begin{proof}
	Let $\pi \colon \fa \rightarrow \fb$ be a trivial finite cover. Then $\mu_\pi$ is an isomorphism between $\aut(\fa)$ and $\aut(\fb)$. Let $b \in B$. We need to show that the fiber group of $\pi$ at $b$ is trivial. Put $S:=\pi^{-1}(b)$. Let us consider the mapping $\varphi \colon \aut(\fb)_b\rightarrow \sym(S)$ given by $h\mapsto \mu_\pi^{-1}(h)|_S$. 
	Then $\varphi$ is clearly a group homomorphism. Let $K$ be the kernel of this homomorphism. Then $K$ is a finite index subgroup of $\aut(\fb)_b$, and thus by our assumption $K=\aut(\fb)_b$. That is, $\varphi$ is the trivial homomorphism, which means that the fiber group of $\pi$ at $b$ is trivial.
\end{proof}

	We now give an explicit description of strongly trivial covers.

\begin{lemma}\label{lem:triv-cov}
	Let $\pi  \colon \cc \rightarrow \ww$ be a strongly trivial covering map. Then for each orbit $O$ of $\ww$ there exists a finite set $F_O$ and a mapping $\psi_O \colon  \pi^{-1}(O)\rightarrow F_O$ such that
\begin{itemize}
\item for every $w \in O$ the restriction of $\psi_O$
to $\pi^{-1}(w)$ is a bijection; 
\item $\psi_O(x) = \psi_O(\mu_\pi^{-1}(\beta)(x))$ for all $x\in \pi^{-1}(O)$ and $\beta\in \aut(\ww)$.
\end{itemize}
\end{lemma}

\begin{proof}
	Let us fix an element $b\in O$ and let $F_O:=\pi^{-1}(b)$. If $x\in \pi^{-1}(O)$ then there exists an automorphism $g$ of $\ww$ such that $g(\pi(x))=b$. Let us define $\psi_O(x)$ to be $\mu_\pi^{-1}(g)(x)$. We claim that $\psi_O(x)\in F_O$ and that its value is well-defined (i.e., it does not depend on our particular choice of $g$). The first claim is clear since by definition $$\pi(\mu_\pi^{-1}(g)(x))=\mu_\pi(\mu_\pi^{-1}(g))(\pi(x))=g(\pi(x))=b$$ 
thus $\mu_\pi^{-1}(g)(x)\in \pi^{-1}(b)=F_O$.
	
In order to show the second claim we need to show that if $h \in \aut(\fb)$ is such that $h(\pi(x))=g(\pi(x))=b$ then $\mu_\pi^{-1}(g)(x)=\mu_\pi^{-1}(h)(x)$. Since $(h^{-1}g)(\pi(x))=\pi(x)$ it follows that $(\mu_\pi^{-1}(h^{-1}g))|_{\pi(x)}$ is in the fiber group at $\pi(x)$. Since $\pi$ is strongly trivial this group is trivial, and hence $\mu_\pi^{-1}(h^{-1}g)(x)=x$. This implies that $$\mu_\pi^{-1}(h)(x)=(\mu_\pi^{-1}(h)\mu_\pi^{-1}(h^{-1}g))(x)=\mu_\pi^{-1}(g)(x).$$
Now the first item follows from the fact that if $w \in O$ is such that $g(w)=b$, then $g$ defines a bijection between $\pi^{-1}(w)$ and $\pi^{-1}(b)$. As for the second item let $x\in \pi^{-1}(O)$ and let $g\in \aut(\ww)$ be such that $g(\pi(x))=b$. If $\beta\in \aut(\ww)$, then 
$$(g\beta^{-1})(\pi(\mu_\pi^{-1}(\beta)(x)))=(g\beta^{-1})(\beta(\pi(x)))=g(\pi(x))=b,$$ and thus 
$$\psi_O(\mu_\pi^{-1}(\beta)(x))=(\mu_\pi^{-1}(g\beta^{-1}))(\mu_\pi^{-1}(\beta)(x))=\mu_\pi^{-1}(g)(x)=\psi_O(x).$$
\end{proof}

\begin{remark}\label{rem:triv-covers}
	Let the sets $F_O$ and the maps $\psi_O$ be defined as in Lemma~\ref{lem:triv-cov} for each orbit $O$ of $\aut(\fb)$. Then there is a natural bijection between $A$ and $\bigcup_O{(F_O\times O)}$ defined as $x \mapsto (\psi_O(x),\pi(x))$ where $O$ is the orbit of $\aut(\fb)$ containing $\pi(x)$. 
	If we identify each element of $A$ with its image under this bijection, then $\aut(\cc)$ consists of those permutations that fix the first coordinate of each element and that act as an automorphism of $\fb$ on the second coordinate. 
\end{remark}

\subsection{Covering reducts}
Let $\fa$ and $\fb$ be structures and let $\pi \colon \fa \rightarrow \fb$ be a finite covering map. A first-order reduct $\fc$ of $\fa$ is a \emph{covering reduct of $\fa$ with respect to $\pi$ (and $\fa$ is called a \emph{covering expansion of $\fc$ with respect to $\pi$}; see~\cite{EvansIvanovMacpherson})} if every $\alpha\in \aut(\dd)$ preserves $\sim_{\pi}$ and $\mu_\pi(\alpha) \in \aut(\ww)$. 

\begin{remark}
We do not need but mention that every finite cover $\pi \colon \fa \to \fb$ is an covering expansion of a free finite covering structure of $\fb$ with respect to $\pi$ (Lemma~2.1.3 in~\cite{EvansIvanovMacpherson}). 
\end{remark}

\begin{definition}
 Let $\pi \colon \fa \to \fb$ be a finite covering map. 
 \begin{itemize}
\item If $\fa$ is 
a covering reduct of a trivial covering structure of $\fb$ with respect to $\pi$, 
then $\pi$ is called a \emph{split cover of $\fb$}~\cite{EvansIvanovMacpherson} (in this case, we also say that \emph{$\pi$ is split}).
\item If $\fa$ is 
a covering reduct of a strongly trivial covering of $\fb$ with respect to $\pi$, 
then $\pi$ is called a \emph{strongly split cover of $\fb$}~\cite{EvansIvanovMacpherson}.
\end{itemize}

\end{definition}

Equivalently (and this motivates the terminology; see~\cite{EvansIvanovMacpherson}), 
a finite cover $\pi \colon \fa \to \fb$ is split if the kernel $K$ of $\mu_{\pi} \colon \aut(\fa) \to \aut(\fb)$ has a closed complement in
$\aut(\fa)$, i.e., there is a closed subgroup $H$ of $\aut(\fa)$ such that $KH = \aut(\fa)$ and $K \cap H = \{1\}$ (so that $\aut(\fa)$ is isomorphic to the  semidirect product $K \rtimes H$). 
Examples~\ref{expl:trivial},~\ref{expl:principal}~\ref{expl:non-free}, and~\ref{expl:neither-free-nor-trivial} 
are examples of split covers of $\set$.
For a non-example, see, e.g.,~\cite{EvansPastori}.
Example~\ref{expl:twisted}, in the case that $|S|=|T|=1$, is an example of a finite split cover of a structure in $\unary$ which is not strongly split.

\subsection{Operations on classes of structures}
Let $\fa$ be a structure, and let $\fb$ be a first-order  reduct of $\fa$. Then we say that $\fb$ is a \emph{finite index (first-order) reduct} of $\fa$ iff the index $[\aut(\fb):\aut(\fa)]$ is finite. We define the following operations on classes of structures.
	
\begin{definition}
	Let $\fa$ be a countable $\omega$-categorical structure. Then
\begin{itemize}
\item $\const(\fa)$ is the class of structures which are interdefinable with an expansion of $\fa$ with finitely many constants,
\item $\mcore(\fa)$ is the class of structures that are interdefinable with the (up to isomorphism unique~\cite{Cores-Journal,BodHilsMartin-Journal}) model-complete core of $\fa$,
\item $\reduct(\fa)$ is the class of first-order reducts of $\fa$,
\item $\reductfin(\fa)$ is the class of finite index first-order reducts of $\fa$,
\item $\fcover(\fa)$ is the class of finite covering structures of $\fa$.
\end{itemize}
	If $\mathcal{C}$ is a class of structures and $\Phi$ is one of the operators above, then we use the notation $\Phi(\mathcal{C})$ for the union of the classes $\Phi(\fa)$ such that $\fa\in \mathcal{C}$.
\end{definition}

\begin{proposition}\label{prop:easy}
The following identities hold. 
\begin{enumerate}
\item $\const\circ \const=\const$,
\item $\mcore\circ \mcore=\mcore$,
\item $\reduct\circ \reduct=\reduct$,
\item $\reductfin\circ \reductfin=\reductfin$,
\item $\fcover\circ \fcover=\fcover$,
\item $\const(\unary_{\nf}) = \unary^*$,
\item $\reduct(\unary) = \reduct(\unary^*)$, 
\item $\kexp=\reduct(\kexp)
$,
\item $\kexpp=\reduct(\kexpp)
$.
\end{enumerate}
\end{proposition}
\begin{proof}
Straightforward from the definitions. 
\end{proof}

	We will show that $\kexp=\reduct(\unary)$ and $\kexpp=(\fcover\circ \reduct)(\unary) = (\reduct \circ \fcover)(\unary)$, 
	and we will give several equivalent descriptions of these classes in Section~\ref{sect:additional}. We also prove Thomas' conjecture for each structure in $\kexpp$ (Theorem~\ref{thomas_strong}).

%% file: RU.tex

\begin{section}{Reducts of Unary Structures}
	In this section we characterise first-order reducts of unary structures in terms of their automorphism groups, and in particular prove Thomas' conjecture for the class $R(\unary)$. We mention that the finite-domain constraint satisfaction tractability conjecture has been shown for all structures in $R(\unary)$~\cite{BodMot-Unary}. 

\begin{lemma}\label{unary}
	Let $\fa$ be a structure. 
	Then $\fa \in \unary$ if and only if there are finitely may sets $O_1,\dots,O_k$ such that 
	$\aut(\fa)=\prod_{i=1}^k{\sym(O_i)}$. 
\end{lemma}

\begin{proof}
	First suppose that $\fa \in \unary$.
	Then $\fa$ is interdefinable with a unary structure $\fa'$; 
	let $O_1,\dots,O_k$ be the minimal non-empty intersections of predicates from $\fa'$;
	clearly, these sets partition $A$. 
	The containment $\aut(\fa) \subseteq \prod_{i=1}^k{\sym(O_i})$ is clear since
	every automorphism of $\fa$ is an automorphism of $\fa'$ and hence preserves
	the sets $O_1,\dots,O_k$. 
	For the reverse containment, let $\alpha \in \sym(A)$ 
	 be such that $\alpha(O_i) = O_i$ for all $i \leq k$. 	
	Since $\aut(\fa) = \aut(\fa')$ 
	 we need to show that $\alpha$ preserves all
	 (unary) relations $U$ of $\fa'$.  
	 Let $x \in U$ and let $i \leq k$ be such that
	 $x \in O_i$. Since $U$ must be a union
	 of orbits of $\aut(\fa')$ we have that
	 $\alpha(x) \in O_i$ also lies in $U$. 
	 
	 Conversely, suppose that $\aut(\fa)=\prod_{i=1}^k{\sym(O_i)}$. Then $\fa$ is first-order interdefinable with the unary structure $\fa=(X;O_1,\dots,O_k)$.
\end{proof}

	The following statement is an easy consequence of Lemma~\ref{unary} and Lemma~\ref{sing_or_inf}.

\begin{corollary}\label{unary_sing_inf}
	Let $\fa \in \unary$. Then the $\nabla(\fa)$-classes are the infinite orbits of $\aut(\fa)$ and the singleton orbits.
\end{corollary}

\begin{lemma}\label{unary_reduct}
	Let $\fa \in R(\unary)$ and let $C_1,\dots,C_k$ be the $\nabla(\fa)$-classes. Then $\prod_{i=1}^k{\sym(C_i)} \subseteq \aut(\fa)$.
\end{lemma}

\begin{proof}
	Let $\fa$ be a first-order reduct of a structure $\fb \in \unary$. 
	Let $O_1,\dots,O_l$ be the orbits of $\fb$. Then $\aut(\fb)=\prod_{i=1}^l{\sym(O_i)}$ by Lemma \ref{unary}. Let us define the binary relation $R$ on $A$ so that $xRy$ iff $x=y$ or the transposition $(xy)$ is contained in $\aut(\fa)$. Then it is easy to see that $R$ is a congruence of $\aut(\fa)$. On the other hand, $\aut(\fa) \supseteq \aut(\fb) = \prod_{i=1}^l{\sym(O_i)}$ implies that each class of $R$ is the union of some of the orbits of $\aut(\fb)$.
	In particular, $R$ has finitely many classes, and so by definition the congruence $\nabla(\fa)$ is finer that $R$. This means that for all $x,y\in C_i$ the transposition $(xy)$ is contained in $\aut(\fa)$. Therefore $\prod_{i=1}^k{\sym(C_i)} \subseteq \aut(\fa)$ since $\aut(\fa)$ is closed.
\end{proof}

\begin{corollary}\label{ru_rfiniteu}
	$\reduct(\unary)=\reductfin(\unary)=\reductfin(\unary^*)$.
\end{corollary}

\begin{proof}
	The containments ``$\supseteq$'' are obvious. Let $\fa\in \reduct(\unary)$. Let $C_1,\dots,C_k$ be the classes of $\nabla(\fa)$. Then the group $\aut(\fa)$ acts on the set $\{C_1,\dots,C_k\}$. By Lemma~\ref{unary_reduct} the kernel of this action is $\prod_{i=1}^k{\sym(C_i)}$. In particular, the index of $\prod_{i=1}^k{\sym(C_i)}$ in $\aut(\fa)$ is finite. On the other hand, $\prod_{i=1}^k{\sym(C_i)}$ is the automorphism group of a unary $\omega$-categorical structure with orbits $C_1,\dots,C_k$. We also know from Lemma~\ref{sing_or_inf} that each class $C_i$ is either a singleton or infinite. Therefore $\fa\in \reductfin(\unary^*)$. 
\end{proof}

The following has been shown in~\cite{BodMot-Unary} (Proposition 6.8); the proof we present here is simpler. 

\begin{corollary}\label{add_constants_ru}
	Let $\fa \in R(\unary)$. Then there exists an expansion of $\fa$ with finitely many constants which is in $\unary^*$.
\end{corollary}

\begin{proof}
	Let $C_1,\dots,C_k$ be the classes of $\nabla(\fa)$, and let us choose elements $c_i\in C_i$. We claim that the structure $(\fa,c_1,\dots,c_k)$ is in $\unary^*$. By Lemma~\ref{sing_or_inf} we know that each class $C_i$ is either a singleton or infinite. Without loss of generality we can assume that $C_i=\{c_i\}$ for $i=1,\dots,l$ and $C_j$ is infinite for $j>l$. We claim that $$\aut(\fa;c_1,\dots,c_k)=\prod_{i=1}^k{\id(\{c_i\})}\times \prod_{i=l+1}^k{\sym(C_i\setminus \{c_i\})}.$$ 
	Then Lemma~\ref{unary} implies that $(\fa,c_1,\dots,c_k)\in \unary^*$.
	To prove the claim, first recall from Lemma \ref{unary_reduct} that $\prod_{i=1}^k{\sym(C_i)} \subseteq \aut(\fa)$,
	and hence $\prod_{i=1}^k{\sym(C_i)}_{\{c_1,\dots,c_k\}} \subseteq \aut(\fa)_{\{c_1,\dots,c_k\}}$. Since every automorphism of $\fa$ 
	that fixes $c_1,\dots,c_k$ must also preserve the sets $C_1,\dots,C_k$ we in fact have equality 
	$\prod_{i=1}^k{\sym(C_i)}_{\{c_1,\dots,c_k\}} = \aut(\fa)_{\{c_1,\dots,c_k\}}$.
	Thus, \begin{align*}
	\aut(\fa;c_1,\dots,c_k) 
	& = \aut(\fa)_{\{c_1,\dots,c_k\}} \\
& = \prod_{i=1}^k{\sym(C_i)}_{\{c_1,\dots,c_k\}} =\prod_{i=1}^k{\id(\{c_i\})}\times \prod_{i=l+1}^k{\sym(C_i\setminus \{c_i\})}.\end{align*} 
\end{proof}

\begin{lemma}\label{unary_reduct2}
	Let $\fa \in R(\unary)$ and let $C_1,\dots,C_k$ be the $\nabla(\fa)$-classes. Then $\aut(\fa)=\prod_{i=1}^k{\sym(C_i)}\rtimes A$, where $A$ is a subgroup of $\aut(\fa)$ acting faithfully on $\{C_1,\dots,C_k\}$. 
\end{lemma}

\begin{proof}
	Recall from Lemma~\ref{sing_or_inf} 
	that every class of $\nabla(\fa)$ is either infinite or a singleton. We can assume that $C_1,\dots,C_l$ are infinite and that for every $j \in \{l+1,\dots,k\}$ there exists $c_j \in A$ such that $C_j=\{c_j\}$. Let $e_i \colon C_i\rightarrow \mathbb{N}$, for $1\leq j\leq l$, be bijections and let $e=\bigcup_{i \in \{1,\dots,l\}} e_i$.  Then it is easy to see that for every $\alpha \in \aut(\fa)$ there exists $\tilde{\alpha} \in \sym(A)$ such that 
\begin{enumerate}
\item $\alpha$ and $\tilde{\alpha}$ have the same action on the set $\{C_1,\dots,C_k\}$,
\item $e(\tilde{\alpha}(x))=e(x)$ for every $x\in \bigcup_{i \in \{1,\dots,l\}} C_i$. 
\end{enumerate}
	The permutation $\alpha^{-1}\tilde{\alpha}$ fixes every class of $\nabla(\fa)$. 
	Since $N := \prod_{i=1}^k{\sym(C_i)} \subseteq \aut(\fa)$ by Lemma~\ref{unary_reduct} it follows that $\alpha^{-1}\tilde{\alpha}\in \aut(\fa)$ and thus $\tilde{\alpha}\in \aut(\fa)$. Let $A:=\{\tilde{\alpha} \mid\alpha\in \aut(\fa)\}$. Then $A \subseteq \aut(\fa)$ is a subgroup of $G$ 
	which acts faithfully on $\{C_1,\dots,C_k\}$. Then it is also clear that $N$ is a normal subgroup of $\aut(\fa)$ since it is the kernel of the action of $\aut(\fa)$ on 
	$\{C_1,\dots,C_k\}$. 
	Therefore, $\aut(\fa)$ can be written as a semidirect product, $\aut(\fa)=\prod_{i=1}^k{\sym(C_i)}\rtimes A$.	
\end{proof}

\begin{corollary}\label{wreath_infinite}
	Let $\fa \in R(\unary)$ be with no finite orbits. Then $\aut(\fa)$ is isomorphic to the wreath product $\sym(\mathbb{N})\wr A$ for some permutation group $A$ on a finite set.
\end{corollary}

\begin{proof}
	Let $C_1,\dots,C_k$ be the classes of $\nabla(\fa)$. Without loss of generality we can assume that $C_i=\{(i,n) \mid n\in \mathbb{N}\}$. Let $A$ be the image of the action of $\aut(\fa)$ on the set $\{C_1,\dots,C_k\}$. Then if we use the bijections $e_i \colon (i,n)\mapsto n$ in the proof of Lemma \ref{unary_reduct2} the statement of the Corollary follows.
\end{proof}

\begin{lemma}\label{unary_reduct4}
	Suppose that $\fa,\fb \in R(\unary)$ have the same domain.
	If $\nabla(\fa)=\nabla(\fb)$ and the actions of the groups $\aut(\fa)$ and $\aut(\fb)$ on the $\nabla(\fa)$-classes are the same, then $\fa$ and $\fb$ are interdefinable.
\end{lemma}

\begin{proof}
	By the $\omega$-categoricity of $\fa$ and $\fb$ it is enough to show that $\aut(\fa)=\aut(\fb)$. Let $C_1,\dots,C_k$ be the classes of $\nabla(\fa)=\nabla(\fb)$. 
	Lemma~\ref{unary_reduct} shows that $\prod_{i=1}^k{\sym(C_i)}\subseteq \aut(\fa)$. Now let $\beta \in \aut(\fb)$. By our assumption 
	about the action of $\aut(\fa)$ and $\aut(\fb)$ on the $\nabla(\fa)$-classes there is a permutation $\alpha \in \aut(\fa)$ such that $\beta \alpha^{-1}$ fixes each class $C_i$. Then $\beta \alpha^{-1}\in \prod_{i=1}^k{\sym(C_i)}\subseteq \aut(\fa)$, and so $\beta \in \aut(\fa)$. Therefore $\aut(\fb)\subseteq \aut(\fa)$. Analogously, $\aut(\fa)\subseteq \aut(\fb)$.  
\end{proof}

\begin{corollary}\label{thomas_unary}
	Every structure in $\unary$ has finitely many first-order reducts.
\end{corollary}

\begin{proof}
	Let $\fa \in \unary$. Then by 
	Lemma~\ref{unary} $\aut(\fa)=\prod_{i=1}^k{\sym(O_i)}$. If $\fb$ is a first-order reduct of $\fa$ then 
	$\nabla(\fb)$ is a union of orbits of $\fa$. This means that there are finitely many choices for the relation $\nabla(\fb)$. If the relation $\nabla(\fb)$ is fixed then there are finitely many possible actions of $\aut(\fb)$ on the classes of $\nabla(\fb)$. By Lemma \ref{unary_reduct4} it follows that $\nabla(\fb)$ and the action of $\aut(\fb)$ on the classes of $\nabla(\fb)$ already determine the structure $\fb$ up to interdefinability. Therefore, 
	$\fa$ has finitely many first-order reducts.
\end{proof}

We also obtain an equivalent description of first-order reducts of unary structures in terms of their automorphism groups.

\begin{corollary}\label{cor:unary_reduct-iff}
	A structure $\fa$ is in $\reduct(\unary)$ if and only if  
	$\prod_{i=1}^k{\sym(C_i)} \subseteq \aut(\fa)$ for some partition of $A$ into classes $C_1,\dots,C_k$.
	\end{corollary}
\begin{proof}
One direction has been shown in Lemma~\ref{unary_reduct}. For the converse implication, suppose that $\aut(\fa)$ contains
$\prod_{i=1}^k{\sym(C_i)}$ for some partition of $A$ into classes $C_1,\dots,C_k$. 
Note that $\aut(A;C_1,\dots,C_k) = \prod_{i=1}^k{\sym(C_i)}$ (see Lemma~\ref{unary}) and that $\fa$ is a first-order reduct of $(A;C_1,\dots,C_k)$. Hence, $\fa \in \reduct(\unary)$.
\end{proof}

\end{section}

%% file: FU.tex

\section{Finite Coverings of Unary Structures}

	In this section we classify the finite coverings of unary structures. First we make the following observation.

\begin{lemma}\label{fu_fustar}
	$\fcover(\unary)=\fcover(\unary^*)$.
\end{lemma}

\begin{proof}
	The containment ``$\supseteq$'' 
	is trivial. In order to show the other direction it is enough to show that $\unary\subseteq \fcover(\unary^*)$ since $\fcover\circ \fcover=\fcover$. So let $\fa \in \unary$ and let $F$ be the union of the finite orbits of $\fa$. Then $F$ is finite. Let us consider the unary structure $\fb$ whose domain is $B:=A\setminus F\cup \{x\}$ for any $x \notin A$, and whose relations are the infinite orbits of $\fa$ and $\{x\}$. Then $\fb\in \unary^*$. Let $\pi \colon \fa\rightarrow \fb$ be defined as $\pi(y)=x$ if $x\in F$, and $\pi(y)=y$ otherwise. Then it is easy to see that $\pi$ is a finite covering map, and hence $\fa\in \fcover(\unary^*)$.
\end{proof}

The following theorem summarises the results from Section~\ref{sect:split} and Section~\ref{sect:cov-reducts-triv-covers}.

\begin{theorem}\label{reduct_trivial_finite1}
	Let $\fb\in \unary^*$ and let $\pi \colon\fa\rightarrow \fb$ be a finite covering map. Then $\fa$ has finitely many covering reducts with respect to $\pi$.
\end{theorem}
\begin{proof}
Proposition~\ref{reduct_trivial2} shows that $\pi$
is strongly split. The statement then follows
from Corollary~\ref{reduct_trivial_finite}.
\end{proof}




\subsection{Finite covers of unary structures split}
\label{sect:split}
	The following series of lemmas is needed to show that every finite covering map 
	of a structure $\fb \in \unary^*$ is strongly split (Proposition~\ref{reduct_trivial2}).
Throughout this subsection, let $\fb\in \unary^*$
and let $\pi \colon\fa\rightarrow \fb$ be a finite covering map.

\begin{remark}\label{weak_strong}
	Observe that $\fb$ satisfies the condition of Lemma \ref{no_finite_index}, that is, $\aut(\fb)_x$ has no finite index subgroup for any $x\in B$. By Lemma \ref{no_finite_index} this implies that every trivial finite cover of $\fb$ is strongly trivial, and hence every split cover of $\fb$ is strongly split.
\end{remark}

\begin{remark} When $\fb$ is taken
from $\unary$ instead of $\unary^*$, then 
there are split covers of $\fb$ that are not strongly split, as illustrated by Example~\ref{expl:twisted} if $|S|=|T|=1$.
\end{remark} 

\begin{lemma}\label{generating_cycle}
Let $F$ be a finite subset of an infinite orbit $O$ of $\fb$. If $|F|$ is large enough then there exists an automorphism $\alpha$ of $\fa$ such that 
\begin{enumerate}
\item $\alpha(x)=x$ for all $x \in E := A \setminus \pi^{-1}(F)$,
\item $\mu_\pi(\alpha)|_F$ is nontrivial.
\end{enumerate}
\end{lemma}

\begin{proof}
	Let $k$ be the maximum of the sizes of the fibers of $\pi$ and 
	let $p>k$ be a prime number. We claim that if $|F|\geq p$ then there is an automorphism $\alpha$ of $\fa$ satisfying Conditions (1) and (2). 

	Let $u_1,\dots,u_p\in F$ be distinct elements. Then the $p$-cycle $(u_1u_2\dots u_p)$ is contained in $\aut(\fb)$ by Lemma \ref{unary}. By the definition of finite covering maps there exists $\beta\in \aut(\fa)$ such that $\mu_\pi(\beta)=(u_1\dots u_p)$. Now let $\alpha := \beta^{k!}\in \aut(\fa)$. Then $\mu_\pi(\alpha)|_F$ is again a $p$-cycle and hence nontrivial. On the other hand, if $u\in B \setminus F$ and $U := \pi^{-1}(u)$ then $\beta(U)=U$, and $\alpha|_U=\beta^{k!}|_U=\id_U$ since $|U|\leq k$. This means that $\alpha|_E=\id_E$. Therefore, $\alpha\in \aut(\fa)$ satisfies the conditions (1) and (2) which proves the lemma.
\end{proof}

Recall that for any finite set $F$ of cardinality at least 5, the alternating group 
$\alt(F)$ is the only non-trivial proper normal subgroup of 
$\sym(F)$ 
(see e.g.~Chapter 8.1 in~\cite{DixonMortimer}). 

\begin{lemma}\label{generating_trans}
Let $F$ be a finite subset of an infinite orbit $O$ of $\fb$. If $|F|$ is large enough then for any pairwise distinct $u_1,u_2,u_3,u_4 \in F$ there exists an automorphism $\alpha$ of $\fa$ such that 
\begin{enumerate}
\item $\alpha(x)=x$ for all $x\in E := A \setminus \pi^{-1}(F)$,
\item $\mu_\pi(\alpha)|_F=(u_1u_2)(u_3u_4)$.
\end{enumerate}

\end{lemma}
\begin{proof}
	Let $K:=\{\mu_\pi(\gamma)|_F \in \sym(F) \mid \gamma\in \aut(\fa)_E\}$.
	We claim that $K$ is a normal subgroup of $\sym(F)$. It is clear that $K$ is a subgroup of $\sym(F)$. Let $\alpha\in K$ and $\beta\in \sym(F)$. We need to show that $\beta\alpha\beta^{-1}\in K$. 
	By the definition of $K$ 
	there exists $\gamma\in \aut(\fa)_E$ so that $\mu_\pi(\gamma)|_F=\alpha$. 
	By Lemma~\ref{unary} 
	there exists $\delta \in \aut(\fb)$ so that 
	$\delta|_F = \beta$. By the definition of finite covers, there exists $\eta \in \aut(\fa)$ such that 
	$\mu_\pi(\eta)=\delta$. 
	Let $\gamma'=\eta \gamma\eta ^{-1} \in \aut(\fa)$. Then on can check that $\gamma'(x)=x$ for all $x\in E$ and 
	\begin{align*}
	\mu_\pi(\gamma')|_F=(\mu_\pi(\eta)\mu_\pi(\gamma)\mu_\pi(\eta)^{-1})|_F=\delta
	|_F\mu_\pi(\gamma)|_F \delta^{-1}|_F=\beta\alpha\beta^{-1}.
	\end{align*}

	We obtained that $K\triangleleft \sym(F)$. By Lemma \ref{generating_cycle} we know that if $F$ is large enough, then $K$ is nontrivial. Therefore if $F$ is large enough, then $K\geq \alt(F)$, and then the statement of the lemma follows.
\end{proof}

\begin{lemma}\label{generating_trans2}
	Let $O$ be an infinite orbit of $\fb$. Then for all distinct $v_1,v_2 \in O$ there exists an $\alpha \in \aut(\fa)$ such that
\begin{enumerate}
\item $\alpha(x)=x$ for all $x\in A \setminus \pi^{-1}(\{v_1,v_2\})$,
\item $\mu_\pi(\alpha)=(v_1v_2)$.
\end{enumerate}
\end{lemma}
\begin{proof}
	 Let $v_1,v_2 \in O$, and let us choose a finite subset $F$ of $O$ which contains the elements $v_1$ and $v_2$ which is large enough so that we can apply Lemma \ref{generating_trans} for $F$.
	Choose $u_3,u_4 \in F$ such that $u_1 := v_1,u_2 := v_2,u_3,u_4$ are pairwise distinct.
	Let $\alpha \in \aut(\fa)$ be as in Lemma \ref{generating_trans}. For each $i \in {\mathbb N}$, choose $\gamma_i\in \aut(\fa)$ so that 
	\begin{align*}
	\mu_\pi(\gamma_i)(v_1) & =v_2, \\
	\mu_\pi(\gamma_i)(v_2) & =v_1, \text{ and } \\
	\mu_\pi(\gamma_i)(F)\cap \mu_\pi(\gamma_j)(F) & =\{v_1,v_2\} \text{ for all } i\neq j.
	\end{align*}
	By Lemma \ref{unary} it follows that such $\gamma_i$'s exist. Let $\beta_i:=\gamma_i\alpha\gamma_i^{-1}\in \aut(\fa)$. Then $\mu_\pi(\beta_i)(v_1)=v_2$, $\mu_\pi(\beta_i)(v_2)=v_1$, and for all $x\in S:= \pi^{-1}(\{v_1,v_2\})$ we have $\beta_i(x)=x$ if $i$ is large enough. Since there are finitely many possible actions of $\beta_i$ on the finite set $S$ there is a subsequence $(\beta_{l(i)})_i$ of $(\beta_{i})_i$ so that $\beta_{l(i)}|_S=\beta_{l(j)}|_S$ for all $i,j\in \mathbb{N}$. Then the sequence $\beta_{l(i)}$ converges to a permutation $\beta$ for which $\beta(x)=x$ for all $x\in A \setminus S$ and $\pi(\beta)=\pi(\beta_{l(1)})=(v_1v_2)$. Since $\aut(\fa)$ is closed it follows that $\beta\in \aut(\fa)$ which finishes the proof of the lemma. 
\end{proof}

\begin{lemma}\label{reduct_trivial11}
	 Let $\mathcal{O}$ be the set of orbits of $\fb$. Then for each $O\in \mathcal{O}$ there exists a finite set $F_O$ and a mapping $\psi_O \colon \pi^{-1}(O) \rightarrow F_O$ such that $\psi_O|_{\pi^{-1}(y)}$ is bijective for every $y\in O$, and $\aut(\fa)$ contains every $\alpha \in \sym(A)$ such that 
\begin{enumerate}
\item $\alpha$ preserves $\sim_{\pi}$,
\item $\mu_\pi(\alpha)\in \aut(\fb)$,
\item $\psi_O(x)=\psi_O(\alpha(x))$ for every $O\in \mathcal{O}$ and $x\in \pi^{-1}(O)$. 
\end{enumerate}
\end{lemma}
\begin{proof}
	If $O$ is finite, then $O=\{u\}$ for some $u\in B$ since $\fb\in \unary^*$. In this case let $\psi_O=\id_{\pi^{-1}(u)}$.
If $O$ is infinite, then we define $\psi_O$ as follows. Let $u\in O$ be arbitrary. 
\begin{itemize}
\item If $x\in \pi^{-1}(u)$ then set $\psi_O(x):=x$.
\item If $x\in \pi^{-1}(O) \setminus \pi^{-1}(u)$ then by Lemma \ref{generating_trans2} there exists a permutation $\alpha\in \aut(\fa)|_{A \setminus \pi^{-1}(\{\pi(x),u\})}$ such that $\alpha(\pi^{-1}(\pi(x)))=\pi^{-1}(u)$. In particular, $\alpha$ defines a bijection between $\pi^{-1}(\pi(x)))$ and $\pi^{-1}(u)$. Set $\psi_O(x):=\alpha(x)$.
\end{itemize}

	We claim that these mappings satisfy the conditions of the lemma. Let $G$ be the permutation group of those $\gamma \in \sym(B)$ for which there exists an automorphism $\alpha$ of $\fa$ with $(\mu_\pi(\alpha))=\gamma$ and satisfying Conditions (1)-(3) of the lemma. Then since $\aut(\fb)=\prod_{O\in \mathcal{O}}{\sym(O)}$ it is enough to show that $(G_{B\setminus O})|_O=\sym(O)$ for all $O\in \mathcal{O}$. If $O$ is a singleton, then the claim is trivial, so we can assume that $O$ is infinite. It is easy to see that $G$ is closed. Thus, $(G_{B\setminus O})|_O$ is also closed. Hence, by Lemma \ref{unary} it is enough to show that $G$ contains for all $u_1,u_2\in O$ the transposition $(u_1u_2)$. For this it is enough to show that $(uv) \in G$ for all $v\in O\setminus \{u\}$, where $u$ is the element of $O$ which is used in the definition of $\psi_O$. But this follows directly from the definition of the mapping $\psi_O$. 
\end{proof}

\begin{proposition}\label{reduct_trivial2}
Any finite covering map $\pi \colon \fa \rightarrow \fb$ for $\fb \in \unary^*$ is strongly split.
\end{proposition}

\begin{proof}
	 Let $F_O$ and $\psi_O$ be defined as in Lemma \ref{reduct_trivial11} for each orbit $O$ of $\fb$. 
	Let $F:=\bigsqcup_O{F_O}$ and $\psi:=\bigcup_O{\psi_O}$. Let
	$\fa'$ be the expansion of $\fa$ obtained by adding to $\fa$ for each $x\in F$ the unary relation 
	$\psi^{-1}(x)$. Then by Lemma \ref{reduct_trivial11} it follows that $\mu_\pi(\aut(\fa'))=\aut(\fb)=\mu_\pi(\aut(\fa))$. Thus $\fa$ is a covering reduct of $\fa'$. We claim that $\pi \colon \fa'\rightarrow \fb$ is a strongly trivial cover. This implies the statement of the proposition. By Remark \ref{weak_strong} it is enough to show that the finite cover $\pi \colon \fa'\rightarrow \fb$ is trivial, i.e., that the kernel of the map $\mu_\pi$ is trivial. Let $\alpha\in \aut(\fa')$ be so that $\mu_\pi(\alpha)=\id_B$ and let $x\in A$. Then $x\sim_{\pi} \alpha(x)$ and $\psi(x)=\psi(\alpha(x))$. It follows from the definition that $\psi$ is injective on $[x]_\pi$. This implies that $x=\alpha(x)$ and hence that $\alpha=\id_A$. Therefore, the kernel of $\mu_\pi$ is trivial.
\end{proof}

\begin{remark}
Proposition~\ref{reduct_trivial2} generalises Theorem 2.4 in~\cite{FiniteCovers}, which states that every finite covering of $\set$ strongly splits.
\end{remark}

\subsection{Covering reducts of trivial coverings}
\label{sect:cov-reducts-triv-covers}
	In this subsection we describe the automorphism groups of covering reducts of a trivial finite covering of a structure in $\unary^*$. In particular, we show that there are always finitely many of them.

\begin{remark}\label{triv_cov_unary}
	Let $\fa$ be a strongly trivial covering of a unary structure $\fb$ with orbits $O_1,\dots,O_k$. Then as in Remark \ref{rem:triv-covers} the elements of the structure $\fa$ can be identified with the elements of $\bigsqcup_{i=1}^k{(F_i\times O_i)}$ for some finite sets $F_i$ so that $\aut(\fa)$ consists of exactly those permutations which preserve the first coordinate and stabilise the sets $O_i$ in the second coordinate. In this case $\aut(\fa)$ can be written as $\prod_{i \in \{1,\dots,k\}} \{\id_{F_i}\} \wr \sym(O_i)$. For convenience we will always assume that the sets $F_i$ are paiwise disjoint.
\end{remark}

\begin{remark}\label{rem:triv-cover-interpret}
It follows from the description of strongly trivial coverings of unary structures in Remark~\ref{triv_cov_unary} 
that every (reduct of a) strongly trivial covering structure of a structure from $\unary$ has a first-order interpretation over $({\mathbb N};=)$. 
\end{remark}

	Throughout this subsection let us fix a structure $\fb\in \unary^*$ and a trivial finite covering map $\pi \colon\fa\rightarrow \fb$. Let $O_1,\dots,O_k$ be the orbits of $\fb$. 
	We identify the elements of the structure $\fa$ with the elements of $\bigsqcup_{i=1}^k{(F_i\times O_i)}$ for some disjoint finite sets $F_1,\dots,F_k$ as explained in Remark \ref{triv_cov_unary}, and define $F:=\bigsqcup_{i=1}^kF_i$.


\begin{definition}\label{group_ns}
\begin{itemize}
\item 
Let $S$ be the group of all permutations of $A$ which fix the sets $F_i\times O_i$ for $i \in \{1,\dots,k\}$ and which preserve the congruence $\sim_{\pi}$.
\item Let $N$ be the group of all permutations of $A$ which fix all fibers setwise (i.e., $N$ is the kernel of the map $\mu_\pi \colon S\rightarrow \aut(\fb)$).
\end{itemize}
\end{definition}

	The following statements are direct consequences of the definitions above.

\begin{proposition}\label{semidir}
\begin{enumerate}
\item A first-order reduct $\fc$ of $\fa$ is a covering reduct of $\fb$ with respect to the covering $\pi$ if and only if $\aut(\fc) \subseteq S$.
\item The group $S$ can be written as a semidirect product $N \rtimes \aut(\fa)$.
\end{enumerate}
\end{proposition}

\begin{proof}
	(1) follows easily from the definition using that $\aut(\fb)=\prod_{i=1}^k{\sym(O_i)}$.
	Since $N$ is the kernel of the homomorphism $\mu_{\pi} \colon S\rightarrow \aut(\fb)$ we have $N\triangleleft S$. It is obvious that 
	$\aut(\fa) \leq S$ and that 
	$S = N \aut(\fa)$.
	Since $\pi$ is a trivial covering map it follows that the kernel of 
	$\mu_{\pi}|_{\aut(\fa)}$ is trivial, 
	that is, $N\cap \aut(\fa)=\{\id_A\}$.
\end{proof}

	The following lemma is a direct consequence of item (2) of Proposition \ref{semidir}.
	Let $H$ and $K$ be subgroups of the same group. Then we say that $H$ \emph{normalises} $K$ if 
	$H$ is a subgroup of the normaliser of $K$, i.e., for every $h \in H$ we have that $$\{h^{-1} k h \mid k \in K\} = K.$$
	
\begin{lemma}\label{semidir2}
	The mapping $G \mapsto G \cap N$ defines a bijection between the closed subgroups of $S$ that contain $\aut(\fa)$
	and the closed subgroups of $N$ which are normalized by $\aut(\fa)$. 
	The inverse map is $H \mapsto H \rtimes \aut(\fa)$.
\end{lemma}
\begin{proof}
If $H$ is a subgroup of $N$ which is normalized by $\aut(\fa)$ then the group generated by $H$ and $\aut(\fa)$ can be written as a product $H\aut(\fa)$. Since $H\cap \aut(\fa) \subseteq N\cap \aut(\fa)=\{\id_A\}$, it follows that this group can be written as a semidirect product $H \rtimes \aut(\fa)$. Then $(H \rtimes \aut(\fa))\cap N=H$.

We claim that if $H$ is closed then so is $H \rtimes \aut(\fa)$. Let $\alpha_1,\alpha_2, \ldots \in H \rtimes \aut(\fa)$ be a sequence converging to some $\alpha\in \sym(A)$. Let $\beta_i$ (and $\beta)$ be the unique element in $\aut(\fa)$ for which $\mu_{\pi}(\beta_i)=\mu_{\pi}(\alpha_i)$ (and $\mu_{\pi}(\beta)=\mu_{\pi}(\alpha)$), that is, $\alpha_i\beta_i^{-1}\in H$ (and $\alpha\beta^{-1}\in H$). Since $\mu_{\pi}$ is continuous it follows that $(\beta_i)_i$ converges to $\beta$. Hence the sequence $(\alpha_i\beta_i^{-1})_i$ converges to $\alpha\beta^{-1}$. Since $\alpha_i\beta_i^{-1}\in H$ and $H$ is closed it follows that $\alpha\beta^{-1}\in H$ and hence $\alpha\in H\aut(\fa)=H \rtimes \aut(\fa)$. Therefore, $H \rtimes \aut(\fa)$ is closed.

Let $G$ be a subgroup of $S$ containing $\aut(\fa)$. We claim that $\aut(\fa)$ normalises $G \cap N$. Indeed, let $g \in \aut(\fa)$. Then $\{g^{-1} h g \mid h \in G \cap N\} = G \cap N$ since $N$ is a normal subgroup of $S$ which contains $\aut(\fa)$. Since $(G\cap N)\cap \aut(\fa)\subseteq N\cap \aut(\fa)=\{\id_A\}$, it follows that $G$ can be written as $G=(G\cap N) \rtimes \aut(\fa)$. Moreover, it is clear that if $G$ is closed, then so is $G\cap N$.

	We have obtained that the mappings defined in the lemma are inverses of each other, which also implies that they both are bijections. 
\end{proof}

Hence, in order to classify the covering reducts of $\fa$ it is enough to classify those closed subgroups of $N$ which are normalized by $\aut(\fa)$.
	If $K$ is a normal subgroup of a group $G$, then we write that two elements $x_1,x_2 \in G$ \emph{are the same modulo $K$} if they represent the same element in the factor group $G/K$, i.e., if $x_1x_2^{-1}\in K$.

\begin{definition}\label{cond_hn}
	Let $H$ be a subgroup of 
	$\prod_{i=1}^k{\sym(F_i)} \subseteq \sym(F)$. 
	For $i \in \{1,\dots,k\}$ let $H_i:=H_{F\setminus F_i}|_{F_i}\subseteq \sym(F_i)$ 
	 and let $N_i$ be a subgroup of $H_i$ normalised by $H|_{F_i}$ for every $i \in \{1,\dots,k\}$ so that if $|O_i|=1$ then $N_i=H_i$.
	 We write $N(H,N_1,\dots,N_k)$ for the group of all permutations $\alpha \in N$ such that 
\begin{itemize}
\item for every $i \in \{1,\dots,k\}$ and for every $x_i\in O_i$ there is a permutation $\gamma \in H$ such that the action of $\alpha$ on the first coordinate of the fiber $F_i\times \{x_i\}$ is exactly the $i$-th coordinate of $\gamma$, and 
\item for every $i \in \{1,\dots,k\}$ and $x,y\in O_i$ the actions of $\alpha$ on the first coordinate of the fibers $F_i\times \{x\}$ and $F_i\times \{y\}$ are the same modulo $N_i$.
\end{itemize}
\end{definition}

	It follows directly from the definition that $N(H,N_1,\dots,N_k)$ is a closed subgroup of $N$ and normalised by $\aut(\fa)$. We will show that subgroups of $N$ with these properties
	are of the form $N(H,N_1,\dots,N_k)$. 

\begin{definition}\label{def:h-of-g}
	Let $G$ be a subgroup of $N$ which is normalised by $\aut(\fa)$. 
\begin{itemize}
\item Let $H(G)$ 
be the subgroup of 
$\prod_{i=1}^k{\sym(F_i)}$ containing all permutations $\gamma$ such that 
there exists a permutation $\alpha \in G$ and 
for every $i \in \{1,\dots,k\}$ an element $x_i\in O_i$ 
such that the action of $\alpha$ on the first coordinate of the fiber $F_i\times \{x_i\}$ is exactly the $i$-th coordinate of $\gamma$.
\item $H_i(G):=H(G)_{F\setminus F_i}|_{F_i}$.
\item Let $N_i(G)$ be the group of all 
$\gamma \in \sym(F_i)$ such that there exists a permutation $\alpha\in G$ and an $x\in O_i$ such that the action of $\alpha$ on the first coordinate of the fiber $F_i\times \{x\}$ equals $\gamma$ and $\alpha$ fixes every element of $A \setminus (F_i\times \{x\})$.
\end{itemize}
\end{definition}

\begin{remark}
	Since $\aut(\fa)=\prod_{i=1}^k{\sym(O_i)}$
	normalises $G$ it does not matter which elements $x_i \in O_i$ we take in the definition of $H(G)$ and $N_i(G)$. It follows that $H(G)$, $H_i(G)$, and $N_i(G)$ are indeed groups.
\end{remark}

\begin{remark}\label{hn_equal}
	It is clear from the definition that $N_i(G)\leq H_i(G)$, and if $|O_i|=1$ then $N_i(G)=H_i(G)$.
\end{remark}

\begin{example}
A simple example for $k=1$ where $H_1 \neq N_1$ is Example~\ref{expl:non-free}, where
$N_1 = \{\id\}$ and $H_1 = {\mathbb Z}_2$.
\end{example}

\begin{lemma}\label{nh_going_back}
	Let $H,N_1,\dots,N_k$ be the groups as introduced in Definition \ref{cond_hn}. Then $H(N(H,N_1,\dots,N_k))=H$ and $N_i(N(H,N_1,\dots,N_k))=N_i$.
\end{lemma}

\begin{proof}
	It follows directly from the definition that $H(N(H,N_1,\dots,N_k))=H$. 
	For the second claim let us assume that $i \in \{1,\dots,k\}$, $x\in O_i$, and that $\alpha \in N$ fixes every element of $A\setminus (F_i\times \{x\})$. Let $H_i$ be the group as introduced in Definition \ref{cond_hn}. We have to show that $\alpha\in N(H,N_1,\dots,N_k)$ if and only if $\alpha$ acts on the first coordinate of the fiber $F_i\times \{x\}$ as a permutation in $N_i$.

	For the ``if'' let us assume that the action of $\alpha$ on the first coordinate of the fiber $F_i\times \{x\}$, which we denote by $\gamma_i$, is an element of $N_i$. We show that $\alpha\in N(H,N_1,\dots,N_k)$. First we show that $\alpha$ satisfies the first bullet of Definition~\ref{cond_hn}. Arbitrarily choose $j \in \{1,\dots,k\}$ and $x_j\in O_j$. Let $\gamma_j$ denote the action of $\alpha$ on the fiber $F_j\times \{x_j\}$, and let $\gamma=\prod_{j=1}^k\gamma_j$. We have to show that $\gamma\in H$. It follows from the definition of $\alpha$ that $\gamma_j=\id(F_j)$ if $j\neq i$. Therefore it is enough to show that $\gamma_i\in H_i=H_{F\setminus F_i}|_{F_i}$. If $x_i=x$ then $\gamma_i\in N_i$, and if $x_i\neq x$ then by definition $\gamma_i=\id(F_i)\in N_i$. In either case we have $\gamma_i\in N_i\leq H_i$. It remains to show that $\alpha$ satisfies the second bullet of Definition~\ref{cond_hn}. In fact we show something stronger. We show that for every $j\in \{1,\dots,k\}$ and $y\in O_j$ the action of $\alpha$ on $F_j\times \{y\}$ is an element of $N_j$. If $y=x$, and thus $j=i$, then this is exactly our assumption on $\alpha$. If $y\neq x$ then this action is trivial since by definition $\alpha$ fixes every element in $A\setminus (F_i\times \{x\})$.

	For the other direction let us assume that $\alpha\in N(H,N_1,\dots,N_k)$. We distinguish two cases. 
	If $|O_i|=1$ then by definition $\alpha$ acts on the first coordinate of the fiber $F_i\times \{x\}$ as a permutation in $H_i(N(H,N_1,\dots,N_k))$, and by Remark~\ref{hn_equal} we have 
\begin{align*}
N_i(N(H,N_1,\dots,N_k)) & =H_i(N(H,N_1,\dots,N_k)) \\
& = (H(N(H,N_1,\dots,N_k))_{F\setminus F_i})|_{F_i}
= H_{F\setminus F_i}|_{F_i}=H_i=N_i
\end{align*}
 where the last equality follows from the assumptions on the groups $H_i$ and $N_i$ in Definition~\ref{cond_hn}.
	
	If $|O_i|>1$ then let $y$ be an element in $O_i$ different from $x$. Then if $\alpha\in N(H,N_1,\dots,N_k)$ then on one hand we know that the action of $\alpha$ on the first coordinate of the fibers $F_i\times \{x\}$ and $F_i\times \{y\}$ are the same modulo $N_i$. On the other hand we know by the definition of $\alpha$ that its action on the fiber $F_i\times \{y\}$ is trivial. Therefore, the action of $\alpha$ on the first coordinate of the fiber $F_i\times \{x\}$ is an element of $N_i$.
\end{proof}
	
\begin{definition}
Let $N_i$ be a subgroup of $\sym(F_i)$ for some $i\in \{1,\dots,k\}$. Then $N^*(N_1,\dots,N_k)$ is defined to be the closure of the group generated by all permutations $\alpha$ for which there exists an $i$ and $x\in O_i$ such that the action of $\alpha$ on the first coordinate of the fiber $F_i\times \{x\}$ is in $N_i$ and $\alpha$ fixes every element of $A \setminus (F_i\times \{x\})$. 
\end{definition}
	
	It follows easily from the definition that $N^*(N_1,\dots,N_k)$ is contained in every closed group $G\leq N$ normalized by $\aut(\fa)$ with $N_i(G)=N_i$. It is also easy to see that in fact $N^*(N_1,\dots,N_k)=N(\prod_{i=1}^k{N_i},N_1,\dots,N_k)$ (using the notation from Definition \ref{cond_hn}).

\begin{example}
Let $\fb := \set$ and $\fa$ a strongly trivial finite
covering structure of $\fb$ with fibers of size four. 
Then the covering structure from Example~\ref{expl:neither-free-nor-trivial}
is a covering reduct $\fc$ of $\fa$. 
Let $G := \aut(\fc)$. 
As $G$ is transitive, we have $S = N$ and $k=1$ in Definition~\ref{def:h-of-g}. 
Then $H(G)=\prod_{i=1}^k {\mathbb Z}_4$ and $N_1(G) = {\mathbb Z}_2$.
\end{example}

\begin{lemma}\label{normal_hn}
	Let $G \leq N$ be normalised by $\aut(\fa)$. Then $N_i(G) \triangleleft H(G)|_{F_i}$. 
\end{lemma}

\begin{proof}
	Let $\alpha\in H(G)|_{F_i}$ and $\beta\in N_i(G)$. Let $\gamma\in G$ be an element witnessing $\alpha\in H(G)|_{F_i}$. Let $x\in O_i$ and let $\delta\in G$ be an element witnessing $\delta\in N_i(G)$ on the fiber $F_i\times \{x\}$. Then the element $\gamma^{-1}\delta\gamma\in G$ witnesses the fact that $\alpha^{-1}\beta\alpha\in N_i(G)$.
\end{proof}

\begin{lemma}\label{generating_n}
	Let $G$ be a closed subgroup of $N$ normalised by $\aut(\fa)$. Let $\alpha\in G$
	and $u,v\in O_i$. Then the actions of $\alpha$ on the first coordinate of the fibers $F_i\times \{u\}$ and $F_i\times \{v\}$ are the same modulo $N_i(G)$.
\end{lemma}

\begin{proof}
	Let $\alpha_u$ and $\alpha_v$ denote the action of $\alpha$ on the first coordinate of the fibers of $u$ and $v$, respectively, so $\alpha_u,\alpha_v \in \sym(F_i)$.
	For $\beta \in \aut(\fb)$, we write $\pi^{-1}(\beta)$ for the unique $\gamma \in \aut(\fa)$ such that
$\mu_{\pi}(\gamma) = \beta$.
	 Let $\beta=(uv)\in \aut(\fb)$. 
	Let $\gamma := \alpha^{-1}(\pi^{-1}(\beta))^{-1}\alpha\pi^{-1}(\beta)$. Then the action of $\gamma$ on the first coordinate of the fiber $F_i\times \{u\}$ is $\alpha_u^{-1}\alpha_v$. On the other hand, $\gamma$ fixes every element of $A \setminus F_i\times \{u,v\}$. Now let $v_1,v_2,\dots$ be pairwise distinct elements of $O_i$. Let $\beta_i:=(uv_i)$, and let $\gamma_i:=(\pi^{-1}(\beta_i))^{-1}\gamma\pi^{-1}(\beta_i)$. Then $\gamma_i$ acts on the first coordinate of the fiber $F_i\times \{u\}$ as $\alpha_u^{-1}\alpha_v$, and it fixes every element outside $F_i\times \{u,v_i\}$. Therefore, the permutations $\gamma_i$ converge to a permutation $\gamma'$ which acts on the first coordinate of the fiber $F_i\times \{u\}$ as $\alpha_u^{-1}\alpha_v$, and fixes every element outside $F_i\times \{u\}$. By our assumption $G$ is closed, so $\gamma'\in G$. By definition this implies that $\alpha_u^{-1}\alpha_v\in N_i(G)$.
\end{proof}

\begin{proposition}\label{desc_nh}
	Let $\fc$ be a covering reduct of $\fa$ and $G := \aut(\fc) \cap N$.
	Then $$G=N(H(G),N_1(G),\dots,N_k(G)).$$
\end{proposition}

\begin{proof}
	By Remark~\ref{hn_equal} we know that $N_i(G)\leq H_i(G)=H(G)_{F\setminus F_i}|_{F_i}$, and if $|O_i|=1$ then $N_i(G)=H_i(G)$. Lemma~\ref{normal_hn} implies that 
	$H(G)|_{F_i}$ normalises $N_i(G)$. Therefore, the group $G=N(H(G),N_1(G),\dots,N_k(G))$ is well-defined.
	
		We first show that $G\leq N(H(G),N_1(G),\dots,N_k(G))$. Let $\alpha \in G$.
	Then the definition of the group $H(G)$ implies that $\alpha$ satisfies the first item in the definition of $N(H(G),N_1(G),\dots,N_k(G))$. 
	By Lemma \ref{generating_n}, $\alpha$ also satisfies the second item of this definition, and
	hence $\alpha \in N(H(G),N_1(G),\dots,N_k(G))$. 

	Now let $\alpha\in N(H(G),N_1(G),\dots,N_k(G))$ be arbitrary. Let $u_i\in O_i$ be arbitrary elements. By Lemma~\ref{nh_going_back} we have $H(N(H(G),N_1(G),\dots,N_k(G))=H(G)$. This implies that 
	there exists an $\alpha' \in G$ such that for every $i \in \{1,\dots,k\}$
	the actions of $\alpha$ and $\alpha'$ agree on $F_i\times \{u_i\}$. 
	For $v\in O_i$ let $\alpha_v$ and $\alpha_v'$ denote the action of $\alpha$ and $\alpha'$, respectively, on the fiber $F_i\times \{v\}$. We claim that for all $v\in O_i$ it holds that $\alpha_v(\alpha_v')^{-1}\in N_i$. 
	By Lemma~\ref{nh_going_back} we have $N_i(N(H(G),N_1(G),\dots,N_k(G))=N_i(G)$, and hence by Lemma~\ref{generating_n} it follows that $\alpha_v\alpha_{u_i}^{-1}\in N_i(G)$, and $\alpha_v'(\alpha_{u_i}')^{-1}=\alpha_v'\alpha_{u_i}^{-1}\in N_i(G)$, and hence $$\alpha_v(\alpha_v')^{-1}=\alpha_v\alpha_{u_i}^{-1}\alpha_{u_i}(\alpha_v')^{-1}=\alpha_v\alpha_{u_i}^{-1}(\alpha_v'\alpha_{u_i}^{-1})^{-1}\in N_i(G).$$


	This implies that $$\alpha(\alpha')^{-1}\in N(\prod_{i=1}^k{N_i}(G),N_1(G),\dots,N_k(G))=N^*(N_1(G),\dots,N_k(G))\subseteq G.$$ Therefore $\alpha=(\alpha(\alpha')^{-1})\alpha'\in G$. 
\end{proof}

\begin{remark}
Proposition~\ref{desc_nh} generalizes 
Theorem 3.1 in~\cite{FiniteCovers} from $\sym(\set)$
to arbitrary automorphism groups of structures in $\unary^*$.
\end{remark}

\begin{corollary}\label{reduct_trivial_finite}
$\fa$ has finitely many covering reducts with respect to $\pi$.
\end{corollary}
\begin{proof}
	By Lemma \ref{semidir2} and item (1) of Proposition~\ref{semidir} it is enough to show that $N$ has finitely many closed subgroups which are normalized by $\aut(\fa)$. By Proposition~\ref{desc_nh} and Remark~\ref{hn_equal} these groups can be characterised by a subgroup $H$ of $\prod_{i=1}^k{\sym(F_i)}$ and a system of subgroups $N_i \leq H_i$.
	Then the statement of the corollary follows from the fact that there are finitely many choices for these group.
\end{proof}

%% file: FRU.tex

\section{Finite Coverings of Reducts of Unary Structures}
In this section we show that every structure in 
$\fcover(\reduct(\unary))$ is a quasi-covering reduct (introduced in Definition~\ref{quasi_cover}) of a strongly trivial covering of some structure in $\unary^*$ (Proposition~\ref{quasi_trivial}), and that there are only finitely many of such reducts for each structure in $R(\unary)$ (Theorem~\ref{quasi_cover2}). 
Moreover, we observe that $\fcover(\reduct(\unary)) \subseteq \reduct(\fcover(\unary))$ (Remark~\ref{fru_rffu}).

\subsection{The Ramsey property and canonical functions}
  Let $\fa,\fb$ be structures.
  A function $f \colon A \to B$ is called \emph{canonical from $\fa$ to $\fb$} if for every $t \in A^n$ and 
$\alpha \in \aut(\fa)$ there exists $\beta \in \aut(\fb)$ such that $f(\alpha(t)) = \beta(f(t))$. 
Hence, a canonical function $f$ induces for every $n$ a function 
from the orbits of $n$-tuples of $\aut(\fa)$
to the orbits of $n$-tuples of $\aut(\fb)$;
these functions will be called the \emph{behavior} of $f$. 
Canonical functions as a tool 
  to classify reducts of homogeneous structures with finite relational signature have been introduced in~\cite{BP-reductsRamsey} and used in~\cite{Poset-Reducts,42,agarwal,AgarwalKompatscher,BodJonsPham,BBPP18}. The existence of certain canonical functions in the automorphism group of a 
  structure $\fa$ 
  is typically shown using Ramsey properties of $\fa$. We will not introduce Ramsey structures here; all that is needed is the well-known fact that $({\mathbb Q};<)$ is Ramsey, and the following result from~\cite{BP-reductsRamsey}.
  A structure is called \emph{ordered} if the signature
contains a binary relation symbol that denotes a (total) linear ordering of the domain.

\begin{lemma}[see~\cite{BodPin-CanonicalFunctions}]\label{gen_can}
	Let $\fd$ be an ordered homogeneous Ramsey structure with finite relational signature and let $f \colon D \to D$ be a function.
	Then there exists a function $$g \in \overline{\{\alpha \circ f \circ \beta \mid \alpha,\beta \in \aut(\fd)\}}$$
	which is canonical as a function from $\fd$ to $\fd$. 
\end{lemma}

The following is an easy consequence of the definitions.

\begin{lemma}\label{gen_behav}
	Let $\fa$ be a homogeneous structure
	with finite relational signature and let $\fb$ be a first-order reduct of $\fa$. 
	If $f$ and $g$ are canonical functions from 
	$\fa$ to $\fb$ with the same behaviour then 
	$\overline{\aut(\fb) \cup \{f\}} = \overline{\aut(\fb) \cup \{g\}}$.
\end{lemma}

The next lemma follows from the observation that if $\fa$ is homogeneous with a relational signature of maximal arity $k$ then the behaviour of a canonical function $f$ from 
$\fa$ to $\fb$ is fully determined by the function 
induced by $f$ on the orbits of $k$-tuples
(see~\cite{BPT-decidability-of-definability}, in particular the comments at the end of Section 4.1). 

\begin{lemma}\label{fin_many_behav}
	Let $\fa$ be a homogeneous structure
	with finite relational signature and let $\fb$ be 
	$\omega$-categorical. 
	Then there are finitely many behaviours of canonical functions from $\fa$ to $\fb$.
\end{lemma}

The structures in $\fcover(\unary^*)$ have homogeneous expansions with finite relational signature which we describe next. 

\begin{lemma}\label{ramsey}
	Let $\fb\in \unary^*$ and let $\pi \colon\fa\rightarrow \fb$ be a strongly trivial finite cover. Then
\begin{enumerate}
\item $\fa$ is interdefinable with a homogeneous structure $\fc$ with finite relational signature, and 
\item $\fa$ is a first-order reduct of an ordered homogeneous Ramsey structure $\fd$ with finite relational signature.
\end{enumerate}
\end{lemma}
\begin{proof}
Let $O_1,\dots,O_k$ be the orbits of $\fb$.
	Following Remark~\ref{triv_cov_unary}, we can assume that $\fa=\bigsqcup_{i=1}^k{(F_i \times O_i)}$ for some finite sets $F_i$, and that $\aut(\fa)$ consists of all permutations which preserve the first coordinate and stabilise the sets $O_i$ on the second coordinate. For each $i \leq k$ and $s\in F_{i}$ we define the unary relation 
	$U_{i,s} :=\{(s,u) \mid u\in O_i\}$. 
	Let $\fc$ be the relational structure with domain $A$ and the relations $U_{i,s}$ and $\sim_{\pi}$. Then $\aut(\fc)=\aut(\fa)$. Hence, $\fa$ and $\fc$ are interdefinable. It is easy to see that $\fc$ is homogeneous. This proves (1). 

	To prove item (2) we define an ordering $<$ on $\fb$ as follows. For each infinite orbit $O_i$
	let us fix an ordering $<_i$ on $O_i$ which is isomorphic to $(\mathbb{Q};<)$. Let us also fix an ordering of $\prec_i$ on $F_{i}$ for all $i$. Then $<$ is defined as follows 
\begin{itemize}
\item If $\pi(x)\in O_i, \pi(y)\in O_j$ and $i<j$, then $x<y$,
\item If $\pi(x),\pi(y)\in O_i$ and $\pi(x)<_i\pi(y)$, then $x<y$,
\item If $\pi(x)=\pi(y) \in O_i$ and $x'$ and $y'$ are the projections of $x$ and $y$ to the first component, then $x<y$ iff $x'
\prec_i y'$. 
\end{itemize}
	
	To show that the expansion $\fd$ of $\fc$ 
	by the ordering $<$ has the Ramsey property,
	we use the fact that if a structure is the disjoint union of substructures induced by definable subsets, and the substructures are Ramsey, 
	then the structure itself is Ramsey (see~\cite{BodirskyRamsey}). 
		For each $i \leq k$, let $\fd_i$ be the substructure of $\fd$ induced by 
	$\pi^{-1}(O_i)$. Note that $\pi^{-1}(O_i)  = \bigcup_{s \in F_i} U_{i,s}$ and hence is definable in $\fd$. If $O_i$ is infinite then $\aut(\fd_i)$ is topologically isomorphic to $\aut({\mathbb Q};<)$. The property of a structure of being Ramsey is a property of the automorphism group of the structure, viewed as a topological group (again, see~\cite{BodirskyRamsey}).
	It follows that $\fd$ is Ramsey.
\end{proof}

\subsection{Reducts of Finite Covers of Reducts of Unary Structures}
Let $\fb \in \reduct(\unary)$ 
and let $\pi \colon \fa \to \fb$ be a finite covering map. 
In this section we study the closed supergroups $G$
of $\aut(\fa)$ that preserve $\sim_\pi$ 
such that $\mu_\pi(G)$ also preserves the congruence $\nabla(\fb)$.

\begin{lemma}\label{minimal_fiber_pres}
	Let $\fb\in \unary^*$ and let $\pi \colon \fa \rightarrow \fb$ be a strongly trivial finite covering map. Let $O_1,\dots,O_k$ be the orbits of $\fb$ and let $\fd$ be the ordered homogeneous finite signature Ramsey expansion of $\fa$ from Lemma~\ref{ramsey}.
	Suppose that $f \in \sym(A)$ 
	preserves $\sim_{\pi}$ and that $\mu_\pi(f)$ preserves the partition $P := \{O_1,\dots,O_k\}$. 
	Then the monoid $M:=\overline{\langle \aut(\fa),f\rangle}$ contains a \emph{surjective} map $h$ which is canonical from $\fd$ to $\fa$ and such that $\mu_\pi(h)$ has the same action on $P$ as $\mu_\pi(f)$.  
\end{lemma}
\begin{proof}
	By Lemma \ref{gen_can} we obtain that there exists a function $$g\in \overline{\{\alpha f \beta \mid \alpha,\beta \in \aut(\fd)\}} \subseteq M$$ which is canonical from $\fd$
	to $\fd$. 
Since $g\in M$ it follows that $g$ preserves $\sim_{\pi}$.
	But note that $g$ is not necessarily surjective. For $m \in M$ 
	define $\mu_\pi(m)$ by $x \mapsto \pi(m(\pi^{-1}(x)))$ as in the case of bijective functions. 
	Since every automorphism of $\fd$ 
	preserves the orbits $O_1,\dots,O_k$, it follows that if $\mu_\pi(f)(O_i)=O_j$ then $\mu_\pi(g)(O_i)\subseteq O_j$. 

	If $\mu_\pi(f)(O_i)=O_j$, then $|F_{O_i}|=|F_{O_j}|$ since $f$ preserves $\sim_{\pi}$  
	and is surjective. 
	Therefore, for every $g' \in M$ and
	$u \in B$ the restriction of $g'$ 
	to $U := \pi^{-1}(u)$ 
	is a bijection between $U$ and 
	$\pi^{-1}(\mu_\pi(g')(u))$. In particular, this holds for $g \in M$. If $O_i=\{u_i\}$ then $g$ defines a bijection between $\pi^{-1}(u_i)$ and $\pi^{-1}(\mu_\pi(f)(u_i))$. If $O_i$ is infinite then $\pi^{-1}(\mu_\pi(g)(O_i))$ is a union of infinitely many classes of $\sim_{\pi}$. 
	Moreover, $O_i$ is infinite if and only if $\mu_\pi(g)(O_i)$ is infinite. 
	Let $e \colon A \to A$ be defined as
	$(s,u) \mapsto (s,\mu_{\pi}(g)(u))$.
	Let $\fc$ be the homogeneous structure 
	from Lemma~\ref{ramsey} which has the property that $\aut(\fc) = \aut(\fa)$. 
	Then $e$ is an isomorphism between $\fc$ and the substructure of $\fc$ induced by $g(C)$.
	Since $\fc$ is homogeneous it follows that $e\in \overline{\aut(\fc)} = \overline{\aut(\fa)}$ and so there is a sequence $e_1,e_2,\ldots\in \aut(\fa)$ which converges to $e$. Then $h_i:= e_i^{-1} g\in M$ converges to $h:=e^{-1}\circ g$  and thus $h \in M$. We claim that the mapping $h$ satisfies the conditions of the lemma. By definition $h(A)=e^{-1}(g(A))=A$, that is, $h$ is surjective. Since $e^{-1}$ preserves the relations of $\fc$ it follows that the mapping $h$ is canonical from $\fd$
	to $\fc$ (and therefore also from $\fd$ to $\fa$). For the same reason 
	 $\mu_\pi(e^{-1})$ 
	preserves all orbits $O_i$. This implies that $\mu_\pi(h)$ and $\mu_\pi(g)$ and therefore also $\mu_\pi(f)$ have the same action on $\{O_1,\dots,O_k\}$. 
\end{proof}

\begin{lemma}\label{minimal_fiber_pres2}
	Let $\fb\in \unary^*$ and let $\pi \colon \fa \rightarrow \fb$ be a finite covering map. Let $O_1,\dots,O_k$ be the orbits of $\fb$. 
	Then $\sym(A)$ has finitely many closed subgroups $G$ such that
	\begin{itemize}
	\item $\aut(\fa) \subseteq G$,  
	\item $G$ preserves $\sim_{\pi}$, and 
	\item $\mu_{\pi}(G)$ preserves the partition 
	$\{O_1,\dots,O_k\}$ of $B$. 
	\end{itemize}
\end{lemma}
\begin{proof}
	By Proposition \ref{reduct_trivial2} we know that $\fa$ is a covering reduct of some strongly trivial covering $\fc$ of $\fb$ (with respect to $\pi$). Then let $\fd$ be the ordered homogeneous finite-signature Ramsey expansion of $\fc$ from Lemma~\ref{ramsey}.
	Let $G$ be a closed subgroup of $\sym(A)$ as in the formulation of the lemma.
	 Then $G$ acts on the set $\{O_1,\dots,O_k\}$. 
	 Let $K$ be the kernel of this action. 
	 Then $K$ is closed and the index of $K$ in $G$ is finite. Also, $K\subseteq S$ where $S$ is the group as in Definition~\ref{group_ns}. 
	Therefore, $K$ is the automorphism group of a covering reduct of $\fa$ (Proposition~\ref{semidir}). 
	Then by Theorem~\ref{reduct_trivial_finite1} there are finitely many possible choices for the group $K$.

	By Lemma \ref{minimal_fiber_pres}, for each $f \in G$ there
	exists a surjective map $h \in \overline{\langle K, f\rangle}$ which is canonical from $\fd$ to $\fa$ such that $f$ and $h$ induce the same permutation $\sigma$ of $\{O_1,\dots,O_k\}$. 
	We claim that $K \cup \{f\}$ and $K \cup \{h\}$ generate the same group. The image of the action of $K\cup \{f\}$ and $K \cup \{h\}$ on $\{O_1,\dots,O_k\}$ is $\langle \sigma \rangle$, and the kernel of these actions is again $K$. 
	Therefore,
	\begin{align}
	[\langle K\cup \{f\}\rangle:K]=[\langle K\cup \{h\}\rangle:K]=l
	\label{eq:index}
	\end{align} where $l$ is the order of the permutation $\sigma$. In particular, the groups $\langle K\cup \{f\}\rangle$ and $\langle K\cup \{h\}\rangle$ are closed. Hence, $h\in K\cup \{f\}$ and thus $\langle K\cup \{f\}\rangle\leq \langle K\cup \{h\}\rangle$. Then by using Equality~$(\ref{eq:index})$ again
	it follows that $\langle K\cup \{f\}\rangle=\langle K\cup \{h\}\rangle$.
 
	Since $[G:K]$ is finite each group $G$ is generated by finitely many (at most $k!$) elements over $K$. By the previous paragraph we can assume that each of these generators are canonical from $\fd$ to $\fa$. There are finitely many possible behaviours of canonical functions from $\fd$ to $\fa$
	(Lemma \ref{fin_many_behav}).
	If two functions have the same behaviour 
they generate the same group over $\aut(\fa)$ 
	(Lemma \ref{gen_behav}). 
	This implies that there are finitely many choices for the group $G$.
\end{proof}

\begin{theorem}\label{fru}
	Let $\fb\in \reduct(\unary)$ and let $\pi \colon \fa\rightarrow \fb$ be a finite covering map. Then $\aut(\fa)$ has finitely many closed supergroups
	$G$ such that $G$ preserves $\sim_{\pi}$  
	and $\mu_\pi(G)$ preserves $\nabla(\fb)$.
\end{theorem}
\begin{proof}
	Let $O_1,\dots,O_k$ be the classes of $\nabla(\fb)$. Then by Lemma \ref{unary_reduct} it follows that $\prod_{i=1}^k{\sym(O_i)} \subseteq \aut(\fb)$. Let $\fb'$ be a structure with $\aut(\fb')=\prod_{i=1}^k{\sym(O_i)}$. Then $\fb'\in \unary^*$ by Lemma~\ref{sing_or_inf} and $\nabla(\fb)=\nabla(\fb')$. The group $\aut(\fa)$ acts naturally on the set 
	$\{O_1,\dots,O_k\}$. Let $K$ be the kernel of this action, and let $\fa'$ be a structure so that $\aut(\fa')=K$. The action of $\aut(\fa')$ on $B$ equals $\aut(\fb')$.
	Therefore, $\pi \colon  \fa'\rightarrow \fb'$ is a finite cover. Then the statement of the theorem follows from Lemma \ref{minimal_fiber_pres2} and from the fact that the orbits of $\fb'$ are exactly the classes of the congruence $\nabla(\fb)$.
\end{proof}

\begin{proposition}\label{fru_rffu}
$\fcover(\reduct(\unary))\subseteq\reductfin(\fcover(\unary^*))$
\end{proposition}
	\begin{proof}
	Let $\fa\in \fcover(\reduct(\unary))$. Following the notation of the proof of  Theorem~\ref{fru} we have that 
	$[\aut(\fa):\aut(\fa')]=[\aut(\fa):K]$ is finite since $K$ is defined as the kernel of the action of $\aut(\fa)$ on the set $\{O_1,\dots,O_k\}$. As we saw in the proof of Theorem~\ref{fru} we have $\fa'\in \fcover(\unary^*)$. Hence, $\fa\in \reductfin(\fcover(\unary^*))$. 
	\end{proof}
	
	Later we will see (in Theorem~\ref{main_kexpp}) that in fact $\fcover(\reduct(\unary))=\reductfin(\fcover(\unary^*))$. The following definition of \emph{quasi-covering reducts} is needed for a model theoretic reformulation of Theorem~\ref{fru}, which is given in Theorem~\ref{quasi_cover2} below.
	
\begin{definition}\label{quasi_cover}
	Let $\fb$ be $\omega$-categorical and 	let $\pi \colon \fa \rightarrow \fb$ be a finite cover. A first-order reduct $\fc$ of $\fa$ is called a \emph{quasi-covering reduct of $\fa$ with respect to $\pi$
	 if 
	$\aut(\fc)$ 
	preserves $\sim_{\pi}$ and 
	$\mu_\pi(\aut(\fc)) \subseteq \sym(B)$
	preserves $\nabla(\fb)$.}
\end{definition}

\begin{theorem}\label{quasi_cover2}
	Let $\fb\in \reduct(\unary)$ and let $\pi \colon \fa\rightarrow \fb$ be a finite cover. Then $\fa$ has finitely many quasi-covering reducts with respect to $\pi$. 
\end{theorem}
\begin{proof}
	The statement follows immediately from Theorem \ref{fru} and Definition \ref{quasi_cover}.
\end{proof}

\begin{proposition}\label{quasi_trivial}
	Let $\fb\in \reduct(\unary)$ and let $\pi \colon \fa\rightarrow \fb$ be a finite cover. Then $\fa$ is a quasi-covering reduct of a strongly trivial covering of some structure in $\unary^*$.
\end{proposition}

\begin{proof}
	Let us define the structures $\fa'$ and $\fb'$ as in the proof of Theorem \ref{fru}. Then $\fa$ is a finite quasi-covering reduct of $\fa'$ with respect to the covering map $\pi$. The map $\pi \colon \fa' \rightarrow \fb'$ is a finite covering, and $\fb'\in \unary^*$. By Proposition \ref{reduct_trivial2}, $\fa'$ is a covering reduct of some strongly trivial covering of $\fb'$. Therefore, $\fa$ is a quasi-covering reduct of a strongly trivial covering of $\fb'\in \unary^*$.
\end{proof}

%% file: Kexp.tex

\begin{section}{Structures with Small Orbit Growth}
In this section we show that $\kexpp=(\fcover\circ \reduct)(\unary)$.
We start in Section~\ref{sect:growth} with some observations from enumerative combinatorics that we need to obtain 
information about $\nabla(\fa)$ if $\fa \in \kexpp$. 
In Section~\ref{sect:stab} we discuss the effect of stabilising a group from $\gexpp$ 
at finitely many constants. Section~\ref{sect:primitive} treats the case that
$\fa$ is primitive; here we rely on work of Macpherson~\cite{Macpherson-Orbits}.
We then focus on permutation groups $G$ in $\gexpp$ where the congruence $\Delta(G)$ is trivial, i.e., $\Delta(G)$ is the identity congruence (Section~\ref{sect:delta-trivial}); the general case is treated in Section~\ref{sect:general}. 
In Section~\ref{sect:thomas} we use these results to prove Thomas' conjecture for all structures in the class $\kexpp$.


\input partition-growth.tex

\subsection{The number of $\nabla$-classes in point stabilisers}
\label{sect:stab}
In this section we examine the possible growth of the number of $\nabla$-classes in stabilisers of finite sets.

\begin{lemma}\label{stabil_finite}
	Let $G \in \gexpp$ be a permutation group on a countably infinite set $X$, that is, $\oi_n(G)\leq c_1n^{dn}$ for some $c_1,d$ with $d<1$.
	Let $F\subset X$ be finite. Then for every $\varepsilon>0$
	\begin{itemize}
	\item  there exists a constant $c_2$ such that $$\oi_n(G_{F}|_{X \setminus F}) < c_2n^{(d+\varepsilon) n}$$
	\item there exists a constant $c_3$ such that $$\oi_n(G_{F}) < c_3n^{(d+\varepsilon) n}.$$
		\end{itemize}
	In particular, $G_F\in \gexpp$ and $G_F|_{X \setminus F} \in \gexpp$. 
\end{lemma}
\begin{proof}
Let $\varepsilon>0$. 
	The orbits of injective $n$-tuples of $G_F$ can be embedded into the orbits of injective $(n+|F|)$-tuples of $G$ by mapping the orbit of a tuple $t$ into the orbit of $(t,t')$ where $t'$ is any $|F|$-tuple such that $(t,t')$ has pairwise distinct entries and all elements of $F$ appear in $(t,t')$. 
	Hence, $$\oi_n(G_F) \leq \oi_{n+|F|}(G) \leq c_1(n+|F|)^{d(n+|F|)} \leq c_2 n^{(d+\varepsilon) n}$$ 
	for an appropriate constant $c_2$. Choosing 
	$\varepsilon>0$ such that $d+\varepsilon<1$
	shows that $G_F \in \gexpp$. The statements for $G_{F}|_{X \setminus F}$ can be shown analogously.
\end{proof}


\begin{definition}
	Let $G$ be an oligomorphic permutation group on a countably infinite set $X$. For every finite set $F\subset X$ let $m_G(F)$ be the number of $\nabla(G_F)$-classes. For $n\in \mathbb{N}$ let $m_G(n):=\max(\{m_G(F) \mid F\subset X, |F|=n\})$.
\end{definition}

\begin{remark}
	If $F_1,F_2 \subset X$ are contained in the same orbit of $n$-subsets of $G$, then $m_G(F_1)=m_G(F_2)$. Hence, the set $\{m_G(F) \mid F\subset X, |F|=n\}$ is finite, and so the maximum of this set always exists.
\end{remark}

\begin{lemma}\label{not_too_many_classes}
	Let $G$ be a permutation group. 
	Suppose that $\oi_n(G)\leq c_1n^{dn}$ for some $c_1$ and $d<1$. 
	Then for every $\varepsilon>0$ we have $m_G(n)\leq c_2n^{d+\varepsilon}$ 
	for some constant $c_2$. 
\end{lemma}
\begin{proof}
	Suppose that $G$ is a permutation group on $X$ and let $F\subset X$ be of size $n$. 
	Suppose that $X_1,X_2,\dots,X_l$ are the infinite classes of the congruence $\nabla(G_F)$, and arbitrarily choose 
	$x_i \in X_i$ for $i \in \{1,\dots,l\}$. 
	Then for each function $f \colon \{1,\dots,k\}\rightarrow \{1,\dots,l\}$ there are pairwise distinct elements $y_1,\dots,y_k$ so that $y_j\in X_{f(j)}\setminus \{x_{f(j)}\}$. Let $t^f:=(x_1,\dots,x_l,y_1,\dots,y_k)$. Then the tuples $t^f$ are injective and lie in pairwise different orbits of $G_{F}|_{X \setminus F}$. Thus $l^n\leq \oi_{n}(G_{F}|_{X \setminus F})\leq c_2n^{(d+\varepsilon) n}$ for some $c_2$ by Lemma \ref{stabil_finite}. Thus $l\leq c_2n^{d+\varepsilon}$, and therefore $m_G(F)\leq c_2n^{d+\varepsilon}$. 
\end{proof}



\begin{lemma}\label{small_classes_2}
	Let $G$ be a permutation group on a countably infinite set $X$ and suppose that $\oi_n(G)\leq cn^{dn}$ for some $c$ and $d<1$. Let $F\subset X$ be finite, let $R$ be a congruence of $G_F$, and let $\frac{k-1}{k}>d$.  Then $R$ has finitely many classes of size at least $k$.
\end{lemma}

\begin{proof}
	By Lemma~\ref{stabil_finite} we can assume that $F=\emptyset$. Suppose to the contrary that $R$ has infinitely many classes of size at least $k$. Let $n$ be arbitrary and let $\mathcal{P}_n^k$ be the set of partitions $P=\{S_1,\dots,S_l\}$ of $\{1,2,\dots,n\}$ such that $|S_i|\leq k$ for all $i=1,\dots,l$. For each $P \in \mathcal{P}_n^k$ we can choose pairwise distinct elements $x_1^P,\dots,x_n^P \in X$ such that $x_i^PRx_j^P$ if and only $x_i^P$ and $x_j^P$ are in the same set in the partition $P$. Then the $n$-tuples $(x_1^P,\dots,x_n^P)$ for $P\in \mathcal{P}_n^k$ are injective and lie in pairwise different orbits of $G$. Therefore $\oi_n(G)\geq |\mathcal{P}_n^k|$. Let us choose $\varepsilon>0$ such that $\frac{k-1}{k}-\varepsilon>d$. Then by Lemma \ref{counting} it follows that $$\oi_n(G)\geq|\mathcal{P}_n^k| = p_k(n) \geq n^{(\frac{k-1}{k}-\varepsilon)n}>cn^{dn}$$ for $n$ large enough. This contradicts our assumption.
\end{proof}

\begin{corollary}\label{small_classes_3}
	Let $G$ be a permutation group on a countably infinite set $X$ and suppose that $\oi_n(G)\leq cn^{dn}$ for some $c$ and $d<1$. Let $F\subset X$ be finite and let $R$ be a congruence of $G_F$. Then $R$ has finitely many infinite classes.
\end{corollary}

\begin{proof}
	Follows directly from Lemma \ref{small_classes_2}.
\end{proof}

\begin{definition}
	Let $(\gexpp)^k$, for $k\in \mathbb{N}$,  denote the class of those groups $G\in \gexpp$ for which the following holds.
\begin{itemize}
\item[$(*^k)$] For every finite $F \subset X$, every congruence of $G_F$ has at most finitely many equivalence classes of size at least $k$.
\end{itemize}
	Let $(\kexpp)^k$ denote the classes of those structures in $\kexpp$ whose automorphism group is in $(\gexpp)^k$.
\end{definition}

	Using the definition above,  Lemma \ref{small_classes_2} immediately implies the following.

\begin{corollary}\label{union_delta_size}
	$\gexpp=\bigcup_{k=1}^\infty(\gexpp)^k$, and $\kexpp=\bigcup_{k=1}^\infty(\kexpp)^k$.
\end{corollary}

\subsection{The primitive case}
	\label{sect:primitive}
	We use the following theorem of Dugald Macpherson~\cite{Macpherson-Orbits}.

\begin{theorem}[Theorem 1.2 in~\cite{Macpherson-Orbits}]\label{dugald_primitive}
	Let $G$ be a permutation group 
	on a countably infinite set $X$ which is  primitive but not highly transitive. Then there is a polynomial $p$ such that $\oi_n(G)\geq \frac{n!}{p(n)}$.
\end{theorem}

	Theorem \ref{dugald_primitive} immediately implies the following.

\begin{lemma}\label{primitive}
	Let $G \in \gexpp$ be primitive. Then $G$ is highly transitive.
\end{lemma}

\begin{proof}
	Stirling's formula implies that 
	for all $c$ and $d<1$ and every polynomial $p$, if $n$ is large enough then $\frac{n!}{p(n)} > cn^{dn}$. Hence, the lemma follows from Theorem \ref{dugald_primitive}.
\end{proof}

\subsection{The case when $\Delta(G)$ is trivial}
\label{sect:delta-trivial}
The result of Macpherson (Theorem~\ref{dugald_primitive}) is used via the following lemma.

\begin{lemma}\label{nice_case}
	Let $G \in \gexpp$ be such that $\Delta(G)$ is trivial and such that $G$ stabilises each class of $\nabla(G)$ setwise. Then $G$ acts highly transitively on each of its orbits. 
\end{lemma}

\begin{proof}
	Let $O_1,\dots,O_m$ be the orbits of $G$. Then $O_1,\dots,O_m$ are also the classes of $\nabla(G)$. We claim that the action of $G$ on $O_i$ is primitive for each $i \in \{1,\dots,m\}$; 
	this suffices, because then the statement of the lemma follows from Lemma~\ref{primitive}.

	Let $R_i$ be a congruence of $G|_{O_i}$. Since $G$ acts transitively on $O_i$ it follows that every class of $R_i$ has the same size. If this size is finite, then let us consider the congruence $R_i^*:=R_i\cup \{(x,x) \mid x\in X\setminus O_i\}$. Then every class of $R_i^*$ is finite and thus $R_i^*$ must be finer than $\Delta(G)$. Since $\Delta(G)$ is the identity congruence, 
	 $R_i^*$ and $R_i$ are the identity congruence, too. 

	Now assume that every class of $R_i$ is infinite. Then by Corollary \ref{small_classes_3} $R_i$ has finitely many classes. Let $C_1,C_2,\dots,C_l$ be these classes. Then $\{O_1,\dots,O_{i-1},O_{i+1},\dots,O_m,C_{1},C_{2},\dots,C_{l}\}$ is an invariant partition. Since $\nabla(G)$ is the finest congruence with finitely many classes, it follows that $l=1$ and thus $R_i$ is again the identity congruence. Therefore, $G|_{O_i}$ is primitive for all $i$.
\end{proof}

	Under the conditions of Lemma~\ref{nice_case} we will show that, in fact, if $G$ is closed, then $G=\prod_{i=1}^m{\sym(O_i)}$ where $O_1,\dots,O_m$ are the orbits of $G$, that is, $G$ is the automorphism group of a unary structure (Lemma \ref{extension_2}). The following lemma is well-known (see e.g.\ Proposition~1.4(2) in~\cite{MoonStalder} or Corollary~2.2 in~\cite{cameron1981normal}).

\begin{lemma}\label{highly_normal}
	Every normal subgroup of a highly transitive permutation group acting on an infinite set is either highly transitive or trivial.
\end{lemma}

\begin{lemma}\label{two_parts}
	Let $G$ be a closed permutation group on a countably infinite set $X$. Let $T$ be an infinite orbit of $G$ such that $G|_T$ is highly transitive and let $S:=X \setminus T$.
	Then one of the following holds.
\begin{enumerate}
\item $\{\id_S\}\times\sym(T) \subseteq G$,
\item 
There exists a surjective homomorphism 
$e \colon G|_S \to G|_T$ such that a permutation $\gamma$ of $X$ is in $G$ if and only if  
	there exists a permutation 
$\alpha \in G|_S$ so that $\gamma|_S=\alpha$ and 
$\gamma|_T=e(\alpha)$.
\end{enumerate}
\end{lemma}
\begin{proof}
	If $\alpha\in \sym(S)$ and $\beta\in \sym(T)$  then we use the notation $(\alpha,\beta)$ for the unique permutation 
	$\gamma \in \sym(X)$ 
	whose restriction to $S$ equals $\alpha$
	and whose restriction to $T$ equals $\beta$. 

\medskip 
	\emph{Case 1. For every $\alpha \in G|_S$ there is a unique $e(\alpha) \in G|_T$ such that $(\alpha,e(\alpha)) \in G$.} 
	
	\medskip 
	It is easy to see that in this case $e$ is a surjective homomorphism from $G|_S$ to $G|_T$, therefore Condition (2) holds.

\medskip 
	\emph{Case 2. For some $\alpha\in G|_S$ there exist at least two distinct permutations $\beta_1,\beta_2 \in G|_T$ such that $\gamma_1:=(\alpha,\beta_1) \in G$ and $\gamma_2:=(\alpha,\beta_2) \in G$.} 
	
	\medskip 
	Let 
	$K:=\{\beta\in \sym(T) \mid (\id_S,\beta)\in G\}$.  
	Then 
\begin{itemize}
\item \emph{$K$ is nontrivial} since $\beta_1\beta_2^{-1}\in K$.
\item \emph{$K$ is closed in $\sym(T)$}, and
\item \emph{$K$ is a normal subgroup of $\sym(T)$}.
\end{itemize}
 To prove normality, let $(\id, \beta)\in G$, and let $\delta \in \sym(T)$ be arbitrary. 
  Since $G|_T$ is dense in $\sym(T)$ there is a sequence $\delta_1,\delta_2,\ldots$ of elements of $G|_T$ which converges to $\delta \in \sym(T)$. By the definition of $G|_T$ we know that there exist elements $\alpha_i\in \sym(S)$ such that $\eta_i:=(\alpha_i,\delta_i)\in G$  for every $i \in {\mathbb N}$. Then $G\ni \eta_i(\id,\beta)\eta_i^{-1}=(\id, \delta_i \beta \delta_i^{-1})$.
   Therefore $\lim_i(\id, \delta_i \beta \delta_i^{-1})=(\id, \delta \beta \delta^{-1})\in G$ since $G$ is closed. By definition, this implies that $\delta \beta\delta^{-1}\in K$ which shows that $K$ is indeed a normal subgroup. By Lemma~\ref{highly_normal}, 
$K=\sym(T)$. Thus, $\{\id_S\}\times \sym(T)\subseteq G$, i.e., Condition (1) holds.
\end{proof}

\begin{lemma}\label{closedness}
Let $G$ be a closed oligomorphic permutation group on a countably infinite set $X$.
Let $O_1,\dots,O_m$ be the orbits of $G$. Let us suppose that each $O_i$ is infinite, and that $G$ acts highly transitively on each $O_i$. 
Let $l \in \{1,\dots,m\}$ and let $S := \bigcup_{i=1}^l O_i$ be such that $\acl_G(S)=X$.
Then $H:=G|_S$ is closed in $\sym(S)$.
\end{lemma}

\begin{proof}
	We first show the statement for $l=m-1$.  
	Let $T:=O_m$ (so that we have the same notation as in Lemma \ref{two_parts}). First, let us assume that Condition (1) of Lemma \ref{two_parts} holds. Let $(\alpha_j)_{j \in {\mathbb N}}$ be a sequence that converges in $G|_S$ to $\alpha \in \sym(S)$. Let $\beta_j\in \{\id_S\}\times\sym(T) \subseteq G$ be such that $\alpha_j|_T=\beta_j|_T$. Let 
	$\alpha_j':=\alpha_j\beta_j^{-1}\in G$. 
	Then $\alpha_j'\rightarrow (\alpha,\id (T))\in G$ since $G$ is closed. In particular $\alpha\in G|_S$ and $H$ is closed.
	
Otherwise, if Condition (1) of Lemma \ref{two_parts} does not hold, then by Lemma \ref{two_parts} we can assume that item (2) of Lemma \ref{two_parts} holds. 
Let 
$e \colon G|_S \to G|_T$ be as in item (2) of Lemma \ref{two_parts}.
	If $F\subset S$ is finite and $\alpha\in G|_S$  then \begin{align*}
	\acl(\alpha(F))\cap T & =\acl((\alpha,e(\alpha))(F))\cap T\\
	& = (\alpha,e(\alpha))(\acl(F))\cap T \\
	& =(\alpha,e(\alpha))(\acl(F)\cap T) 
	=e(\alpha)(\acl(F)\cap T).
	\end{align*}
	 By assumption $\acl(F)\cap T$ is nonempty for some $F$ (and it is always finite). Let $k:=|{\acl}(F)\cap T|$. Then by our previous observation and the fact that $G|_T=e(G|_S)$ is highly transitive it follows that for any subset $F'$ of $T$ of size $k$ there exists a finite subset $F''$ of $S$ such that $\acl(F'')\cap T=F'$.

	We claim that for all $x\in T$ there is a finite set $F$ of $S$ such that for all $\alpha\in (G|_S)_F$ it holds that $e(\alpha)(x)=x$. Let $F_1'$ and $F_2'$ be subsets of $T$ of size $k$ such that $F_1'\cap F_2'=\{x\}$. Then as we have seen there exist finite subset $F_1''$ and $F_2''$ of $S$ such that $\acl(F_i'')\cap T=F_i'$ for $i=1$ and $i=2$. Now let $F:=F_1''\cup F_2''$. Then if $\alpha\in (G|_S)_F$, then $\alpha\in (G|_S)_{F_i''}$, so $e(\alpha)(F_i')=F_i'$. Therefore $e(\alpha)(x)\in F_1'\cap F_2'=\{x\}$, that is $e(\alpha)(x)=x$.

	Now let $(\alpha_j)_j$ be a convergent sequence in $G|_S$. We want to show that the sequence $(e(\alpha_j))_j$ is also convergent, i.e., for all $x\in T$ we have $e(\alpha_j)(x)=e(\alpha_{j+1})(x)$ if $j$ is large enough. By our claim it follows that there is a finite set $F\subset S$ such that for all $\alpha\in (G|_S)_F$ we have $e(\alpha)(x)=x$. Since $(\alpha_j)_j$ a convergent there is an index $l$ such that $\alpha_j(y)=\alpha_{j+1}(y)$ for all $y\in F$ and $l\leq j$. Then if $j\geq l$ it follows that $\alpha_j\alpha_{j+1}^{-1}\in H|_F$, hence $e(\alpha_j)(e(\alpha_{j+1}))^{-1}(x)=e(\alpha_j\alpha_{j+1}^{-1})(x)=x$, and thus $e(\alpha_j)(x)=e(\alpha_{j+1})(x)$. Therefore $(e(\alpha_j))_j$ is convergent, which shows that $G|_S$ is closed.
	
	For $l < m-1$, note that 
	$S \subseteq P := O_1 \cup \cdots \cup O_{m-1}$. Hence, 
	$\acl(P) = X$ and we can apply the above argument for $P$ instead of $S$. We obtain that $G|_P$ is closed. Hence, the group $G|_P$
	satisfies all the assumptions for $G$
	but has fewer orbits, so by 
	induction we finally obtain that $G|_S$ is closed.
\end{proof}

In the proof of the next lemma it will be convenient to use a recent general result of Paolini and Shelah. 
A closed subgroup $G$ of $\sym(X)$ has the
\begin{itemize}
\item \emph{small index property} if every subgroup of $G$ of index less than $2^{\aleph_0}$ is open, i.e., contains the pointwise stabiliser of a finite set $F \subset X$.
\item \emph{strong small index property} if every subgroup of $G$ of index less than $2^{\aleph_0}$ lies between the pointwise and the setwise stabiliser of a finite set $F \subset X$.
\end{itemize}

The strong small index property of $\sym(X)$ 
itself has been shown in~\cite{DixonNeumannThomas}. (In fact, all automorphisms of $\omega$-categorical $\omega$-stable structures, and thus, by the results that we are about to prove, all groups in $\gexpp$, have the small index property~\cite{HodgesHodkinsonLascarShelah}.) 
On the other hand, already $\reduct(\unary)$ contains structures whose automorphism groups do not have the strong small index property (take e.g.~an equivalence relation with two infinite classes).
A permutation group $G$ on a set $X$ is said to have \emph{no algebraicity} if $\acl_G(Y) = Y$ for every $Y \subseteq X$.
The following has been proved in~\cite{PaoliniShelahReconstructing} (Corollary 2).

\begin{theorem}[\cite{PaoliniShelahReconstructing}]\label{paolini_shelah}
Let $X_1$ and $X_2$ be countable
and let $G_1 \leq \sym(X_1)$ and $G_2 \leq \sym(X_2)$ be closed oligomorphic groups
which 
have the strong small index property and have no algebraicity. Let us suppose that $\xi$ is a topological isomorphism from $G_1$ to $G_2$.  Then there exists
a bijection $b$ from $X_1$ to $X_2$ that \emph{induces $\xi$}, i.e., for all $x\in X_2$ and $\alpha\in G_1$ we have 
$$(\xi \alpha)(x) = b(\alpha(b^{-1}(x))).$$
 \end{theorem}


The small index property for $G$ implies that every homomorphism $h \colon G \to \sym(Y)$, for a countable set $Y$, is continuous (this is easy to see and well-known, see e.g.~\cite{Rosendal,Melleray}).
It follows from~\cite{Gaughan} that the image
of a continuous homomorphism from $\sym(X)$ to $\sym(Y)$ is closed in $\sym(Y)$ (see Theorem 1.3 in~\cite{YaacovTsankov} for a much more general recent result which also implies this).

\begin{lemma}\label{extension}
	Let $G$ be a closed oligomorphic permutation group on a countably infinite set $X$.
	Let $O_1,\dots,O_m$ be the orbits of $G$. Let us suppose that each $O_i$ is infinite, $G$ acts highly transitively on each $O_i$, and
for some $l < m$
	and $S := \bigcup_{i=1}^lO_i$
	we have $\acl_G(S)=X$ and $G|_S=\sym(O_1)\times \dots \times \sym(O_l)$. Then $\Delta(G)$ is not trivial.
\end{lemma}

\begin{proof}
	Let $j \in \{l+1,\dots,m\}$ 
	and $T:=O_j$. Let $G_j:=G|_{S\cup T}$. 
Then $G_j$ is closed by Lemma~\ref{closedness}
and we can apply Lemma \ref{two_parts} to the group $G_j$ with respect to the partition of $S \cup T$ into $S$ and $T$. Since $T \subseteq \acl_G(S)$ it follows that Condition (1) of \ref{two_parts} cannot hold. 
	Thus by Lemma \ref{two_parts} there exists a homomorphism $e_j \colon G|_S\rightarrow \sym(O_j)$ so that $G_j=\{(\alpha,e_j(\alpha)) \mid \alpha\in G|_S\}$. 
	For $i \in \{1,\dots,l\}$ and $\alpha \in G|_{O_i} = \sym(O_i)$ let $\hat{\alpha}_i$ denote the unique permutation of $S$ for which $\hat{\alpha}_i|_{O_i}=\alpha$ and $\hat{\alpha}_i|_{O_{k}}=\id_{O_{k}}$ if $k \neq i$.
	Then define the homomorphisms $$e_{ij} \colon \sym(O_i)\rightarrow \sym(O_j), \alpha \mapsto e_j(\hat{\alpha}_i).$$ 
	As mentioned before the lemma, 
	the map $e_{ij}$ is continuous.
	Let $H_i:=\{\hat{\alpha}_i \mid \alpha\in \sym(O_i)\}$. Then $H_i \triangleleft G|_S$, and so $e_j(H_i)\triangleleft e_j(G|_S)$.
	 By definition it follows that $e_j(G|_S)=(G_j)|_{O_j}=G|_{O_j}$. In particular $e_j(G|_S) \leq \sym(O_j)$ is highly transitive. 
Thus by Lemma \ref{highly_normal} it follows that either $e_j(H_i)$ is also highly transitive or it is trivial. If $e_{ij}$ is trivial for every $i \in \{1,\dots,l\}$, then $G$ fixes every element of $O_j$ contradicting the fact that $G$ acts transitively on $T=O_j$. Thus, there is an $i\in \{1,\dots,l\}$ such that the image $I\leq \sym(O_j)$ of $e_{ij}$ is highly transitive. As we have mentioned before the statement of the lemma, $I$ is a closed subset of $\sym(O_j)$, so we can apply Theorem \ref{paolini_shelah}
and obtain a bijection $b_{j}$ between $O_i$ and $O_j$ which induces $e_{ij}$ (we could as well have derived this from the argument in Example 2 on page 224 of~\cite{HodgesLong}). 
Now let $i'\neq i$, and let $\alpha$ be a nontrivial permutation of $O_{i'}$. Then $\hat{\alpha}_{i'}$ commutes with every element of $H_i$, and so $e_{i'j}(\alpha)$ commutes with every element of $e_j(H_i)=\sym(O_j)$. Therefore $e_{i'j}(\alpha)=\id_{O_j}$.

	We have obtained that for  every $j \in \{l+1,\dots,m\}$ there is a unique $i(j) \in \{1,\dots,l\}$ and a bijection $b_j \colon O_{i(j)} \rightarrow O_j$ such that for all $g\in G$ and $x\in O_j$ we have $g(x)=b_j(g(b_j^{-1}(x)))$. In fact, the same is true in the case when $j\leq l$: we can choose $i(j)$ to be $j$, and $b_j$ to be the identity map. Let $b$ be the union of the functions $b_1,\dots,b_m$ and define the relation $\sim$ by $x\sim y\Leftrightarrow b(x)=b(y)$. Then $\sim$ is a congruence of $G$ all of whose classes are finite. Moreover, $\sim$ is nontrivial since $m>l$. This implies that also $\Delta(G)$ is nontrivial.
\end{proof}

\begin{lemma}\label{two_cuts}
	Let $H$ be an oligomorphic permutation group on a countably infinite set $X$ with two 
	infinite orbits $Y$ and $Z$. Let us assume that $H$ acts 2-transitively on $Z$ and that there exists $y\in Y$ so that $|\nabla((H_y)|_{Z})|\geq 2$. Then for every $n \in {\mathbb N}$ there exist $y_1,\dots,y_n\in Y$ such that $|\nabla((H_{y_1,\dots,y_n})|_{Z})|\geq n+1$.
\end{lemma}

\begin{proof}
	By transitivity we know that for all $y'\in Y$ it holds that $|\nabla(H_{y'}|_{Z})|\geq 2$.

	We show the statement of the lemma by induction on $n$. For $n=1$ the statement is trivial. Now suppose that we know there exist $y_1,\dots,y_{n-1}\in Y$ such that $\nabla((H_{y_1,\dots,y_{n-1}})|_{Z})$ has at least $n$ classes. Let $C_1,\dots,C_m$ be the classes of $\nabla((H_{y_1,\dots,y_{n-1}})|_{Z})$. Let $z_1,z_2\in C_1$. Since $H$ acts 2-transitively on $Z$, it follows that there exists a $y_n\in Y$ such that $(z_1,z_2)\not\in \nabla(H_{y_n}|_{Z})$. Let $$R:=\nabla((H_{y_1,\dots,y_{n-1}})|_{Z})\cap \nabla(H_{y_n}|_{Z}).$$ Then $R$ is a congruence of $H_{y_1,\dots,y_n}$ which is strictly finer than $(H_{y_1,\dots,y_{n-1},y_n})$, and has finitely many classes. Therefore, $\nabla((H_{y_1,\dots,y_{n}})|_{Z})$ is also finer than $\nabla((H_{y_1,\dots,y_{n-1}})|_{Z})$. In particular,  $|\nabla((H_{y_1,\dots,y_{n}})|_{Z})|\geq m+1\geq n+1$.
\end{proof}

\begin{lemma}\label{gen_sym}
		Let $G \in \gexpp$ be closed. 
	Suppose that $\Delta(G)$ is trivial and that $G$ stabilises each class of $\nabla(G)$ setwise. Let $O_1,\dots,O_l$ be the orbits of $G$. 
	Suppose that each $O_i$ is infinite, and that $S := X\setminus O_l$ is algebraically closed. Then 
	$\{\id_S\}\times \sym(O_l) \subseteq G$.
\end{lemma}

\begin{proof}
	 By our assumptions, the orbits $O_1,\dots,O_l$ are the classes of $\nabla(G)$.
	 By Lemma \ref{nice_case} it follows that $G$ acts highly transitively on each orbit $O_i$. 
	 We apply Lemma~\ref{two_parts} to 
	 $S$ and $T := O_l$. 
	 If Item (1) of Lemma~\ref{two_parts} applies
	 then we are done. 
	 We claim that item (2) of Lemma~\ref{two_parts} cannot hold. 
	 For this, it suffices to show that 
	 $G$ contains a nontrivial permutation $\gamma$ such that $\gamma|_S=\id_S$. 	
	 Indeed, if 
	 there exists a homomorphism 
	 $e \colon G|_S \to G|_T$ as in item (2) of Lemma \ref{two_parts}, then (using the notation from the proof of Lemma~\ref{two_parts})
	 $\gamma = (\id_S,e(\id_S)) = (\id_S,\id_T)$
	 is trivial. 
	 	
	\medskip 
	{\bf Claim 1.}  For every finite $F\subseteq S$ and $L := G_F|_T$, the congruence $\nabla(L/\Delta(L))$ is universal.

\begin{proof}[Proof of Claim 1]
	Suppose to the contrary 
	that there exists a finite $F \subseteq S$
	such that $\nabla(L/\Delta(L))$ is not universal; choose a minimal $F$ with this property. 
	For some $y\in F$, put $F' := F \setminus \{y\}$ and 
	\begin{align*}
	E & := G_{F'}|_T \\ 
	Z & :=T/\Delta(E) \\
	\text{ and } K & := E/\Delta(E) \leq \sym(Z).
	\end{align*} 
	Let us consider the mapping $\pi \colon T\rightarrow Z$ which maps each element to its $\Delta(E)$-class. Then if $u, v \in Z^n$ are in different orbits of $K$, 
	then all tuples in $\pi^{-1}(u)$ are in different orbits of $G$ than the tuples from $\pi^{-1}(v)$. Moreover, if $u\in Z^n$ is injective, then so are all the tuples in $\pi^{-1}(u)$. This means that the number of injective $n$-orbits of $K$ is at most $\oi_n(G_{F'})$. By Corollary~\ref{stabil_finite} it follows that $\oi_n(G_{F'}) \leq cn^{dn}$ 
	for some constants $c,d$ with $d<1$. Therefore, $\oi_n(K)\leq cn^{dn}$ and thus $K \in \gexpp$. From the definition of $K$ it is clear that $\Delta(K)$ is the identity congruence. 
	It follows from the minimal choice of $F$ that the congruence $\nabla(K)$ is the universal congruence. 
	 Therefore, Lemma \ref{nice_case} implies that $K$ is highly transitive.

	Now let $Y:=O_j\setminus F'$ where 
	$j$ is such that $y \in O_j$. 
	Then $G_{F'}$ acts naturally on $Y \sqcup Z$. Let $H$ be the image of this action (as a subgroup of $\sym(Y \sqcup Z)$); note that $H$ is highly transitive on $Z$. We claim that the group $H$, the orbits $Y,Z$ and the element $y\in Y$ satisfy the conditions of Lemma \ref{two_cuts}. The only nontrivial fact that we have to check is that the congruence $\nabla(H_y|_Z)$ is not universal. We know that $\nabla(L/\Delta(L))$ is not universal. By Lemma \ref{delta_union_nabla} this is equivalent to the fact that the congruence generated by $\nabla(L)$ and $\Delta(L)$ is not universal. The congruence $\Delta(E)$ is also a congruence of $L$ with finite classes, hence $\Delta(E)$ is finer than $\Delta(L)$. Hence, the congruence generated by $\nabla(L)$ and $\Delta(E)$ is finer than the congruence generated by $\nabla(L)$ and $\Delta(L)$, and thus it is also not universal. Using Lemma \ref{delta_union_nabla} again it follows that $\nabla(L/\Delta(E))$ is not universal. Since $L=(G_{F'})_{y}|_T$ this means that the congruence $\nabla(H_y|_Z)$ is not universal. 

	If we apply Lemma \ref{two_cuts} we obtain that for every $n$ there exist $y_1,\dots,y_n\in Y$ so that $|\nabla(H_{y_1,\dots,y_{n}}|_{Z})| \geq n+1$. This also implies that $|\nabla(G_{F'\cup \{y_1,\dots,y_{n}\}})| \geq n+1$, that is, $m_G(F'\cup \{y_1,\dots,y_n\}) \geq n+1$. In particular $m_G(|F'|+n)\geq n+1$. This contradicts Lemma~\ref{not_too_many_classes} for  $0<\varepsilon<1-d$ if $n$ is large enough.
\end{proof}

	Claim 1 implies that for every finite $F \subseteq S$ the group 
	$G_F$ acts transitively on $T/\Delta(G_F|_T)$. This also implies that all $\Delta(G_F)$-classes contained in $T$ have the same size. For a finite set $F\subset S$, let $k(F)$ denote this size. By Corollary \ref{union_delta_size} we have $G\in (\gexpp)^k$ for some $k$. This implies that $k(F)\leq k$ for every finite subset $F$ of $S$. This also implies that there exists a finite set $F$ so that $k(F)$ is maximal. So let us choose $F\subset S$ so that $k(F)$ is maximal, and, similarly as above, let $L := G_{F}|_T$,  $Z:=T/\Delta(L)$, and let $\pi \colon T \rightarrow Z$ map each element to its $\Delta(L)$-class.
	
	\medskip
	{\bf Claim 2.} For any finite $F' \subset S$ that contains $F$ the group $G_{F'}$ acts highly transitively on $Z$. 
	
	\begin{proof}[Proof of Claim 2.] 
Let 
$E := G_{F'}|_T$ and let $K:=E/\Delta(E)$.
By the maximality of $k(F)$ it follows that 
$\Delta(G_F)$ and $\Delta(G_{F'})$ agree on $T$,
and hence $T/\Delta(E) = T/\Delta(L)$ and thus $K \leq \sym(Z)$.   
	As in the proof of Claim 1, 
	it follows that $K \in \gexpp$, and Claim 1 implies that $\nabla(K)$ is the universal congruence. So it is enough to show that 
	$\Delta(K)$ is the identity congruence. In order to show this, let us consider the relation on $X$ 
	$$R:=\{(x,y) \in T^2 \mid (\pi(x),\pi(y)) \in \Delta(K)\} \cup \{(x,x) \mid x\in S\}.$$ 
Then $R$ is a congruence of $G_{F'}$ with finite classes. Therefore, $R$ is finer than $\Delta(G_{F'})$. As $\Delta(G_F)$ and $\Delta(G_{F'})$ agree on $T$, this is only possible if $\Delta(K)$ is the identity congruence. Hence, the conditions of Lemma \ref{nice_case} hold for $K$, and thus by Lemma \ref{nice_case} it follows that $K$ is highly transitive.
	\end{proof}

	Now let us choose a prime $p>k(F)$ and let $z_1,\dots,z_p\in Z$. Let 
	\begin{align*}
	F & =S_0\subset S_1\subset \cdots \\
	\text{ and } \quad 
	\pi^{-1}(z_1)\cup \cdots \cup \pi^{-1}(z_p) & = T_0\subset T_1 \subset \cdots
	\end{align*} be sequences of finite subsets of $S$ and $T$, respectively, so that $\bigcup{S_i}=S$ and $\bigcup{T_i}=T$. By Claim 2, the stabiliser $G_{S_i}$ acts highly transitively on $Z=T/\Delta(L)$. In particular, there is a permutation $\gamma_i' \in G$ which fixes every element in $S_i$ and which acts on $T_i/\Delta(L)$ as $(z_1z_2\dots z_p)$. Now let $\gamma_i := (\gamma_i')^{k(F)!}$. Then $\gamma_i|_{T_0}$ is nontrivial, but $\gamma_i|_{T_i\setminus T_0}=\id_{T_i\setminus T_0}$. 

	Since the permutations $\gamma_i$ have finitely many possible actions on the set $T_0$, we can assume, by choosing a subsequence if necessary, that $\gamma_i|_{T_0}$ are the same for all $i$. Then the permutations $\gamma_i$ converge to a permutation $\gamma$ for which $\gamma|_{S\cup T\setminus T_0}$ is trivial, but $\gamma|_{T_0}$ is not trivial. Since $G$ is closed it follows that $\gamma\in G$. This finishes the proof of the lemma.
\end{proof}

\begin{corollary}\label{gen_sym_2}
	Let $G \in \gexpp$ be closed. 
	Suppose that $\Delta(G)$ is trivial and that $G$ stabilises every class of $\nabla(G)$ setwise. Let $O_1,\dots,O_l$ be the orbits of $G$. Suppose that each $O_i$ is infinite and that $X\setminus O_i$ is algebraically closed. Then $G=\sym(O_1)\times \dots \times \sym(O_l)$.
\end{corollary}

\begin{proof}
	Direct consequence of Lemma \ref{gen_sym}.
\end{proof}

\begin{lemma}\label{extension_2}
			Let $G \in \gexpp$ be closed and such that $\Delta(G)$ is trivial. Suppose that $G$ fixes every class of $\nabla(G)$ setwise. Let $O_1,\dots,O_m$ be the orbits of $G$ and suppose that each $O_i$ is infinite. Then $G=\sym(O_1)\times \dots \times \sym(O_m)$.
\end{lemma}
\begin{proof}
	Let $I \subseteq \{1,\dots,m\}$ be minimal so that $\acl_G(S)=X$ for $S :=\bigcup_{i\in I}{O_i}$. Without loss of generality we can assume that $I=\{1,\dots,l\}$ for some $l \leq m$. By Lemma \ref{closedness} the group $G|_S \in \gexpp$ is closed. 
	By the minimality of $I$ it follows that the sets $S\setminus O_i$, for $i \in \{1,\dots,l\}$,  are algebraically closed with respect to $G$.
	By Corollary \ref{gen_sym_2} it follows that $G|_S=\sym(O_1)\times \dots \times \sym(O_l)$. By Lemma \ref{nice_case} it follows that $G$ acts highly transitively on each orbit $O_i$.  Thus, if $l<m$ then Lemma \ref{extension} implies that $\Delta(G)$ is nontrivial. This means that $l=m$, and thus $G=G|_S=\sym(O_1)\times \dots \times \sym(O_m)$.
\end{proof}

	In order to drop the condition that $G$ fixes every $\nabla(G)$-class we need the following observation about finite index subgroups of oligomorphic groups.

\begin{proposition}\label{easy}
	Let $G$ be an oligomorphic 
	permutation group on a countably infinite set $X$, and let $H$ be a finite-index subgroup of $G$. 
	Then $$\oi_n(H) \leq [G:H] \cdot \oi_n(G).$$ In particular, $H$ is oligomorphic.
\end{proposition}
\begin{proof}
	Choose elements $\gamma_1,\dots,\gamma_{[G:H]}\in G$ such that $G=\bigcup_{i=1}^{[G:H]}{\gamma_i H}$. If the tuples $t_1,t_2,\dots,t_l$ represent all injective $n$-orbits of $G$, then the tuples $\gamma_it_j$ for $1\leq i\leq {[G:H]}$ and 
	$1\leq j\leq l$ represent all injective $n$-orbits of $H$. Therefore, $\oi_n(H)\leq [G:H] \cdot \oi_n(G)$.
\end{proof}

\begin{lemma}\label{finite_index_acl}
	Let $G \leq \sym(X)$ be an oligomorphic permutation group and let $H$ be a finite-index subgroup of $G$. Then $\acl_H(x)\subseteq \acl_G(x)$ for all $x\in X$.
\end{lemma}

\begin{proof}
	By Proposition \ref{easy} the permutation group $H$ is oligomorphic. 
	Let $x\in X$. We claim that the index $k := [G_x:H_x]$ is finite.  
	Let $K:=\bigcap_{g\in G} g^{-1} H g$, 
that is, $K$ is the kernel of the left action of $G$ on the left cosets of $H$ in $G$. Then $K \leq H$, and $|G/K|=[G:K]$ is finite. 
Then by the second isomorphism theorem we have $G/K \supseteq G_xK/K\simeq G_x/(K\cap G_x)=G_x/K_x$. In particular, $[G_x:K_x]=|G_x/K_x|$ is finite. Since $K\leq H$ we have $K_x\leq H_x$. Therefore, $[G_x:H_x]$ is also finite.

	We have shown that $k := [G_x:H_x]$ is finite. Choose elements $\gamma_1,\dots,\gamma_k \in \aut(\fa)$ such that $G_x=\bigcup_{i=1}^{k}{\gamma_i H_x}$. Let $y\in \acl_H(x)$. Since $[H:H_x]$ is finite, $H_x(y)$ is finite. Therefore, $$G_x(y)=\left(\bigcup_{i=1}^{k}{\gamma_i H_x}\right)(y)=\bigcup_{i=1}^{k}{\gamma_i H_x(y)}$$ is finite, that is, $y\in \acl_G(x)$. This proves that $\acl_H(x) \subseteq \acl_G(x)$.
\end{proof}

\begin{lemma}\label{finite_index_delta}
	Let $G$ be an oligomorphic permutation group on a countably infinite set $X$ and let $H$ be a finite-index subgroup of $G$. Then $\Delta(G)=\Delta(H)$.
\end{lemma}

\begin{proof}
	Clearly, $\Delta(G)$ is a congruence of $H$ with finite classes. Thus, $\Delta(G) \subseteq \Delta(H)$. Now let $(x,y)\in \Delta(H)$. Then $y\in \acl_H(x)$ and $x\in \acl_H(y)$ by Lemma \ref{delta_alt}. By Lemma \ref{finite_index_acl} this implies that $y\in \acl_G(x)$ and $x\in \acl_G(y)$. Again by Lemma \ref{delta_alt} we obtain $(x,y)\in \Delta(G)$, showing that $\Delta(G)=\Delta(H)$. 
\end{proof}


\begin{theorem}\label{extension_3}
		Let $G \in \gexpp$ be closed such that that $\Delta(G)$ is trivial. Let $O_1,\dots,O_m$ be the classes of $\nabla(G)$. Then $\sym(O_1)\times \dots \times \sym(O_m) \subseteq G$.
\end{theorem}

\begin{proof}[Proof of Theorem \ref{extension_3}]
	Let $K$ be the kernel of the action of $G$ on $\{O_1,\dots,O_m\}$. Then $[G:K]$ is finite, and thus by Proposition \ref{easy} it follows that $K\in \gexpp$. Without loss of generality we can assume that 
	$O_1,\dots,O_l$ are the infinite orbits of $K$.
	Let $Y=O_1 \cup \dots \cup O_l$. Then $X\setminus Y$ is finite, and $K$ fixes each element in $X\setminus Y$. By Lemma \ref{stabil_finite} it follows that the group $K|_{Y}$ is in $\gexpp$. By Lemma \ref{finite_index_delta} it follows that $\Delta(K)$ is trivial, and thus $\Delta(K|_{Y})$ is also trivial. Moreover, $K|_{Y}$ fixes every class of $\nabla(K|_{Y})$ setwise and all orbits of $K|_{Y}$ are infinite. Hence, we can apply Lemma~\ref{extension_2} and obtain that $K|_{Y}=\sym(O_{1})\times \dots \times \sym(O_l)$. Therefore, $\sym(O_1)\times \dots \times \sym(O_m) = K \subseteq G$.
\end{proof}

\begin{corollary}\label{extension_4}
	Let $\fa\in \kexpp$ be such that $\Delta(\fa)$ is trivial. Then $\fa\in \reduct(\unary)$.
\end{corollary}
\begin{proof}
	Apply Theorem \ref{extension_3} to $\aut(\fa)$ and combine with Corollary~\ref{cor:unary_reduct-iff}. 
\end{proof}

\subsection{The general case}
\label{sect:general}

\begin{lemma}\label{desc_kexpp}
	$\kexpp \subseteq \fcover(\reduct(\unary))$.
\end{lemma}
\begin{proof}
	Let $\fa\in \kexpp$ and let us consider the factor mapping $\pi \colon \fa\rightarrow \fb$ where $\fb := \fa/\Delta(\fa)$. Recall that by $\fa/\Delta(\fa)$ we denote any structure with $\aut(\fa/\Delta(\fa))=\aut(\fa)/\Delta(\fa)$. If $u,v \in B$ are in different orbits of $\aut(\fb)$ then the tuples in $\pi^{-1}(u)$ lie in different orbits  of $\aut(\fa)$ than the tuples in $\pi^{-1}(v)$. Moreover, if $u\in B^n$ is injective, then so are the tuples in $\pi^{-1}(u)$. This means that the number of injective $n$-orbits of $\fb$ is at most $\oi_n(\fa)$ and thus $\fb \in \kexpp$. Then $\Delta(\fb)$ must be trivial:
	otherwise, $\aut(\fb)$ has a nontrivial congruence all of whose classes are finite, 
	contradicting the definition of $\Delta(\fa)$.
	By Corollary \ref{extension_4} it then follows that $\fb \in \reduct(\unary)$, and thus $\fa \in \fcover(\reduct(\unary))$.
\end{proof}

	The reverse containment holds as well. 

\begin{theorem}\label{main_kexpp}
	$\kexpp=\fcover(\reduct(\unary))=\reductfin(\fcover(\unary^*))$
\end{theorem}

\begin{proof}
	We already know that $ \kexpp \subseteq \fcover(\reduct(\unary))$ (Lemma~\ref{desc_kexpp})
	and that $\fcover(\reduct(\unary)) \subseteq \reductfin(\fcover(\unary^*))$ (see Remark \ref{fru_rffu}). So we have to show that $\reductfin(\fcover(\unary^*))\subseteq \kexpp$. Proposition \ref{prop:easy} implies  that $\reductfin(\kexpp)=\kexpp$. Therefore it is enough to show that $\fcover(\unary^*)\subseteq \kexpp$. 
	So let $\fb\in \unary^*$ and let $\pi \colon \fa\rightarrow \fb$ be a finite covering. Lemma \ref{reduct_trivial2} shows that 
	$\pi$ is strongly split. 
	Therefore we can assume that $\pi$ is a strongly trivial covering map. 
	
	It follows from the description of trivial coverings given
	in Remark~\ref{triv_cov_unary} 
	that the orbit of an injective $n$-tuple $t=(t_1,\dots,t_n)$ of a trivial covering of a unary structure is uniquely determined by the orbits of $t_1,\dots,t_n$ and by the partition of the set $\{t_1,\dots,t_n\}$ defined by the congruence $\sim_{\pi}$. This means that the number of injective orbits of $\fa$ it at most $m^n \cdot p_k(n)$ where 
	\begin{itemize}
	\item 
	$m$ is the number of orbits of $\aut(\fa)$, 
	\item $k$ is the maximal size of the classes of $\sim_{\pi}$, and 
	\item $p_k(n)$ is the number of  partitions of $\{1,\dots,n\}$ with parts of size at most $k$ (see Section~\ref{sect:growth}). 
	\end{itemize}
	Let us choose $d > \frac{k-1}{k}$. 
	Then by Lemma \ref{counting_2} we have 
	$p_k(n) < c_1 n^{dn}$ for some $c_1$. Thus $\oi_n(\fa)\leq m^n c_1 n^{dn} \leq c_2 n^{dn}$ for some $c_2$. Therefore $\fa\in \kexpp$.
\end{proof}

\begin{remark}\label{rem:triv-cover-interpretation}
Recall from Proposition~\ref{reduct_trivial2}
that every finite cover $\pi \colon \fa \to \fb$ for $\fb \in \unary^*$ is strongly split,
and hence all structures in $\kexpp = \reduct(\fcover(\unary^*))$ have a first-order interpretation in $({\mathbb N};=)$ (Remark~\ref{rem:triv-cover-interpret}).  
Since $({\mathbb N};=)$ is $\omega$-stable and first-order interpretations preserve $\omega$-stability,
it follows that all structures in $\kexpp$ are $\omega$-stable.
\end{remark}

\subsection{Thomas' conjecture for the class $\kexpp$}
\label{sect:thomas}
	Let $k,m\in \mathbb{N}$. Let $\group(k,m)$ be the class of all oligomorphic permutation groups $G$ for which 
	the classes of $\Delta(G)$ have size at most $k$ 
	and where $\nabla(G/\Delta(G))$ has at most $m$ classes.
	Let $\struct(k,m)$ be the class of all structures whose automorphism group is in $\group(k,m)$.

\begin{lemma}\label{skm_classes}
	Let $k,m\in \mathbb{N}$. Let $\fb\in \unary^*$, let $\pi \colon\fa\rightarrow \fb$ be a finite covering, and let $\fc$ be a quasi-covering reduct of $\fb$. Then $\fc\in \struct(k,m)$ iff $\fa\in \struct(k,m)$.
\end{lemma}

\begin{proof}
	By definition $\Delta(\fa)=\Delta(\fc)={\sim_{\pi}}$, and $$\nabla(\aut(\fc/\Delta(\fc))=\nabla(\aut(\fc/\Delta(\fa))=\nabla(\aut(\fa/\Delta(\fa))=\nabla(\fb).$$
\end{proof}

\begin{lemma}\label{kexpp_limited}
	Let $k,m\in \mathbb{N}$. There are finitely many structures in $\kexpp\cap \struct(k,m)$ up to bi-definability. 
\end{lemma}
\begin{proof}
	By Theorem \ref{main_kexpp} we know that $\kexpp=(\fcover\circ \reduct)(\unary)$. Proposition~\ref{quasi_trivial} implies that every structure $\fc$ in $\kexpp$ is a quasi-covering reduct of a finite covering structure $\fa$ of some structure in $\unary^*$. By Lemma \ref{skm_classes} we know that $\fc\in \struct(k,m)$ if and only if $\fa\in \struct(k,m)$. 
	If $\fa$ is a trivial covering of some structure in $\unary^*$, then by Theorem \ref{quasi_cover2}
	it has finitely many quasi-covering reducts. Therefore it is enough to show that there are finitely many structures in $\struct(k,m)$ up to bi-definability which are strongly trivial covering structures of some structure in $\unary^*$.

	Let $\fb\in \unary^*$ and let $\pi \colon\fa \rightarrow \fb$ be a strongly trivial finite covering map. 
	Let $O_1,\dots,O_l$ be the orbits of $\fb$. Then $l\leq m$. 
	Following Remark \ref{triv_cov_unary} we can assume without loss of generality that $A=\bigsqcup_{i=1}^l{F_i\times O_i}$ for some finite sets $F_i$, and $$\aut(\fa)=\bigsqcup_{i=1}^l{\id_{F_i}\wr \sym(O_i)}.$$ Since $\sim_{\pi}$ is a congruence with finite classes it follows that $|F_i|\leq k$. Then there are finitely many options for $l$, the sizes of the orbits $O_i$ (they are all either one or infinite), and the sizes of the sets $F_i$, and if we fix these parameters, then the group $\aut(\fa)$ is uniquely determined up to isomorphism. This implies that there are finitely many structures in $\struct(k,m)$ up to bi-definability which are a trivial covering structure of a structure in $\unary^*$.
\end{proof}

\begin{lemma}\label{kexpp_km}
	Let $\fa\in \kexpp$ and let $\fb \in \reduct(\fa)$. Let $k$ be the size of the largest $\Delta(\fa)$-class and let $m$ be the number of $\nabla(\fa)$-classes. Then $\fb\in \struct(k,m)$.
\end{lemma}
\begin{proof}
	If $R$ is a congruence of $\aut(\fb)$, then it is also a congruence of $\aut(\fa)$. Therefore, the size of every class of $\Delta(\fb)$ is at most $k$.
	Similarly, the number of $\nabla(\fb)$-classes
	is at most the number of $\nabla(\fa)$-classes. The number of 
	$\nabla(\fb)$-classes is an upper bound for the number of $\nabla(\fb/\Delta(\fb))$-classes. This proves the lemma.
\end{proof}

	Lemmas \ref{kexpp_limited} and \ref{kexpp_km} immediately imply the following weak version of Thomas' conjecture for the class $\kexpp$.

\begin{theorem}\label{thomas_weak}
	Let $\fa\in \kexpp$. Then $\fa$ has finitely many first-order reducts up to bi-definability.
\end{theorem}

	Then the (standard version of) Thomas' conjecture follows as follows. 
	First we state an important well-known link between
infinite descending chains of first-order reducts and infinite signatures. We say that 
a structure $\fb$ has \emph{essentially
infinite signature} if there does not exist
a structure $\fb'$ with finite signature
such that $\aut(\fb) = \aut(\fb')$. 

\begin{lemma}\label{lem:inf-sign}
Let $\fa$ be an $\omega$-categorical structure. Then 
there exists an infinite sequence $\fb_1,\fb_2,\dots$ of first-order reducts of $\fa$ such that $\aut(\fb_1) \supsetneq \aut(\fb_2) \supsetneq \cdots$ 
if and only if $\fa$ has a reduct with essentially infinite signature. 
\end{lemma}
\begin{proof}
Assume that the reduct $\fb = (B;R_1,R_2,\dots)$ of $\fa$ has essentially infinite signature. 
By assumption, 
$\fb$ and $\fb_n := (B;R_1,\dots,R_n)$ are not first-order interdefinable. 
Moreover, for every $n \in {\mathbb N}$ there exists an $f(n) \in {\mathbb N}$ such that $\fb_n$ and $\fb_{f(n)}$ are not first-order interdefinable (otherwise, every relation in $\fb$
would be first-order definable in $\fb_n$, contradicting our assumptions). 
So $\fb_1,\fb_{f(1)},\fb_{f(f(1))},\dots$ provides an infinite strictly descending
chain of first-order reducts of $\fa$. 

Suppose conversely that $\fb_1,\fb_2,\dots$ is an infinite strictly descending chain of first-order reducts of $\fa$.  Define $\fc$ as the first-order reduct of $\fa$ whose relations are precisely the relations of all the $\fb_i$.  Assume for contradiction that there exists a finite-signature structure $\fc'$ with $\aut(\fc') = \aut(\fc)$.
Let $i \in {\mathbb N}$ be such that all relations used in the definitions of the relations of $\fc'$ in $\fc$ already
appear in the signature of $\fb_i$. 
Then $\aut(\fb_i) = \aut(\fc') = \aut(\fc) = \aut(\fb_j)$ for all $j \geq i$, contradicting the assumption that $(\fb_i)_{i \in {\mathbb N}}$ is strictly decreasing. 
\end{proof} 

\begin{proposition}\label{prop:interdef-bidef}
Let $\fa$ be an $\omega$-categorical structure.
Then $\fa$ has finitely many first-order reducts up to interdefinability if and only if $\fa$ has finitely many first-order reducts up to bidefinability. 
\end{proposition}
\begin{proof}
If $\fb$ is a first-order reduct of $\fa$ with essentially infinite signature, then $\fb$ has
an infinite strictly descending chain of first-order reducts. Note that if $\aut(\fb_1) \subsetneq \aut(\fb_2)$ then for some $n$ there are strictly more orbits of $n$-tuples in $\aut(\fb_1)$ than in $\aut(\fb_2)$, so $\fb_1$ and $\fb_2$ are not bidefinable  (if two reducts are bidefinable then they have the same number of orbits of $n$-tuples for all $n$). So $\fb$ and $\fa$ have 
infinitely many first-order reducts up to bi-definability,  so the statement is trivially true in this case. 

Therefore it
suffices to show that
every first-order reduct $\fb$ of $\fa$ with finite signature is bidefinable to at most finitely many reducts of $\fb$ up to interdefinability. 
The equivalence class of $\fb$ with respect to interdefinability is given by 
its orbits of $n$-tuples, for some finite $n$ (since $\fb$ has finite signature), and thus the same holds for any structure which is bidefinable $\fb$. 
Since $\fa$ is $\omega$-categorical,
there are finitely many orbits of $n$-tuples in $\fa$, which implies that there are finitely many first-order reducts of $\fa$ up to interdefinability that
are bidefinable with $\fb$. 
\end{proof}

\begin{theorem}\label{thomas_strong}
	Let $\fa \in \kexpp$. Then $\fa$ has finitely many first-order reducts.
\end{theorem}

\begin{proof}
	Follows from Theorem \ref{thomas_weak} and Proposition~\ref{prop:interdef-bidef}.
\end{proof}

\begin{corollary}
$\kexpp$ contains countably many structures up to interdefinability. It contains no structure with  essentially infinite signature. 
\end{corollary}
\begin{proof}
The first statement is implied by Lemma~\ref{kexpp_limited} 
in combination with Proposition~\ref{prop:interdef-bidef}. The second statement follows from Theorem~\ref{thomas_strong} 
and Lemma~\ref{lem:inf-sign}. 
\end{proof}
\end{section}

%% file: partition-growth.tex
\subsection{Growth rates for partitions}
\label{sect:growth}
For $n,k \in {\mathbb N}$, let $p_k(n)$ be the number of partitions of the set $\{1,\dots,n\}$ with parts of size at most $k$; this is the Sloane integer sequence A229223. 
Asymptotic formulas for $p_k(n)$ are known for 
$k \in \{1,\dots,4\}$ (called \emph{allied Bell numbers} in a letter of John Riordan). 
We need an upper and a lower bound for all $k \in {\mathbb N}$.

\begin{lemma}\label{counting}
	Let $\varepsilon>0$. Then 
	$p_k(n) \geq n^{(\frac{k-1}{k}-\varepsilon)n}$ if $n$ is large enough.
\end{lemma}

\begin{proof} 
	Let $s_k(n)$ be the number of partitions of $\{1,\dots,kn\}$ where all the parts contain exactly $k$ elements. Clearly, $s_k(1)=1$ for all $k \in {\mathbb N}$. 
	To form a partition of $\{1,\dots,kn\}$ for $n>1$ we first choose the class containing the number $kn$, and then we choose a partition of the remaining elements. Hence, $s_k$ satisfies the recursion $$s_k(n)={{kn-1}\choose k-1} {s_{k}(n-1)}.$$
	Since ${{kn-1}\choose k-1} \geq n^{k-1}$ 
	we obtain by induction that $$s_k(n) \geq n^{k-1}(n-1)^{k-1} \cdots 2^{k-1}=(n!)^{k-1}.$$
	Stirling's formula ($n! \sim \sqrt{2\pi n} (\frac{n}{e})^n$ for $n$ tending to infinity) implies that
	$$(n!)^{k-1} \geq n^{(k-1)(1-\varepsilon')n}$$ for any $\varepsilon'>0$ if $n$ is large enough. Hence,  \begin{align*}
	p_k(n) & \geq s_k(\lfloor \frac{n}{k} \rfloor)
	 \geq \lfloor \frac{n}{k} \rfloor^{(k-1)(1-\varepsilon')\lfloor \frac{n}{k} \rfloor} \\
	& \geq \Bigl(\frac{n}{k}-1\Bigr)^{(k-1)(1-\varepsilon')(\frac{n}{k}-1)}\geq n^{(k-1)(1-\varepsilon')(1-\varepsilon'')\frac{1}{k}n}\geq n^{(\frac{k-1}{k}-\varepsilon)n}
	\end{align*}
	for an appropriate choices of $\varepsilon',\varepsilon''>0$ if $n$ is large enough.
\end{proof}

\begin{lemma}\label{counting_2}
	Let $n,k \in {\mathbb N}$. If $d>\frac{k-1}{k}$, then $p_k(n) < cn^{dn}$ for some $c$.
\end{lemma}

\begin{proof}
	To form a partition of $\{1,\dots,n\}$ for $n>1$, we first choose the class containing the number $n$, and then we choose a partition of the remaining elements. 
	We thus have the following recursion formula: 
	\begin{align}
	p_k(n) =& \sum_{i=0}^{k-1}{{n-1\choose i}p_k(n-1-i)}. 
	\label{recursion}
	\end{align}
 We claim that the following inequality holds if $n$ is large enough. 
 \begin{align}
 \sum_{i=0}^{k-1}{{n-1\choose i}(n-1-i)^{d(n-1-i)}}< n^{dn}.\label{main_ineq} 
 \end{align}
 In order to prove this it is enough to show that if $i\leq k-1$ and $n$ is large enough, then \begin{align}{n-1\choose i}(n-1-i)^{d(n-1-i)}<\frac{1}{k}n^{dn},\end{align} 
 that is, 
 \begin{align}{n-1\choose i}<\frac{1}{k}\bigg(\frac {n^{n}}{(n-1-i)^{n-1-i}}\bigg)^d.
 \label{what_we_need} 
 \end{align}
We have \begin{align*}\frac{n^n}{(n-1-i)^{n-1-i}} & =\prod_{j=0}^{i}\frac{(n-j)^{n-j}}{(n-j-1)^{n-j-1}}\\
 & =\prod_{j=0}^{i}\bigg((n-j)\Big(\frac{n-j}{n-j-1}\Big)^{n-j-1}\bigg) \\
 & \geq \prod_{j=0}^{i}(n-j)=n(n-1) \cdots (n-i).\end{align*}
This implies that in order to show Inequality (\ref{what_we_need}) it is enough to show that $${n-1\choose i}<\frac{1}{k}(n(n-1)\cdots (n-i))^d$$ if $n$ is large enough. By rearranging the inequality above we obtain that it is equivalent to the following inequality. \begin{equation}\frac{1}{i!}((n-1) \cdots (n-i))^{1-d}<\frac{1}{k}n^d
\label{final_ineq} 
\end{equation}
The LHS of the Inequality~(\ref{final_ineq}) is asympotically $\frac{1}{i!}n^{i(1-d)}$. 
	By our assumption $d>\frac{k-1}{k}$, thus $i(1-d)<\frac{i}{k}\leq \frac{k-1}{k}<d$. This implies Inequality~(\ref{final_ineq}), and hence Inequality~(\ref{main_ineq}) if $n$ is large enough.

	Now let us choose an $N$ so that Inequality~(\ref{main_ineq}) holds for all $n> N$, and then let us choose a $c$ so that $p_k(n)< cn^{dn}$ holds for $n\leq N$. Then we show that $p_k(n) < cn^{dn}$ also holds for $n>N$ by induction on $n$. Suppose that we already know that $p_k(m) < cm^{dm}$ holds for all $m<n$. Then by using the recursion formula (\ref{recursion}) and Inequality~(\ref{main_ineq}) we obtain $$p_k(n)=\sum_{i=0}^{k-1}{{n-1\choose i}p_k(n-1-i)}<c\sum_{i=0}^{k-1}{{n-1\choose i}(n-1-i)^{d(n-1-i)}}<cn^{dn}.$$ 
\end{proof}

%% file: RU2.tex

\begin{section}{Exponential Orbit Growth and Reducts of Unary Structures}

	In this section we show that $\kexp=\reduct(\unary)$.

\begin{lemma}\label{u_sub_kexp}
	$\reduct(\unary)\subseteq \kexp$. 
\end{lemma}

\begin{proof}
	Let $\fa\in \unary$, and let $O_1,\dots,O_k$ be the orbits of $\aut(\fa)$. Then $\aut(\fa)=\prod_{i=1}^k{\sym(O_i)}$ by Lemma \ref{unary}. Hence, the number of injective $n$-orbits of $\aut(\fa)$ is at most $k^n$. This implies that $\fa\in \kexp$, and hence $\unary\subseteq \kexp$. The statement follows from the fact that the class $\kexp$ is closed under taking first-order reducts.
\end{proof}

\begin{lemma}\label{delta_trivi}
	Let $\fa\in \kexp$. Then $\Delta(\fa)$ is trivial on each infinite orbit of $\aut(\fa)$.
\end{lemma}
\begin{proof}
	Let us apply Lemma \ref{small_classes_2} for $G=\aut(\fa)$, $R=\Delta(\fa)$, some $c$ such that $0<c<\frac{1}{2}$, and $k=2$. 
	We obtain that $\Delta(\fa)$ has at most finitely many classes of size at least $2$. 
	If $O$ is an infinite orbit of $\aut(\fa)$, then every class of $\Delta(\fa)$ contained in $O$ has the same size. Therefore, $\Delta(\fa)$ must be trivial on each infinite orbit $O$.  
\end{proof}

\begin{lemma}\label{delta_trivi2}
	Let $\fa\in \kexp$. Then $\fa$ is a first-order reduct of some structure $\fb\in \kexp$ for which $\Delta(\fb)$ is trivial.
\end{lemma}

\begin{proof}
	Let $F$ be the union of finite orbits of $\fa$. Then $F$ is finite. Let $\fb$ be a structure obtained from $\fa$ by adding a constant for each element of $F$. Then $\aut(\fb)=\aut(\fa)|_F$, and it is easy to see that $\aut(\fa)|_F\in \gexp$. Therefore $\fb\in \kexp$. Now let $C$ be a class of $\Delta(\fb)$. If $C$ is contained in a finite orbit of $\fa$, then by definition $|C|=1$, and if $C$ is contained in an infinite orbit of $\fa$, then $|C|=1$ by Lemma \ref{delta_trivi}. Therefore $\Delta(\fb)$ is trivial.
\end{proof}

\begin{theorem}\label{main_kexp}
	$\kexp=\reduct(\unary)$.
\end{theorem}

\begin{proof}
	The containment ``$\supseteq$'' is Lemma \ref{u_sub_kexp}.
	Now assume that $\fa\in \kexp$. By Lemma \ref{delta_trivi2} $\fa$ is a first-order reduct of some $\fb\in \kexp\subseteq \kexpp$ such that $\Delta(\fb)$ is trivial. By Corollary \ref{extension_4} we obtain that $\fb\in \reduct(\unary)$. Hence, $\fa\in (\reduct\circ \reduct)(\unary)=\reduct(\unary)$.
\end{proof}

\end{section}

%% file: additional.tex

\section{Additional Descriptions of the Classes $\kexp$ and $\kexpp$}
\label{sect:additional}
In this section we present additional descriptions of the class $\kexpp$ that follow from our main results. 

\subsection{Generating $\kexpp$ from $\sets$}
	In this subsection we show that $\kexpp$ is the smallest class that contains $\sets$ and is closed under taking first-order reducts, finite covering structures, and adding constants.
	
\begin{lemma}\label{gen_from_set}
The following inclusions hold. 
\begin{enumerate}
\item $\unary^* \subseteq (\reduct\circ \fcover\circ \const)(\sets)$.
\item $\unary_{\nf} \subseteq (\reduct\circ \fcover)(\set)$.
\end{enumerate}
\end{lemma}
\begin{proof}
	Let $\fa\in \unary^*$ and let $O_1,\dots,O_m$ be the orbits of $\aut(\fa)$ so that $O_1=\{y_1\},\dots,O_l=\{y_l\}$ are the finite orbits. Let $F=\{y_1,\dots,y_l\}$. Pick a bijection $b_i$ between $O_i$ and $\mathbb{N}$ for each $i \in \{l+1,\dots,m\}$. Let $b=\bigcup_{i=l+1}^{m}b_i$  and let $E:=\{(x,y) \mid x,y\in \bigcup_{i=l+1}^{m}O_i, b(x)=b(y)\}$. Let $\fc$ be the structure $\fa$ expanded by the relation $E$. Then $\Delta(\fc)=E\cup \{(x,x) \mid x\in F\}$ and $\aut(\fc/\Delta(\fc))=\sym(C/\Delta(\fc))_F$. Therefore, $\fc/\Delta(\fc)\in \const(\sets)$, which shows (1). If $\fa$ has no finite orbits, then $l=0$ and $\fc/\Delta(\fc)$ is bi-definable with $\set$, which shows (2). 
\end{proof}
\begin{lemma}\label{stab_finite}
	The classes $\kexp$ and $\kexpp$ are closed under $\const$.
\end{lemma}
\begin{proof}
	We need to show that if a permutation group $G$ on $X$ is in $\gexp$ or in $\gexpp$, then so is $G_F$ for any finite $F \subset X$. In the case of $\gexp$ this is clear. For the class $\gexpp$ this is stated in Lemma \ref{stabil_finite}.
\end{proof}
\begin{lemma}\label{fu_uf}
	For any class $\mathcal C$ of structures 
	\begin{align}
	(\fcover\circ \reduct) ({\mathcal C}) & \subseteq (\reduct\circ \fcover) ({\mathcal C}). 
	\label{eq:fu_uf} \\
(\fcover\circ \reductfin) ({\mathcal C}) & \subseteq (\reductfin\circ \fcover) ({\mathcal C}). \label{eq:fu_uf2}
\end{align}
\end{lemma}
\begin{proof}
	Let $\fc$ be a structure, let $\fb$ be a first-order reduct of $\fc$, and let $\pi \colon \fa \rightarrow \fb$ be a finite cover. Let $G:=\aut(\fa) \cap \mu_\pi^{-1}(\aut(\fc))$. Then $G$ is closed. So $G$ is the automorphism group of some first-order expansion $\fd$ of $\fa$ and $\pi \colon \fd \rightarrow \fc$ is a finite cover. Hence, $\fa \in (\reduct\circ \fcover)(\fc)$. Moreover, if $\fc$ is $\omega$-categorical and $[\aut(\fb):\aut(\fc)]$ is finite, that is, $\aut(\fb)=g_1\aut(\fc)\cup \dots \cup g_n\aut(\fc)$ for some $g_1,\dots,g_n\in \aut(\fb)$, then $\aut(\fa)=h_1\aut(\fc)\cup \dots \cup h_n\aut(\fc)$ for some $h_i$ so that $\mu_\pi(h_i)=g_i$. In particular, $[\aut(\fa):\aut(\fd)]$ is finite.
\end{proof}

If $\fa$ has no finite orbits, then every first-order reduct of $\fa$ does not have finite orbits, too, so $\reduct({\mathcal C}_{\nf}) \subseteq (\reduct(\mathcal C))_{\nf}$ for any class $\mathcal C$. For $\mathcal C = \unary$ we even get that
\begin{align}
\reduct(\unary_{\nf}) = \reduct(\unary)_{\nf}
\label{eq:r-unf_ru-nf}
\end{align}
(see Corollary~\ref{wreath_infinite}). 
Since a finite covering structure $\fa$ of $\fb$ has finite orbits if and only if $\fb$ has finite orbits
we have for any class ${\mathcal C}$ of structures
that 
\begin{align}
\fcover({\mathcal C}_{\nf}) = (\fcover(\mathcal C))_{\nf}. 
\label{eq:fnf_nff}
\end{align}


\begin{theorem}\label{gen_from_set_final}
The following equalities hold. 
\begin{enumerate}
\item $\kexpp=(\reduct\circ \fcover\circ \const)(\sets)$, 
\item $(\kexpp)_{\nf}=(\reduct \circ \fcover)(\set)$.
\end{enumerate}
\end{theorem}

\begin{proof}
	Clearly, $\set \in \sets \subset \kexp \subset \kexpp$ and
	$\kexp$ and $\kexpp$ are closed under $\reduct$.	
	The closure of $\kexpp$ under $\fcover$ follows from
	$\kexpp = (\fcover \circ R)(\unary^*)$ 
	(Theorem~\ref{main_kexpp}) 
	and the closure under $\const$ is stated in Lemma~\ref{stab_finite}. 
	Finally, $(\kexpp)_{\nf}$ is  closed under 
	$\reduct$ since 
	$\reduct((\kexpp)_{\nf}) \subseteq (\reduct(\kexpp))_{\nf} = (\kexpp)_{\nf}$. 
	This shows the inclusions $\supseteq$
	in $(1)$ and in $(2)$. 
	For the converse containments observe that  
	\begin{align*}
	\kexpp & =(\fcover \circ \reduct)(\unary^*) && \text{(by Theorem~\ref{main_kexpp})} \\
	& \subseteq (\fcover\circ \reduct \circ \reduct \circ \fcover \circ \const)(\sets) && \text{(by Lemma~\ref{gen_from_set} (1))} \\
	& = (\fcover \circ \reduct \circ \fcover \circ  \const)(\sets)  \\
		& \subseteq (\reduct\circ \fcover \circ \fcover \circ \const)(\sets) && \text{(by~(\ref{eq:fu_uf}) in Lemma~\ref{fu_uf})} \\
	& = (\reduct\circ \fcover \circ \const)(\sets) 
	 \end{align*} 
	 which shows $(1)$. 
	 Moreover, 
	 \begin{align*}
	(\kexpp)_{\nf} & =  ((\fcover \circ \reduct)(\unary^*))_{\nf} && \text{(by Theorem~\ref{main_kexpp})} \\
	& = \fcover((\reduct(\unary))_{\nf}) && \text{(by~(\ref{eq:fnf_nff})} \\
	& = (\fcover \circ \reduct) (\unary_{\nf}) && \text{(\ref{eq:r-unf_ru-nf})} \\
	& \subseteq (\fcover \circ \reduct \circ \reduct \circ \fcover)(\set) && \text{(by Lemma~\ref{gen_from_set} (2))} \\
	& \subseteq (\reduct \circ \fcover)(\set) && \text{(as above)}
	 \end{align*}
	 which shows $(2)$. 
\end{proof}

\subsection{Model-complete cores}
	The model-complete core of an $\omega$-categorical structure has already been defined in the introduction. 
	In this section we show that 
	$\kexpp$ is the smallest
	class of structures that contains $\set$ and is closed under taking first-order reducts, finite covers, and model-complete cores.

\begin{lemma}\label{mc_core_orbit}
	Let $\fa$ be an $\omega$-categorical structure and $\fb$ its model-complete core. 
	Then $\orb_n(\fb)\leq \orb_n(\fa)$ and $\oi_n(\fb)\leq \oi_n(\fa)$ for all $n \in {\mathbb N}$.
\end{lemma}

\begin{proof}
	For $\orb_n$ this is Proposition~3.6.24.\ in~\cite{Bodirsky-HDR}. The statement for $\oi_n$ can be shown analogously. 
\end{proof}

\begin{corollary}\label{kexpp_closed_under_m}
	The classes $\kexp$ and $\kexpp$ are closed under $M$. 
\end{corollary}

\begin{remark}
Analogous statements hold for the \emph{model companion} instead of the model-complete core.
\end{remark}

\begin{definition}
	Let $\fa$ be a structure with signature $\tau$ and let $F \subseteq A$. Then let $\fa(F)$ denote the following $\tau$-structure.
\begin{itemize}
\item The domain of $\fa(F)$ is $A(F) := (F \times \mathbb{N}) \sqcup ((A \setminus F) \times \{0\})$.
\item For each $R \in \tau$ of arity $k$ the relation $R^{\fa(F)}$ is defined as $\{((x_1,n_1),\dots,(x_k,n_k)) \mid (x_1,\dots,x_n) \in R\}$.
\end{itemize}
\end{definition}

\begin{remark}
	The map $f \colon A(F) \to A$ defined by $(x,n) \mapsto x$ is a homomorphism from $\fa(F)$ to $\fa$. 
	Conversely, the mapping $g \colon A \to A(F)$ defined by $x \mapsto (x,0)$ is a homomorphism from $\fa$ to $\fa(F)$ (in fact it is an embedding). Therefore, $\fa$ and $\fa(F)$ are homomorphically equivalent.
\end{remark}

\begin{remark}\label{reduct_split}
	It follows directly from the definition that if $\fa$ is a first-order reduct of $\fb$, and $F \subseteq A$, then $\fa(F)$ is a first-order reduct of $\fb(F)$ (since we can use the same definitions).
\end{remark}

\begin{lemma}\label{get_rid_of_finite}
	Let $\fa\in \unary^*$ and let $F$ be the union of the finite orbits of $\fa$. Then $\fa(F) \in \unary_{\nf}$.
\end{lemma}
\begin{proof}
	Let $O_1,\dots,O_k$ be the orbits of $\aut(\fa)$. Then $\aut(\fa)=\prod_{i=1}^k {\sym(O_i)}$ by Lemma \ref{unary}. Let $O_1=\{y_1\},\dots,O_l=\{y_l\}$ be the finite orbits of $\aut(\fa)$. Then $\aut(\fa(F))=\prod_{i=1}^l{\sym(\{y_i\}\times 
	{\mathbb N})} \times \prod_{i=l+1}^k{\sym(O_i \times \{0\})}$ has no finite orbits and therefore $\fa(F)\in \unary_{\nf}$. 
\end{proof}

\begin{lemma}\label{get_rid_of_finite2}
	Let $\fa\in \reduct(\unary)$ and let $F$ be the union of the finite orbits of $\fa$. Then $\fa(F)\in \reduct(\unary)_{\nf}$.
\end{lemma}
\begin{proof}
	Let $C_1,\dots,C_n$ be the classes of $\nabla(\fa)$. Then $\prod_{i=1}^n{\sym(C_i)} \subseteq \aut(\fa)$ by Lemma \ref{unary_reduct}, and 
	hence, $\fa$ is a first-order reduct of $\fb \in \unary^*$. Lemma~\ref{get_rid_of_finite} implies
	that $\fb(F) \in \unary_{\nf}$. 
	Remark \ref{reduct_split} implies 
	that $\fa(F)$ is a first-order reduct of $\fb(F)$. 
	Hence, 
	$\fa(F) \in \reduct(\unary_{\nf}) = \reduct({\unary})_{\nf}$.
\end{proof}

\begin{corollary}\label{ru_model}
	Every structure $\fa\in \reduct(\unary)$ is interdefinable with a model-complete core of a structure in $\reduct(\unary)_{\nf}$, i.e.,
	$$\reduct(\unary) \subseteq M(\reduct(\unary)_{\nf}).$$ 
\end{corollary}

\begin{proof}
	Let $\fa^*$ be the expansion of $\fa$ by all first-order definable relations. Then $\fa^*$ is a  model-complete core and interdefinable with $\fa$. Let $F$ be the union of the finite orbits of $\aut(\fa^*) = \aut(\fa)$. By Lemma \ref{get_rid_of_finite2} we know that $\fa^*(F)\in \reduct(\unary)_{\nf}$. Since $\fa^*$ and $\fa^*(F)$ are homomorphically equivalent it follows that $\fa^*$ is the model-complete core of $\fa^*(F)$.
	Hence, $\fa \in M(\reduct(\unary)_{\nf})$. 
\end{proof}

\begin{corollary}\label{fru_model}
	Let $\fb\in \reduct(\unary^*)$ and let $\pi \colon\fa\rightarrow \fb$ be a finite cover. Then $\fa$ is interdefinable with a model-complete core of a structure in $\fcover(\reduct(\unary^*))_{\nf}$, i.e.,
	$$F(\reduct(\unary^*)) \subseteq M(\fcover(\reduct(\unary^*))_{\nf}).$$ 
\end{corollary}

\begin{proof}
	As in the previous proof let $\fa^*$ be the expansion of $\fa$ by all relations that are first-order definable in $\fa$, and let $F$ 
	be the union of the finite orbits of $\aut(\fa) = \aut(\fa^*)$. Then $\pi(F)$ is the union of finite orbits of $\aut(\fb)$. By Corollary \ref{get_rid_of_finite2} we know that $\fb(\pi(F))\in \reduct(\unary)$. Let $\pi' \colon A(F) \rightarrow B(\pi(F))$ be defined as 

	\begin{align*}
	\pi'(x,n) := \begin{cases}
	(\pi(x),n) & \text{if } x\in F \\
	(\pi(x),0) & \text{otherwise.}
	\end{cases}
	\end{align*}
	
	Then it is easy to see that $\pi' \colon \fa(F) \rightarrow \fb(\pi(F))$ is a finite covering map. Hence, $\fa(F)\in \fcover(\reduct(\unary^*))$. By 
	Lemma~\ref{get_rid_of_finite} the structure
	$\fa(F)$ has no finite orbits, and as before we can conclude that $\fa$ is the model-complete core of $\fa(F)$.
\end{proof}

\begin{lemma}\label{gen_mccore}
The following identities hold. 
\begin{enumerate}
\item $\kexp=\mcore((\kexp)_{\nf})$,
\item $\kexpp=\mcore((\kexpp)_{\nf})$,
\item $\kexpp=(\mcore\circ \reduct\circ \fcover)(\set)$.
\end{enumerate}
\end{lemma}

\begin{proof}
	The containments ``$\supseteq$'' in item (1) and item (2) follow from Corollary \ref{kexpp_closed_under_m}. By Theorems \ref{main_kexp} and \ref{main_kexpp} we know that $\kexp=\reduct(\unary)$ and $\kexpp=\fcover(\reduct(\unary))$. Then the containments ``$\subseteq$'' in items (1) and item (2) follow from Corollaries \ref{ru_model} and \ref{fru_model}.
	To show Item (3), observe that
	\begin{align*}
	\kexpp & = \mcore((\kexpp)_{\nf}) && \text{(by item (2) of the lemma)} \\
	& = \mcore(\reduct(\fcover(\set))) && \text{(by item (2) of 
	Theorem~\ref{gen_from_set_final})}
	\end{align*}
\end{proof}


\subsection{Summary}
\label{sect:summary}
	The following theorem summarizes some of the equivalent characterizations of the classes $\kexp,\kexpp,(\kexp)_{\nf},(\kexpp)_{\nf}$.

\begin{theorem}\label{summary}
\begin{align}\kexp & =\reduct(\unary)=\reductfin(\unary^*) \label{eq:kexp} \\
(\kexp)_{\nf}& =\reduct(\unary_{\nf})=\reductfin(\unary_{\nf}) \label{eq:kexpnf} \\
\kexpp & =(\fcover\circ \reduct) (\unary)=(\reduct\circ \fcover)(\unary) \label{eq:kexpp} \\
& =(\fcover\circ \reductfin)(\unary^*)=(\reductfin\circ \fcover)(\unary^*) \nonumber \\ 
& =(\reduct\circ \fcover\circ \const)(\sets)=(\mcore\circ \reduct\circ \fcover)(\set) \nonumber  \\
(\kexpp)_{\nf}
& =(\reduct\circ \fcover)(\set)
=(\fcover\circ \reduct)(\unary_{\nf})  \label{eq:kexppnf}  \\
& =(\fcover\circ \reductfin)(\unary_{\nf})
=(\reductfin\circ \fcover)(\unary_{\nf}) \nonumber
\end{align}
\end{theorem}

\begin{proof}
	(\ref{eq:kexp}): Corollary~\ref{ru_rfiniteu} states that $R(\unary) = \reductfin(\unary^*)$ 
	and Theorem \ref{main_kexp} that $\kexp=\reduct(\unary)$.

	(\ref{eq:kexpnf}): 
	we have $\reductfin(\unary_{\nf})  \subseteq \reduct(\unary_{\nf}) = \reduct(\unary)_{\nf} = (\kexp)_{\nf}$ by (\ref{eq:r-unf_ru-nf}) and
(\ref{eq:kexp}), and $\reduct(\unary)_{\nf} \subseteq \reductfin(\unary_{\nf})$ can be shown as in the proof of Corollary~\ref{ru_rfiniteu}.
	


	(\ref{eq:kexpp}): By Theorem~\ref{main_kexpp} we know that $\kexpp=(\fcover\circ \reduct)(\unary)=(\reductfin\circ \fcover)(\unary^*)$. 
	This also implies that the class $\kexpp$ is closed under $\fcover$, and it is obviously closed under $\reduct$, so $$\kexpp = (\reductfin\circ \fcover)(\unary^*) \subseteq  (\reduct \circ \fcover)(\unary) \subseteq \kexpp.$$ 
	The equality $(\fcover\circ \reduct)(\unary)=(\fcover\circ \reductfin)(\unary^*)$ follows from the fact that $\reduct(\unary)=\reductfin(\unary^*)$ (Corollary~\ref{ru_rfiniteu}). 
	The equality 
	$\kexpp = (R \circ F \circ C)(\sets)$ is item (1) of Theorem~\ref{gen_from_set_final} 
	and the equality 
	$\kexpp = (M \circ R \circ F)(\set)$
	is item (3) of Lemma~\ref{gen_mccore}.
		
	(\ref{eq:kexppnf}): 
	the proof of Theorem~\ref{gen_from_set_final} (2)
	shows the following equalities: $(\kexpp)_{\nf}=(\fcover\circ \reduct)(\unary_{\nf}) = (\reduct\circ \fcover)(\set)$. 
	Finally, 
	\begin{align*}
	(\kexpp)_{\nf}  = (\fcover\circ \reduct)(\unary_{\nf})
		&  = 
	(\fcover\circ \reductfin)(\unary_{\nf}) && \text{(as in Corollary~\ref{ru_rfiniteu})} \\
	& \subseteq (\reductfin\circ \fcover)(\unary_{\nf}) && \text{(by (\ref{eq:fu_uf}))} \\
	& \subseteq (\kexpp)_{\nf}
	\end{align*} 
	and thus $(\kexpp)_{\nf}=(\fcover\circ \reductfin)(\unary_{\nf}) =(\reductfin\circ \fcover)(\unary_{\nf})$. 
\end{proof}

\section{Consequences for Constraint Satisfaction} 
In the introduction we have already mentioned 
that for finite structures $\fa$ there is a complexity dichotomy for $\csp(\fa)$: these problems are in P or NP-complete. Such a complexity dichotomy has also been conjectured for the much larger class of
 first-order reducts of \emph{finitely bounded} homogeneous structures.  A structure $\fb$ with finite relational signature $\tau$ is called \emph{finitely bounded} if there exists a finite set of finite $\tau$-structures $\mathcal F$ such that a finite $\tau$-structure $\fa$ embeds into $\fb$ if and only if no structure from $\mathcal F$ embeds into $\fa$. For first-order reducts of finitely bounded homogeneous structures there is also a more specific 
 \emph{infinite-domain tractability conjecture}~\cite{BPP-projective-homomorphisms}: assuming that $\fa$ is a model-complete core the conjecture says that $\csp(\fa)$ is in P if and only if $\fa$ has a pseudo-Siggers polymorphism (for a definition of pseudo-Siggers polymorphisms and a proof that the conjecture can be phrased like this, see~\cite{BartoPinskerDichotomy}).
 
Let $\fa$ be a structure from $\kexpp$ with finite relational signature. The next lemma shows that the question whether $\csp(\fa)$ is in P or NP-complete falls into the scope of this conjecture. 

\begin{lemma}\label{lem:finitely-bounded}
Every structure in $\kexpp$ is a first-order reduct of a finitely bounded homogeneous structure. 
\end{lemma}
\begin{proof}
Let $\fa \in \kexpp$. 
By Theorem~\ref{main_kexpp} we have 
$\kexpp=\reduct(\fcover(\unary^*))$, so
$\fa$ is a first-order reduct of a structure $\fa' \in \fcover(\unary^*)$. By Proposition~\ref{reduct_trivial2}, every finite cover
of a structure in $\unary^*$ is strongly split, so we can assume that $\fa'$ is a strongly trivial covering structure of a structure $\fb \in \unary^*$. 
Let $\fc$ be the structure from the proof of Lemma~\ref{ramsey}, and let $\tau := (\{U_{i,s} \mid i \leq k, s \in F_i\} \cup \{\sim_\pi\})$ be the signature of $\fc$. 
Then it is easy to specify a finite set of forbidden finite $\tau$-structures such that in any finite $\tau$-structure that avoids these structures 
\begin{itemize}
\item the 
relation $\sim_\pi$ is an equivalence relation,
\item 
the sets denoted by the unary relations $U_{i,s}$
are pairwise disjoint and cover all of $C$, 
\item for all $i,s$ if $x \sim_\pi y$ and $x,y \in U_{i,s}$ then $x=y$, 
\item for all $i,s$ the cardinality of $U_{i,j}$ is at most the cardinality of $U_{i,j}$ in $\fc$. 
\end{itemize}
These are precisely the finite structures that embed into $\fc$. 
\end{proof}
Let $\fa$ be a structure from $\reduct(\fcover(\unary))$. 
In this section we discuss the consequences of our results 
for classifying the computational complexity of $\csp(\fa)$. First, since $\reduct(\fcover(\unary)) = \kexpp$ is closed under $\mcore$ as discussed above we can assume that $\fa$ is a model-complete core. The following lemma shows that we can even assume that $\fa \in \fcover(\unary)$. 
\begin{lemma}
Let $\fa \in \reduct(\fcover(\unary))$. Then there exists a model-complete core $\fc$ in $\fcover(\unary^*)$ such that 
\begin{itemize}
\item $\csp(\fa)$ and $\csp(\fc)$ are polynomial-time equivalent;
\item the $\nabla(\fc)$-classes are the orbits of $\aut(\fc)$ and they are primitively positively definable in $\fc$. 
\end{itemize}
\end{lemma}
\begin{proof}
Let $\fc'$ be the model-complete core of $\fa$.
Then $\fc'$ is in 
$M(\reduct(\fcover(\unary))) = \kexpp = \fcover(\reduct(\unary))$. 
So suppose that $\pi \colon \fc' \to \fb'$ is a finite covering for $\fb' \in \reduct(\unary)$. 
Add a constant $c$ from each $\nabla(\fc')$-equivalence class to $\fc'$ and let $\fc$ be the resulting structure. 
Then $\fc$ is still $\omega$-categorical and a model-complete core. 
Moreover, $\fc$ and $\fa$ are polynomial-time equivalent~\cite{Cores-Journal}.  

Add a constant $\pi(c)$ to $\fb'$
for each of the new constants $c$ and let 
 $\fb$ be the structure obtained in this way. 
  The proof of Corollary~\ref{add_constants_ru} shows that $\fb \in \unary^*$. Then $\pi \colon \fc \to \fb$ is a finite cover. Therefore, $\fc \in \fcover(\unary^*)$.  
Moreover, the $\nabla(\fc)$-classes are the orbits of $\fc$ and orbits in model-complete cores are primitive positive definable~\cite{Cores-Journal}. 
\end{proof}
It can be shown using the universal-algebraic approach to constraint satisfaction that  
if $\fc \in \fcover(\reduct(\unary))$ 
then $\csp(\fc)$ is either in P or NP-complete.
This lies beyond the scope of this article, but will appear elsewhere.
